\documentclass[letterpaper,11pt,twoside]{article}
\usepackage{amsmath,amsfonts,amssymb,amsthm,mathabx,etoolbox}
\usepackage[utf8]{inputenc}
\usepackage[margin=2cm]{geometry}
\usepackage[scr=boondoxo]{mathalpha}
\usepackage{booktabs}

\usepackage{xcolor}

\usepackage{secdot}
\usepackage{graphicx}
\graphicspath{{images/}}

\usepackage{caption}
\usepackage{tikz-cd}
\usepackage{titlecaps}
\Addlcwords{a is with of for on in are and its to the by under around}

\newcommand{\mytitle}{ 
Graph Parameterization for Invariant Manifolds Attached to Quasi-Periodic Tori
}

\newcommand{\myshorttitle}{
Graph Parameterization of Reducible Invariant Tori
}

\usepackage{hyperref}
\hypersetup{
hidelinks=true,
bookmarksnumbered=true,
bookmarksopen=true,
pdftitle={\myshorttitle{}},
pdfauthor={M. Barcelona, J. Gimeno, and M. Jorba-Cusco},
}
\usepackage{orcidlink}

\usepackage{comment} 

\newtheorem{thm}{Theorem}[section]
\newtheorem{lem}[thm]{Lemma}
\newtheorem{cor}[thm]{Corollary}
\newtheorem{pro}[thm]{Proposition}

\theoremstyle{remark}
\newtheorem{rmk}[thm]{Remark}
\AtBeginEnvironment{rmk}{%
  \pushQED{\qed}%
}
\AtEndEnvironment{rmk}{%
  \popQED
}
\theoremstyle{definition}
\newtheorem{dfn}[thm]{Definition}
\newtheorem{ex}[thm]{Example}

\newtheorem{palg}[thm]{Pseudo-Algorithm}
\newtheorem{test}[thm]{Test}

\newcommand{\R}{\mathbb{R}}
\newcommand{\C}{\mathbb{C}}
\newcommand{\N}{\mathbb{N}}
\newcommand{\Q}{\mathbb{Q}}
\newcommand{\Z}{\mathbb{Z}}
\newcommand{\T}{\mathbb{T}}

\newcommand{\flux}{\Phi}
\newcommand{\Nseg}{N_{\mathrm{seg}}}

\newcommand{\nt}{n _L}
\newcommand{\nh}{n _N}

\newcommand{\st}{s}
\newcommand{\sh}{r}

\newcommand{\Et}{E _L}
\newcommand{\Eh}{E _N}

\newcommand{\Vt}{V_{L}}
\newcommand{\Vh}{V_{N}}

\newcommand{\Lt}{\Lambda _{L}}
\newcommand{\Lh}{\Lambda _{N}}

\newcommand{\tor}{\varphi}
\newcommand{\etor}{\mathcal E _{\textrm{tor}}}
\newcommand{\ered}{\mathcal E _{\textrm{red}}}

\newcommand{\graphfun}{W}
\newcommand{\poincmap}{P}
\newcommand{\expoincmap}{\mathcal{P}}
\newcommand{\remainL}{F}
\newcommand{\remainN}{G}
\newcommand{\param}{K}
\newcommand{\torinv}{\mathcal T}
\newcommand{\mani}{\mathcal K}
\newcommand{\adframe}{Q}
\newcommand{\intdyn}{f}

\newcommand{\cR}{\Upsilon}

\newcommand{\I}{\mathbf{i}}

\newcommand{\dop}{\mathrm{D}} 
\newcommand{\idop}{id}        
\newcommand{\idmat}{Id}       
\newcommand{\kop}{\mathfrak{K}}          

\DeclareMathOperator{\diag}{diag}
\DeclareMathOperator{\Spec}{Spec}

\DeclareMathOperator{\Range}{Range}
\DeclareMathOperator{\vecop}{vec}

\newcommand{\innerp}[2]{\langle {#1}, {#2} \rangle}
\newcommand{\tol}{\mathtt{tol}}

\newcommand{\bydef}{\,\stackrel{\mbox{\tiny\textnormal{\raisebox{0ex}[0ex][0ex]{def}}}}{=}\,}

\title{\mytitle{}}

\author{
  Miquel Barcelona\,\orcidlink{0000-0002-0036-5240}\textsuperscript{(1)} \and 
  Joan Gimeno\,\orcidlink{0000-0002-8707-6379}\textsuperscript{(1)} \and 
  Marc Jorba-Cusc\'o\,\orcidlink{0000-0003-0308-3756}\textsuperscript{(2,3,\textasteriskcentered)}
}

\begin{document}

\maketitle
{\small
\begin{itemize}
\renewcommand{\itemsep}{0pt}
\item [(1)] Departament de Matem\`atiques i Inform\`atica, %
Universitat de Barcelona, %
Gran Via de les Corts Catalanes, 585, 08007 Barcelona, %
Spain, \texttt{miquel.barcelona@ub.edu} and  {\texttt{joan@maia.ub.es}}
\item[(2)] Departament de Matem\`atiques, Universitat Polit\`ecnica de Catalunya, 
Av. Diagonal 647, 08028, Barcelona, Spain, 
\verb+marc.jorba@upc.edu+ 
\item[(3)] Centre de Recerca Matem\`atica (CRM), 
Edifici C, Campus UAB, 0 Floor, 08193 Bellaterra, Barcelona, Spain
\item[\textsuperscript\textasteriskcentered] {\footnotesize
  Corresponding Author}
\end{itemize}
}

\begin{abstract}
We develop a high-order Fourier–Taylor graph-parameterization method for computing local invariant manifolds attached to reducible quasi-periodic invariant tori of discrete dynamical systems, with particular emphasis on Poincar\'e maps. Assuming that the invariant torus and its Floquet reduction are known, we represent the manifold as a graph over a selected spectral subspace and compute its coefficients recursively through jet transport and cohomological equations on the torus. The formulation accommodates stable, unstable, and center spectral selections. For center manifolds, it simultaneously determines the nonlinear reduced dynamics, providing a low-dimensional representation from which nearby invariant objects can be lifted to the original phase space.

\smallskip

For real systems with complex-conjugate Floquet multipliers, we derive an equivalent real formulation based on Hermitian symmetry and Krawtchouk transformations. This allows real-valued jet transport while retaining the advantages of complex spectral coordinates in the cohomological solver. The implementation also incorporates inverse dynamics for stable directions, multiple shooting, coordinate scaling, and a posteriori validation tests.

\smallskip

Numerical experiments with a quasi-periodically forced pendulum, an oscillating solar sail in the Earth–Moon system, and a simplified solar system model illustrate the construction and its computational properties. In the solar-sail benchmark, the real formulation reduces computation time by approximately a factor of 2.2 compared with the complex implementation.
\end{abstract}

\medskip


\pagestyle{myheadings}
\markboth{\myshorttitle{}}
{M. Barcelona, J. Gimeno, and M. Jorba-Cusc\'o}
\newpage
\phantomsection\pdfbookmark[1]{\contentsname}{toc}
{\small \tableofcontents}
\newpage

\section{Introduction} \label{sec.intro}

Local invariant manifolds are a standard way to organize the dynamics near an
invariant object. Besides separating stable, unstable, and center behavior, a
high-order parameterization provides coordinates in which the dynamics can be
studied at much lower dimension and then lifted back to the original phase
space. A typical approach is to express some of these coordinates as 
a function of the others, namely, the graph-parameterization method. 
The original idea can be traced back to Hadamard's seminal paper \cite{Hadamard1901}.
See also \cite{HirschPS77} for comprehensive text on classical view 
of the transform method, mainly devoted to theoretical issues. 

The advent of computers and development of numerical methods extended the 
applicability of the method also to the semi-analytical computation 
of high-order invariant manifolds. For equilibria, Farr\'es and Jorba developed a high-order center-manifold
algorithm based explicitly on a graph-parameterization and recursive expansions
\cite{FarresJ10b}. Related parameterization methods for invariant manifolds are
developed in \cite{CabreFLl05,HaroCLMF16}, while Floquet coordinates lead to
closely related Fourier--Taylor constructions around periodic orbits
\cite{CastelliLMJ15}. Recent work has also extended these ideas to center manifolds,
connections, and invariant tori arising in celestial mechanics
\cite{BarcelonaHM2024,FernandezHM2024,FernandezHLM2026}.

In the particular case of quasi-periodic invariant tori, an additional
difficulty arises: the linearized normal dynamics depends on the phase on the
torus. When this quasi-periodic linear cocycle is reducible, a phase-dependent
Floquet transformation turns it into a constant matrix. The normal spectral
splitting can then be treated much as in the case of a fixed point, while the
phase dependence is represented by Fourier series. Numerical algorithms that
compute a torus together with its Floquet transformation and Floquet matrix are
described in \cite{JorbaO09,CanadellH2014,CanadellH2017}; recent work also uses
parameterization methods to continue reducible tori through saddle-node
bifurcations \cite{FiguerasGP2026}. The computation of whiskered tori, their
invariant bundles, and their associated manifolds have been developed in several
parameterization and flow-map formulations; see, for example,
\cite{HuguetLlS2012,HuguetLlS2013,KumarAL2022,KumarAL2025,FernandezHM2024,
FernandezHLM2026}. The closest numerical predecessor to the construction of this
manuscript is \cite{GimenoJNO21}, where high-order Fourier--Taylor expansions of
stable and unstable manifolds of reducible tori are computed in Poincar\'e maps.

In this work, we formulate the local manifold explicitly as a graph over a
chosen spectral subspace of the constant Floquet matrix. This point of view
allows the same construction to be used for stable, unstable, and center
selections: the chosen spectral variables are kept as independent coordinates
and the remaining normal directions are solved as a graph. Expanding that graph
in Fourier--Taylor series reduces the computation, order by order, to jet
transport followed by linear cohomological equations on the torus. Unstable
manifolds are treated with the forward map, and stable manifolds are treated
with the inverse
map. We also give a parallel-shooting version of the construction, in the spirit
of recent multiple-shooting parameterization methods
\cite{gonzalez2017,Kumar2026}.

Complex-conjugate Floquet multipliers naturally occur even when the original
system is real. A novelty of this work is the following: 
rather than evaluating the original dynamics at complex states,
we exploit the Hermitian symmetry of Fourier--Taylor coefficients and derive an
equivalent real formulation using Krawtchouk transformations. When the selected
spectral subspace contains the center directions, the parameterization also
produces an internal map on the manifold. This gives a reduced model together
with an explicit lifting to the full phase space, connecting the geometric
construction with invariant-manifold-based model reduction
\cite{JainHaller2022,LiHaller2022}.

We assume throughout this work that the invariant torus and a smooth Floquet
reduction have already been computed to the requested accuracy. The scalar
implementation used in the numerical examples also assumes that the Floquet
matrix is diagonalizable over $\C$. Repeated multipliers are compatible with
this assumption, but a nontrivial Jordan structure leads naturally to block
Sylvester equations rather than the scalar divisions used by the diagonal
solver. Reducibility itself may also fail along families of invariant tori; this
can be detected numerically and can have genuine dynamical consequences
\cite{Jorba01,VitoloBS2011}. These nonreducible and genuinely Jordan-block cases
are not part of the numerical study of the current work. A rather general
classification of cocycles into reducible, non-reducible and almost reducible
can be found in the works by Avila and Krikorian\cite{AvilaK2006, AvilaK2015}.

The method is illustrated on several models: a quasi-periodically forced
pendulum, on an oscillating solar-sail model in the Moon--Earth system and a
simplified solar system model. The examples are used to test the Fourier--Taylor
graph construction, the real treatment of complex conjugate multipliers, the
reduced dynamics obtained from a center spectral selection and the
multiple-shooting strategy.

\subsection{Organization of the Paper}

Section~\ref{sec.invariant-tori-normal} introduces quasi-periodic skew-product
maps, their construction from angular Poincar\'e maps, and the normal dynamics
and reducibility of invariant tori.

Section~\ref{sec.graph} formulates the graph-parameterization in Floquet-adapted
coordinates and derives the corresponding invariance equation. 

Section~\ref{sec.GT-construction} develops the recursive Fourier--Taylor
construction and presents the computational algorithm, while
Section~\ref{sec.solve.cohom} analyzed the solvability of the resulting
cohomological equations and the associated nonresonance conditions.

Section~\ref{sec.hermitian} established the Hermitian structure of
complex-conjugate Floquet coordinates and develops the equivalent real formulation
using Krawtchouk transformations. 

Section~\ref{sec.center-reduction} describes the reduction to center dynamics
and the construction and lifting of higher-dimensional invariant tori.
Section~\ref{sec.numtools} discusses the numerical implementation, coordinate
scaling, and validation tests, and Section~\ref{sec.examples} presents the
numerical examples and celestial-mechanics models. 

Finally, Appendix~\ref{app.multiple-shooting} develops the multiple-shooting
formulation, while Appendix \ref{app.stable} discusses the backward skew-product
used to compute stable manifolds.

\section{Invariant Tori and Their Normal Dynamics}
\label{sec.invariant-tori-normal}

The starting point for the graph-parameterization is an invariant
torus together with a reducible description of its transverse linear dynamics.
This
section fixes the notation for both. We first describe the skew-product and its
invariant torus, including the Poincar\'e map construction used in the examples,
and then introduce the normal cocycle and its Floquet reduction. The last part
records the small divisors that appear when the torus and the Floquet change of
variables are computed in Fourier series. 

\subsection{Quasi-periodic skew-products and invariant tori}

Let $\mathcal{U}\subset\R^n$ be open and let $\T^d=\R^d/\Z^d$. We consider a
smooth map $\poincmap\colon\mathcal{U}\times\T^d\to\R^n$ and the associated
skew-product
\begin{equation}
    \label{eq.SP}
\begin{cases}
\bar{x}=\poincmap(x,\theta),\\
\bar{\theta}=\theta+\rho \pmod{1},
\end{cases}
\qquad
\expoincmap(x,\theta;\rho)\bydef
\bigl(\poincmap(x,\theta),\theta+\rho\bigr).
\end{equation}
Here, $x\in\mathcal U$ is the fiber variable and $\theta\in\T^d$ is the phase
of the quasi-periodic forcing. The phase evolves independently of $x$ through
the rigid rotation $\theta\mapsto\theta+\rho$.

The torus translation $R_\rho(\theta)=\theta+\rho$ is ergodic with respect
to normalized Haar (Lebesgue) measure if, and only if, 
\[
    \innerp{k}{\rho}\notin\Z,
    \qquad 
    \text{for all } k\in\Z^d\setminus\{0\}.
\]
Equivalently, $1,\rho_1,\ldots,\rho_d$ are linearly independent over $\Q$.
We will refer to a vector satisfying this condition as \emph{totally
irrational modulo one}, or simply totally irrational when there is no risk of
confusion. 
Ergodicity alone, however, gives no
quantitative control on the size of the Fourier divisors; the classical cohomological
equations used require Diophantine estimates when small divisors are present.

An invariant torus of \eqref{eq.SP} is represented by a smooth graph
$\tor\colon\T^d\to\mathcal U$ satisfying
\begin{equation}
    \label{eq.invariance_torus}
    \tor(\theta+\rho)=\poincmap\bigl(\tor(\theta),\theta\bigr).
\end{equation}
The corresponding invariant set in the extended phase space is
$
    \torinv\bydef
    \bigl\{(\tor(\theta),\theta):\theta\in\T^d\bigr\}
$. 
The dynamics on $\torinv$ is therefore conjugate to the rigid rotation by
$\rho$. Notice that the graph map
$\theta\mapsto(\tor(\theta),\theta)$ is automatically an embedding, even if
$\tor$ is not an embedding into $\R^n$ by itself.

The existence and persistence of quasi-periodic invariant tori can be
established in several ways. In Hamiltonian and near-integrable problems this
is the domain of KAM theory; see the foundational works
\cite{Kolmogorov54,Arnold63a,Moser62} and parameterization-based formulations
such as \cite{LlaveGJV05,HaroCLMF16}. Related results are available for other
classes, including conformally symplectic systems
\cite{CallejaCL11,CallejaCGL22b}. Normal-hyperbolicity and invariant-graph
arguments give another route in dissipative and forced systems; see, for
example, \cite{Fenichel1971,Stark97a}. We will not study the existence problem. The torus and its rotation vector are assumed to have already been
computed.

\subsubsection{Quasi-periodic skew-products from Poincar\'e maps}

The skew-product \eqref{eq.SP} also arises naturally from a
quasi-periodically forced differential equation. Consider
\begin{equation}
    \label{eq.nonaut}
    \dot x=X(x,t),
\end{equation}
where $X:\mathcal U \times \R \rightarrow \R^n$ is quasi-periodic in time. We
write
\[
    X(x,t)=\widetilde X(x,\tilde\omega t),
\]
with $\widetilde X\colon\mathcal U\times\T^{1+d}\to\R^n$ periodic with period
one in each angular variable and
$\tilde\omega=(\tilde\omega_0,\ldots,\tilde\omega_d)\in\R^{1+d}$ with
rationally independent components. Introducing
$\tilde\theta=\tilde\omega t$ gives the autonomous system
\begin{equation}
    \label{eq.qpvf}
\begin{cases}
    \dot x=\widetilde X(x,\tilde\theta),\\
    \dot{\tilde\theta}=\tilde\omega.
\end{cases}
\end{equation}
This is the elementary autonomous extension of a quasi-periodic forcing: the
forcing already lives in the compact torus $\T^{1+d}$, so no further
compactification is needed. More general recurrent non-autonomous equations are
often organized as skew-product flows or semi-flows over a compact base (for
example, the hull of an almost-periodic forcing); see
\cite{NovoOS2007} for a general treatment of skew-product semi-flows over
minimal bases. Here we only need the finite-dimensional torus translation in
\eqref{eq.qpvf}.

One angular variable can then be removed when the forcing phases themselves are
used to define the section; see \cite[Section 4]{GimenoJJZ25}.
Choose a component with $\tilde\omega_0\neq0$. After reversing that angular
coordinate if necessary, we may assume $\tilde\omega_0>0$. Write
$\tilde\theta=(\tilde\theta_0,\theta)$, where
$\theta=(\tilde\theta_1,\ldots,\tilde\theta_d)$, and take the angular section
\[
    \Sigma=\bigl\{(x,\tilde\theta)\in\mathcal U\times\T^{1+d}:
    \tilde\theta_0=0\bigr\}.
\]
Since $\dot{\tilde\theta}_0=\tilde\omega_0$ is constant, every return to
$\Sigma$ takes the same time
$
    T=\frac{1}{\tilde\omega_0}
$. 
Let $\flux(t;x,\tilde\theta)$ denote the flow of \eqref{eq.qpvf}, with values
in $\mathcal U\times\T^{1+d}$, and let $\pi_x$ be the projection onto the
$x$-variables. The fiber part of the return map is then
\[
    \poincmap(x,\theta)
    \bydef
    \pi_x\!\left(\flux\bigl(T;x,(0,\theta)\bigr)\right).
\]
During this return, the remaining phases advance by
\[
    \rho_j=\frac{\tilde\omega_j}{\tilde\omega_0},
    \qquad j=1,\ldots,d,
\]
modulo integers. Hence, the full return map on $\Sigma$ has exactly the form
\eqref{eq.SP}, with rotation vector $\rho=(\rho_1,\ldots,\rho_d)$. Because the
components of $\tilde\omega$ are rationally independent, this $\rho$ is totally
irrational and, therefore, these two vectors describe the same forcing after one
angular phase has been chosen as the time map.

A $d$-dimensional invariant graph satisfying \eqref{eq.invariance_torus}
suspends under the flow to a $(1+d)$-dimensional invariant torus of
\eqref{eq.qpvf}, provided the flow torus is transverse to the angular section.
Conversely, intersecting such a flow torus with $\Sigma$ removes the phase used
as the flight time and produces a $d$-torus for the return map. When $d=0$, the
forcing is periodic and the construction reduces to the familiar stroboscopic
map; the invariant torus of the map is then simply a fixed point corresponding
to a periodic orbit of the flow. See \cite{DomingoG27} for a related work to ours
for fixed points.

The constant return time is a consequence of using an angular section.
A Poincar\'e section defined by an event in the physical variables generally
has a state-dependent return time, and the induced map need not have the rigid
phase translation in \eqref{eq.SP}, see \cite{GimenoJJMZ23}. Our algorithms use the skew-product
structure explicitly, so we either start from a map of the form
\eqref{eq.SP} or obtain one from a flow through the angular construction just
described. Related numerical formulations for reducible invariant tori and
their invariant manifolds in Poincar\'e maps can be found in
\cite{JorbaO09,GimenoJNO21}.

\subsection{Normal variational dynamics and reducibility}
\label{sec.torus-reducibility}

Let $\tor$ satisfy \eqref{eq.invariance_torus} and define
$
    A(\theta)\bydef
    \dop_x\poincmap\bigl(\tor(\theta),\theta\bigr)
$. 
Linearizing the fiber equation along the torus gives
\begin{equation}
\begin{cases}
\bar v=A(\theta)v,\\
\bar\theta=\theta+\rho \pmod{1}.
\end{cases}
\label{eq.LSP}
\end{equation}
This linear skew-product is the quasi-periodic cocycle associated with the
torus.

For a graph torus, there is a natural way to identify this cocycle with the
normal dynamics. At the point $(\tor(\theta),\theta)$, tangent vectors to
$\torinv$ have the form
$\bigl(\dop\tor(\theta)\eta,\eta\bigr)$, while the vertical space
$\mathcal V_\theta=\R^n\times\{0\}$ is a complementary subspace. The
skew-product structure preserves this vertical bundle because
\[
    \dop\expoincmap\bigl(\tor(\theta),\theta\bigr)(v,0)
    =\bigl(A(\theta)v,0\bigr).
\]
Thus the vertical bundle provides a canonical transverse representation of
the normal bundle, and $A(\theta)$ is the corresponding normal variational
cocycle. Here ``normal'' means transverse to the torus; no orthogonality is implied. For a torus that is not represented as a graph, one would instead
choose a transverse bundle or work directly with the quotient normal bundle, which is not in the scope of this work.

For a periodic orbit, Floquet theory uses a periodic change of variables to
replace the linearized dynamics by constant coefficients. For a
quasi-periodic cocycle this simplification is an additional property, which
we state explicitly.

\begin{dfn}[Reducibility]\label{def:reducible}
The cocycle \eqref{eq.LSP} is \emph{reducible} if there exist a smooth map
$C\colon\T^d\to\mathrm{GL}(n,\R)$ and a constant matrix
$B\in\R^{n\times n}$ such that the change of variables
$v=C(\theta)u$ gives $\bar u=Bu$. Equivalently,
\[
    A(\theta)C(\theta)=C(\theta+\rho)B,
    \qquad
    B=C^{-1}(\theta+\rho)A(\theta)C(\theta),
\]
and the last expression is independent of $\theta$. We call $C$ a Floquet
transformation, $B$ the Floquet matrix, and the eigenvalues of $B$ the
Floquet multipliers. We call the invariant torus reducible when its normal
cocycle is reducible.
\end{dfn}

Once the cocycle has been reduced, the phase dependence is carried entirely
by $C(\theta)$ and the linear dynamics is encoded by the single matrix $B$.
If $E\subset\C^n$ is a $B$-invariant subspace, then
$C(\theta)E$ is an invariant bundle along the torus. In particular, let
$E_B^s$, $E_B^c$, and $E_B^u$ be the sums of generalized eigenspaces of $B$
associated with multipliers of modulus smaller than, equal to, and larger
than one. The corresponding normal bundles are
\[
    \mathcal E^{\star}(\theta)=C(\theta)E_B^{\star},
    \qquad \star\in\{s,c,u\}.
\]
For a real system, nonreal multipliers occur in conjugate pairs and the two
complex eigendirections together determine a real invariant subbundle. This
is the reason complex spectral coordinates appear later, even though the
original map and the Floquet transformation are real. The normal center
bundle $\mathcal E^c$ should also be distinguished from the tangent bundle of
$\torinv$: it consists of transverse directions whose multipliers lie on the
unit circle.

This representation is especially convenient numerically. The invariant
splitting and its linear rates can be obtained from the spectrum of $B$, while
$C(\theta)$ transports those constant spectral subspaces along the torus.
Algorithms that compute the torus and its Floquet reduction simultaneously
are described in
\cite{JorbaO09,CanadellH2014,CanadellH2017,GimenoJNO21}; see
also \cite{HaroCLMF16} for the parameterization method more broadly.

In the computations of this manuscript, we assume that the torus and its Floquet reduction
have already been obtained. Their residuals are
\begin{equation}
\etor(\theta)
    \bydef \tor(\theta+\rho)-\poincmap(\tor(\theta),\theta),
\qquad
\ered(\theta)
    \bydef C^{-1}(\theta+\rho)
       \dop_x\poincmap(\tor(\theta),\theta)C(\theta)-B.
\label{eq.error0}
\end{equation}
Both vanish for an exact reducible torus. Numerically, we check that their
norms are below the prescribed tolerance before using $C$ and $B$ in the
graph-parameterization.

\subsubsection{Nonresonance and diophantine conditions in the torus and Floquet corrections}
\label{sec.torus-reducibility-conditions}

Definition~\ref{def:reducible} is a statement about the existence of a smooth
change of variables $C(\theta)$. It does not, by itself, impose a Diophantine
condition. Small divisors enter when the torus and the Floquet transformation
are found by solving linearized equations in Fourier series, as in Newton or
KAM schemes.

For a Fourier expansion
\[
    u(\theta)=\sum_{k\in\Z^d}\widehat u_k
    e^{2\pi\I\innerp{k}{\theta}},
\quad \text{
the shift by $\rho$ gives} \quad 
    u(\theta+\rho)=\sum_{k\in\Z^d}
    \bigl(e^{2\pi\I\innerp{k}{\rho}}\widehat u_k\bigr)
    e^{2\pi\I\innerp{k}{\theta}}.
\]
Thus each Fourier mode carries the factor
$e^{2\pi\I\innerp{k}{\rho}}$. Suppose first that $B$ is diagonalizable over
$\C$, with Floquet multipliers $\lambda_1,\ldots,\lambda_n$. The linearized
equation for the torus correction contains divisors of the form
$e^{2\pi\I\innerp{k}{\rho}}-\lambda_j$. A typical quantitative hypothesis is
\begin{equation}
\label{eq.torus-correction-nonresonance}
    \left|e^{2\pi\I\innerp{k}{\rho}}-\lambda_j\right|
        \ge \frac{\gamma}{|k|^{\tau}},
    \qquad k\in\Z^d\setminus\{0\}, \: j = 1, \dotsc, n.
\end{equation}
In a persistence KAM theorem this is one of the hypotheses used to solve for the
torus, together with regularity, smallness, and nondegeneracy assumptions, see \cite{Llave01}. For
the numerical Newton step, \eqref{eq.torus-correction-nonresonance} is the condition that controls the nonzero
Fourier modes of the torus correction. The zero mode is a finite-dimensional
solvability question: for example, $1\notin\Spec(B)$ makes it invertible,
whereas a neutral multiplier may require a parameter, a phase condition, or a
separate nondegeneracy argument. Jorba and Olmedo make the role of these
divisors explicit in \cite{JorbaO09}.

Correcting the Floquet transformation produces a different family of
divisors. A corresponding nonresonance condition is
\begin{equation}
\label{eq.floquet-reduction-nonresonance}
    \left|e^{2\pi\I\innerp{k}{\rho}}
        -\frac{\lambda_j}{\lambda_\ell}\right|
        \ge \frac{\gamma}{|k|^{\tau}},
    \qquad k\in\Z^d\setminus\{0\}, \: j, \ell = 1, \dotsc, n.
\end{equation}
Equivalently, the mode equations contain differences of the form
$\lambda_\ell e^{2\pi\I\innerp{k}{\rho}}-\lambda_j$, up to an interchange of
the two spectral indices. Condition~\eqref{eq.floquet-reduction-nonresonance} controls the Fourier
equation for the change of variables $C(\theta)$. It is therefore an additional
nonresonance hypothesis in perturbative or Newton schemes for Floquet
reduction, not a stronger existence condition for the torus itself. If a
numerical scheme computes a reducible torus and its Floquet transformation
simultaneously, both \eqref{eq.torus-correction-nonresonance} and
\eqref{eq.floquet-reduction-nonresonance} have to be controlled. The two
families nevertheless belong to different linear equations: a torus can
persist after the Floquet nonresonance fails, and its normal cocycle may then
be nonreducible. This distinction is also made explicit in \cite{JorbaO09};
numerical detection and dynamical consequences of nonreducibility are studied,
for example, in
\cite{Jorba01,VitoloBS2011,JorbaT08,JorbaNOT07, FiguerasH15, FagellaJJT19}. For a
concrete Floquet example where reducibility and nonreducibility are examined in
a quasi-periodically forced Hill equation, see \cite{BroerS98}.

The modulus of the multipliers already eliminates many possible small
divisors. If $|\lambda_j|\neq1$, then the left-hand side of
\eqref{eq.torus-correction-nonresonance} is uniformly separated from zero.
Likewise, if $|\lambda_j|\neq|\lambda_\ell|$, then
$\lambda_j/\lambda_\ell$ is off the unit circle and
\eqref{eq.floquet-reduction-nonresonance} is uniformly separated from zero.
The delicate cases are therefore multipliers on the unit circle and pairs of
multipliers with the same modulus. 

Our numerical examples use a diagonalizable
Floquet matrix with no multiple eigenvalues, since reducibility with multiple
eigenvalues is a separate subject in the perturbative literature; see, e.g. \cite{ZhangXWX2020}.

For the continuous-time system \eqref{eq.qpvf}, the analogous homological
equations are additive. If $\chi_j$ are Floquet exponents, the divisors involve
combinations such as
\[
    2\pi\I\innerp{k}{\tilde\omega}+\chi_j-\chi_\ell.
\]
After taking the angular Poincar\'e map, exponentiation turns these into the
multiplicative relations in
\eqref{eq.floquet-reduction-nonresonance}. This is the precise sense in which
the rotation vector and the normal Floquet spectrum both enter the
reducibility calculation. Classical perturbative results for continuous-time
quasi-periodic systems are given in \cite{JorbaS92, Eliasson01,Puig2006}; see also
\cite{Lang91,Llave01} for Diophantine estimates.

The graph-parameterization introduces another set of resonances at higher Taylor
orders. Those divisors contain products of the Floquet multipliers selected as
independent variables and control the nonlinear invariant-manifold equations.
They are derived in Section~\ref{sec.solve.cohom}. Thus there are three related
but distinct questions: solving for the torus, reducing its linear normal
cocycle, and solving the nonlinear graph equations over the reduced torus. Each of them infer stronger constraints.

\section{Graph-Parameterization for a Reducible Invariant Torus}\label{sec.graph}

We now construct a local invariant manifold attached to a reducible torus.
When the fiber map in \eqref{eq.SP} is locally invertible, the forward and
backward dynamics can be handled with the same equations. This is useful
numerically: unstable manifolds are naturally computed forward, whereas
stable manifolds can be computed using the inverse map. A sign
$\sigma\in\{+1,-1\}$ will indicate the chosen orientation.

The construction is a graph-form parameterization of the invariant manifold;
see, for example, \cite{HaroCLMF16}. We choose a $B$-invariant spectral subspace as the space of independent reduced variables, and seek an invariant manifold that is locally represented as a graph over this subspace.

After the Floquet reduction, the spectral calculations are performed in
$\C^n$ so that complex-conjugate eigendirections can be treated directly. The
original map remains real; the corresponding reality conditions are developed
in Section~\ref{sec.hermitian} separately.

\subsection{Parameterization and Directed Invariance Equation}\label{sec.inv-eq}

For $\sigma \in \{+1,-1\}$, assume that $\poincmap$ in \eqref{eq.SP} is locally invertible in a neighborhood of the invariant torus $\tor$. We define the directed fiber map and rotation vector by 
\begin{equation*}
 \rho _\sigma \bydef \sigma \rho, \qquad 
 \poincmap _\sigma (x, \theta) \bydef 
 \begin{cases}
  \poincmap (x, \theta) & \text{if }\sigma = +1 \\
  \bigl(\poincmap(\cdot, \theta - \rho) ^{-1}\bigr)(x) & \text{if }\sigma = -1.
 \end{cases}
\end{equation*}
The directed skew-product is then $\expoincmap _\sigma (x,\theta) \bydef \bigl(\poincmap _\sigma (x,\theta), \theta + \rho _\sigma \bigr)$. Notice that $\poincmap _{-1}$ incorporates the phase shift associated with the inversion of the full skew-product. 

A torus $\tor$ is invariant for either orientation when $\poincmap _\sigma \bigl(\tor(\theta), \theta \bigr) = \tor(\theta + \rho _\sigma)$. If $A _\sigma(\theta) \bydef \dop _x \poincmap _\sigma (\tor(\theta), \theta)$, the reducibility identity is read as:
\begin{equation}
    \label{eq.GT.directed-reducibility}
 A _\sigma (\theta) C(\theta) = C(\theta + \rho _\sigma) B ^\sigma.
\end{equation}
For $\sigma=-1$, this identity follows by evaluating the forward relation at $\theta-\rho$ and then inverting it; see Appendix~\ref{app.stable}. Thus the same Floquet transformation $C(\theta)$ is used in both directions, while $B$ is replaced by $B^{-1}$ for the inverse dynamics. Computationally, we work with $B$ and avoid a direct invertibility expression to prevent numerical stability issues.

\smallskip

We seek a local invariant manifold attached to the torus. Let $\mathcal{D} \subset \C ^{\nt}$, with $\nt \leq n$, be a neighborhood of the origin and let $\st \in \mathcal{D}$ denote the reduced coordinates. We consider an embedding $\param \colon \mathcal{D} \times \T ^d \to \C ^n$.

\begin{dfn}
    \label{dfn:GT.directed-PM}
The image of an embedding $\param$ is invariant under the directed skew-product $\expoincmap _\sigma$ if, and only if, there exists an internal map $\intdyn _\sigma \colon \mathcal{D} \times \T ^d \to \C ^{\nt}$ such that
\begin{equation} 
    \label{eq.GT.directed-PM}
 \poincmap _\sigma \bigl(\param(\st, \theta), \theta\bigr) = 
 \param \bigl(\intdyn _\sigma(\st, \theta), \theta + \rho _\sigma \bigr), \qquad \text{for all } (\st,\theta).
\end{equation}
\end{dfn}

\begin{dfn}
    \label{dfn:GT.zero-section}
Let $\param\colon \mathcal{D} \times \T ^d \to \C ^n$ be an embedding with internal map $\intdyn_\sigma$. The manifold $\param(\mathcal{D} \times \T ^d)$ is attached to an invariant torus $\tor$ when 
\begin{equation}
    \label{eq.GT.zero-section}
 \param(0, \theta) = \tor(\theta) \qquad \text{and} \qquad 
 \intdyn _\sigma(0, \theta) = 0, \qquad \text{for all } \theta.
\end{equation}
\end{dfn}

\smallskip

Equation \eqref{eq.GT.directed-PM} is the parameterization equation in the directed skew-product dynamics.

\subsection{Reducible Adapted Coordinates and Spectral Selection}\label{sec.adapt-frame}

The constant Floquet matrix $B$ contains the linear normal dynamics of the torus $\tor$. We assume that the fiber space admits a $B$-invariant splitting:
\begin{equation}
  \label{eq.GT.LN-splitting}
 \C ^n = \Et \oplus \Eh, \qquad 
 \dim \Et = \nt, \qquad 
 \dim \Eh = \nh \bydef n  - \nt.
\end{equation}
where $\Et$ is the spectral subspace over which the desired manifold will be written as a graph and $\Eh$ is its complete complementary spectral subspace.
Since the Floquet-reduced fiber already contains only directions
transverse to the torus, the subscript $N$ denotes the \emph{complement of the
selected spectral subspace} within this fiber, not the full normal bundle of the
torus.

Let $V_L\in\C^{n\times \nt}$ and
$V_N\in\C^{n\times \nh}$ be bases of these subspaces.  We write
\begin{equation}
    \label{eq.GT.spectral-blocks}
    V=
    \begin{bmatrix}
     V _L & V _N
    \end{bmatrix},
    \qquad
    V^{-1}BV
    =\Lambda
    =\begin{pmatrix}
        \Lambda_L & 0\\
        0 & \Lambda_N
      \end{pmatrix}.
\end{equation}
The splitting only requires that $\Et$ and $\Eh$ are $B$-invariant. In the
complex formulation, we choose eigen-coordinates in $\Et$, so that
$\Lambda_L$ is diagonal; $\Lambda_N$ may remain a matrix block. We do not
treat a non-diagonalizable selected block, since that would require a block
version of the homological equations and is not covered by the numerical
experiments. If $\Lambda_N$ is also diagonal, the cohomological equations
decouple completely into scalar Fourier equations; see
Section~\ref{sec.solve.cohom}.

\begin{rmk}[Choice of spectral coordinates]
The construction needs a $B$-invariant splitting, not a particular
normalization of the eigenvectors. When diagonalization is poorly conditioned,
a complex Schur decomposition can provide orthonormal bases for invariant
spectral blocks. QR or SVD factorizations can be useful to recondition a basis
inside an already selected invariant subspace, but they do not by themselves
produce the invariant spectral splitting required in
\eqref{eq.GT.spectral-blocks}.
\end{rmk}

\smallskip

The adapted frame along the reducible invariant torus is
\begin{equation}
    \label{eq.GT.adapted-frame}
    \adframe(\theta)\bydef C(\theta)V, \qquad \text{ that is} \qquad 
    \begin{cases}
     \adframe_L(\theta) = C(\theta)V_L \\
     \adframe_N(\theta) = C(\theta)V_N.
    \end{cases}
\end{equation}
Combining \eqref{eq.GT.directed-reducibility} and
\eqref{eq.GT.spectral-blocks} gives
\begin{equation}
    \label{eq.GT.directed-adapted-frame}
    A_\sigma(\theta)\adframe(\theta)
    =\adframe(\theta+\rho_\sigma)\Lambda_\sigma,
    \qquad
    \Lambda_\sigma\bydef \Lambda^\sigma,
\end{equation}
where
\begin{equation}
    \label{eq.GT.directed-blocks}
    \Lambda_\sigma
    =\begin{pmatrix}
        \Lambda_{L,\sigma} & 0\\
        0 & \Lambda_{N,\sigma}
      \end{pmatrix},
    \qquad
    \Lambda_{L,\sigma}\bydef \Lambda_L^\sigma,
    \qquad
    \Lambda_{N,\sigma}\bydef \Lambda_N^\sigma.
\end{equation}
Thus the same adapted frame works in both directions; only the phase shift and
the constant linear dynamics are reversed.

\smallskip

The splitting \eqref{eq.GT.LN-splitting} is functional rather than merely
spectral notation.  If the Floquet spectrum of $B$ is classified as stable, center,
and unstable,
\begin{equation*}
    E^s\oplus E^c\oplus E^u,
    \label{eq.GT.scu-splitting}
\end{equation*}
then the two blocks in \eqref{eq.GT.LN-splitting} are assembled according to
the manifold being computed.  Typical choices are
\begin{equation*}
    \label{eq.GT.standard-choices}
\begin{array}{@{}c||c|c|c||c|c@{}}
\textsc{manifold} & \sigma & \Et & \Eh & \poincmap _\sigma & \Lambda _{L,\sigma}\\ \hline \hline
\text{unstable} & +1 & E^u & E^s\oplus E^c & \poincmap & \Lambda _s\\[1mm]
\text{stable}   & -1 & E^s & E^c\oplus E^u & \poincmap^{-1} & \Lambda _u^{-1}\\[1mm]
\text{center}   & +1 & E^c & E^s\oplus E^u & \poincmap & \text{unknown } \intdyn _c
\end{array}
\end{equation*}
whose invariant equations are:
\begin{equation*}
 \begin{aligned}
  \poincmap _\sigma \bigl(\param ^\sigma(\st, \theta), \theta \bigr) &= \param ^\sigma( \Lambda _{L,\sigma}\st, \theta + \rho _\sigma), & \sigma &\in \{+1,-1\} \\
  \poincmap \bigl(\param^c(\st, \theta), \theta \bigr) &= \param^c( \intdyn _c (\st, \theta), \theta + \rho), &\param^c &\text{ center and $\intdyn _c$ unknown.}
 \end{aligned}
\end{equation*}

More generally, $\Et$ may be any selected invariant spectral subspace for
which the required non-resonance and spectral-gap conditions hold.  For
example, one may select a slow stable bundle while leaving stronger stable
directions in $\Eh$ \cite{HaroLl05c,HaroLl06b}.

\begin{rmk}[On choice of $\sigma$]
 For a stable manifold, we take $\sigma=-1$, so that the stable
directions become expanding directions of the directed map.  This convention is numerically motivated for jet transport, since it prevents the coefficients of a stable parameterization from being contaminated by forward amplification of transverse unstable errors.  The use of the inverse Poincar\'e map for this
purpose is standard in high-order computations; see \cite{CapinskiM2017,GimenoJNO21}.  The role of $\sigma$ is therefore organizational: it allows all formulas to be written once. Mathematically, one can just keep $\sigma = +1$.
\end{rmk}

\subsection{Graph Representation and Normalization}\label{sec.norma-cond}

For each $\theta\in\T^d$, introduce the affine adapted coordinate map
\begin{equation*}
    \Psi_\theta(\st,\sh)
    \bydef\tor(\theta)+\adframe(\theta) 
    \begin{pmatrix}
     \st \\ \sh
    \end{pmatrix}
    =\tor(\theta)+\adframe _L(\theta) \st +\adframe _N(\theta)\sh,
    \qquad
    (\st,\sh)\in\C^{\nt}\times\C^{\nh}.
    \label{eq.GT.Psi}
\end{equation*}
The fiber dynamics of $\poincmap _\sigma$ in these coordinates is
\begin{equation}
    \label{eq.GT.reduced-map}
    \cR_\sigma(\st,\sh,\theta)
    \bydef \adframe(\theta+\rho_\sigma)^{-1}
      \left[
        \poincmap _\sigma\bigl(\Psi_\theta(\st,\sh),\theta\bigr)
        -\tor(\theta+\rho_\sigma)
      \right].
\end{equation}
We split in components 
$
    \cR_\sigma
    =(
        \cR_{L,\sigma},
        \cR_{N,\sigma}
      )
$.
By torus invariance and frame reducibility \eqref{eq.GT.directed-adapted-frame},
$$
    \cR_\sigma(0,0,\theta)=0,
    \qquad
    \dop _{(\st,\sh)}\cR_\sigma(0,0,\theta)=\Lambda _\sigma, \qquad \text{for all } \theta.
$$
Consequently, there exist nonlinear remainders $\remainL _\sigma$ and $\remainN _\sigma$
with vanishing constant and linear terms described by a discrete system:
\begin{equation}
\label{eq.GT.adapted-dynamics}
\begin{cases}
    \bar \st
      =\cR_{L,\sigma}(\st,\sh,\theta)
       =\Lambda_{L,\sigma}\st+ \remainL _\sigma(\st,\sh,\theta),
       \\
    \bar \sh
      =\cR_{N,\sigma}(\st,\sh,\theta)
       =\Lambda_{N,\sigma}\sh+ \remainN _\sigma(\st,\sh,\theta),
       \\
    \bar\theta
      =\theta+\rho_\sigma,
\end{cases}
\end{equation}
where for all $\theta$, $\remainL _\sigma(\st,\sh,\theta) $ and $\remainN _\sigma(\st,\sh,\theta)$ are $\mathcal O\bigl(\|(\st,\sh)\|^2\bigr)$.

\smallskip

We represent the desired invariant manifold as a graph
\begin{equation}
    \label{eq.GT.graph}
    \sh=\graphfun(\st,\theta),
    \qquad
    \graphfun\colon \mathcal D\times\T^d\to \C^{\nh}.
\end{equation}
Thus, $\st$ denotes independent variables and $\sh$ the dependent ones. The corresponding parameterization in the original fiber variables is
\begin{equation}
    \label{eq.GT.K-graph}
    \param(\st,\theta)
    =\tor(\theta)
     +\adframe _L(\theta)\st
     +\adframe _N(\theta)\graphfun(\st,\theta).
\end{equation}
Because the $\st$-coordinate is fixed by the graph representation, the internal
dynamics is not prescribed a priori.  It is induced by the $L$-component of
the transformed map:
\begin{equation}
    \label{eq.GT.induced-dynamics}
    \intdyn _\sigma(\st,\theta)
    \bydef\cR_{L,\sigma}\bigl(\st,\graphfun(\st,\theta),\theta\bigr)
    =\Lambda_{L,\sigma}\st
     +\remainL _\sigma\bigl(\st,\graphfun(\st,\theta),\theta\bigr).
\end{equation}
The graph is invariant if, and only if, $\bar \sh = \graphfun (\bar \st, \bar \theta)$. That is, 
\begin{equation}
    \label{eq.GT.graph-invariance-abstract}
    \cR_{N,\sigma}\bigl(\st,\graphfun(\st,\theta),\theta\bigr)
    =\graphfun\bigl(\intdyn _\sigma(\st,\theta),\theta+\rho_\sigma\bigr).
\end{equation}
Using \eqref{eq.GT.adapted-dynamics}, this becomes the explicit graph equation
\begin{equation}
    \label{eq.GT.explicit-graph-equation}
    \Lambda_{N,\sigma}\graphfun(\st,\theta)
    +\remainN _\sigma\bigl(\st,\graphfun(\st,\theta),\theta\bigr)
    =\graphfun\!\left(
        \Lambda_{L,\sigma}\st
        +\remainL _\sigma\bigl(\st,\graphfun(\st,\theta),\theta\bigr),
        \theta+\rho_\sigma
      \right).
\end{equation}
Equations \eqref{eq.GT.K-graph}--\eqref{eq.GT.explicit-graph-equation}
then imply the directed parameterization equation
\eqref{eq.GT.directed-PM} identically.

The normalization of the graph is
\begin{equation}
    \label{eq.GT.normalization}
    \graphfun(0,\theta)=0,
    \qquad
    \dop _\st \graphfun(0,\theta)=0,
    \qquad \text{for all }\theta.
\end{equation}
The first condition anchors the graph at the prescribed torus.  The second
ensures that the tangent space of the manifold along the torus is precisely the
selected bundle:
$
    \dop _\st \param(0,\theta)=\adframe _L(\theta)
$. 
In particular, $
    \dop _\st \intdyn _\sigma(0,\theta)=\Lambda_{L,\sigma} $.

\begin{rmk}[On the normalization condition choice]
The normalizations \eqref{eq.GT.normalization} also set the high-orders.  A general parameterization of the same geometric manifold may be
changed by a near-identity reparameterization of $\st$.  Our graph normalization choice fixes
this freedom by requiring that the $L$-component in adapted coordinates be
exactly $\st$:
\begin{equation*}
    \pi_L \adframe (\theta)^{-1}\bigl(\param(\st,\theta)-\tor(\theta)\bigr)=\st,
    \label{eq.GT.graph-normalization}
\end{equation*}
where $\pi_L$ denotes projection onto the first $\nt$ adapted coordinates.
Thus all nonlinear corrections to the parameterization lie in the complete
complementary bundle $\Eh$, while the nonlinear terms of the tangential
equation are absorbed into the reduced dynamics $f_\sigma$.
\end{rmk}

\section{Construction of the Graph}
\label{sec.GT-construction}

We now derive the order-by-order equations to compute $\graphfun$ in \eqref{eq.GT.graph}.  The
derivation is independent of the particular representation chosen for the
coefficients.  Fourier series in $\theta$ and multivariate Taylor series in
$\st$ will be used in the implementation, while the nonlinear action of a
Poincar\'e map is evaluated by jet transport.

\subsection{Representation of the Unknowns}
\label{sec.GT-fourier-taylor}

For an analytic function, or sufficiently differentiable, in $\st$ at the
origin, normalization \eqref{eq.GT.normalization} gives the formal
expansion
\begin{equation}
    \label{eq.GT.W-FPS}
    \graphfun(\st,\theta)
    =\sum_{l\ge2}\graphfun_{(l)}(\st,\theta),
    \qquad
    \graphfun_{(l)}(\st,\theta)
    =\sum_{|\alpha|=l}w_{l,\alpha}(\theta)\st^\alpha,
\end{equation}
where each $\graphfun_{(l)}$ is a homogeneous polynomial of degree $l$ from
$\C^{\nt}$ to $\C^{\nh}$.  With multi-index notation
$\alpha=(\alpha_1,\ldots,\alpha_{\nt})\in\N^{\nt}$,
$|\alpha|=\alpha_1+\cdots+\alpha_{\nt}$, and
$\st^\alpha=\st _1^{\alpha_1}\cdots \st _{\nt}^{\alpha_{\nt}}$.
There are
\begin{equation*}
    \label{eq.GT.number-monomials}
    N_l=\binom{\nt+l-1}{l}
\end{equation*}
monomials of total degree $l$.  Each coefficient is quasi-periodic in the
base angle and is represented by a Fourier series
\begin{equation}
    \label{eq.GT.W-Fourier}
    w_{l,\alpha}(\theta)
    =\sum_{k\in\mathbb Z^d}
      \widehat w_{l,\alpha,k}\,
      e^{2\pi\I\,\innerp{k}{\theta}}.
\end{equation}
In numerical computations both sums \eqref{eq.GT.W-FPS} and \eqref{eq.GT.W-Fourier} are truncated.

The internal dynamics has the corresponding expansion
\begin{equation}
    \label{eq.GT.f-series}
    \intdyn _\sigma(\st,\theta)
    =\Lambda_{L,\sigma}\st
     +\sum_{l\ge2}f_{\sigma,(l)}(\st,\theta).
\end{equation}
Unlike in a parameterization with prescribed internal dynamics, the
higher-order terms $f_{\sigma,(l)}$ are determined by the graph invariance
equation \eqref{eq.GT.graph-invariance-abstract} together with the normalization
conditions \eqref{eq.GT.normalization}. They need not to 
vanish, and in center or mixed problems they contain the reduced nonlinear
dynamics of interest.

For $l\geq 2$, denote by
\begin{equation*}
    \label{eq.GT.W-less-j}
    \graphfun^{<l}(\st,\theta)\bydef\sum_{q=2}^{l-1}\graphfun _{(q)}(\st,\theta)
\end{equation*}
the already computed part of the graph.  The recursive step assumes that the
invariance equation is satisfied through degree $l-1$ and determines
$\graphfun_{(l)}$ and $\intdyn _{\sigma,(l)}$ from quantities depending only on lower-order
coefficients.
The recursive construction is obtained by power matching homogeneous Taylor terms in
\eqref{eq.GT.graph-invariance-abstract}.

\subsection{The Online Algorithm}

The following algorithm presents how to compute the terms of the
parameterization \eqref{eq.GT.W-Fourier} and the internal dynamics
\eqref{eq.GT.f-series} from a torus and its Floquet reduction provided by
\eqref{eq.error0}.

\begin{palg}[Computation of a graph-parameterization for a reducible invariant torus] 
\label{palg.graph-tor} \
\begin{enumerate}
\itemsep=0em
    \item [$\star$] \texttt{Input:} Directed fiber map $\poincmap _\sigma$ and rotation $\rho_\sigma$ defined in Section~\ref{sec.inv-eq}, with $\sigma \in \{-1,1\}$;
      invariant torus $\tor(\theta)$; Floquet transformation $C(\theta)$; Floquet matrix $B$; and rotation vector $\rho$.
      A nonzero scaling factor $\mathfrak{p}$ may be prescribed; its default value is $1$. 
    \item [$\star$] \texttt{Output:} Graph function $\graphfun(\st,\theta)$ in the form \eqref{eq.GT.W-FPS} and internal dynamics $\intdyn _{\sigma}(\st,\theta)$, satisfying the invariance equation \eqref{eq.GT.graph-invariance-abstract} and the normalization \eqref{eq.GT.normalization}.
\end{enumerate}
\begin{enumerate}
\itemsep=0em
\renewcommand{\theenumi}{\arabic{enumi}}
\setcounter{enumi}{-8}

\item Choose a $B$-invariant splitting $BV=V\Lambda$ as in
\eqref{eq.GT.spectral-blocks}, with the selected block $\Lambda_L$ \hfill {\footnotesize $\rhd$ for real see Sec.~\ref{sec.canonical-conjugate-pairs}}

\item \label{palg.graph-tor.scaling} $\Vt \gets \mathfrak{p} \Vt $ \hfill {\footnotesize $\rhd$scaling factor}

\item $\rho _\sigma \gets \sigma \rho$ and $\Lambda _\sigma \gets \Lambda ^\sigma$

\item $\adframe (\theta) \gets C(\theta) V$

\item Factorize $Q_\sigma^+(\theta) \bydef
C(\theta+\rho_\sigma)V$ to solve repeated linear systems

\item $\graphfun (\st,\theta) \gets 0$ and $R(\st, \theta) \gets 0$ \hfill {\footnotesize $\rhd$\sl normalization conditions on $\graphfun$ and auxiliary jet}

\item $\intdyn _{\sigma} (\st,\theta) \gets \Lambda _{L,\sigma} s$ \hfill {\footnotesize $\rhd$\sl internal dynamics initialization}

\end{enumerate}
For $l = 2,3, \dotsc $\hfill {\footnotesize $\rhd$\sl only operations up to order $l$ are needed}
\begin{enumerate}
\itemsep=0em
\renewcommand{\theenumi}{\arabic{enumi}}

\item $\graphfun _{(l)}(\st,\theta) \gets 0$ \hfill {\footnotesize $\rhd$\sl the degree-$l$ coefficient is the unknown}

\item $\param (\st,\theta)
\gets \tor(\theta) + \adframe _L(\theta)\st +\adframe _N(\theta) \graphfun(\st,\theta)$\hfill {\footnotesize $\rhd$\sl $\graphfun ^{<l}$ already known}

\item \label{palg.graph-tor.jet} Solve $\cR$ in $Q ^+_\sigma(\theta) \cR (\st,\theta) = \poincmap _\sigma \bigl(\param (\st, \theta), \theta\bigr) -\tor(\theta + \rho _\sigma)$ \hfill {\footnotesize $\rhd$\sl jet transport step and linear solver up to order $l$}

\item $
\begin{bmatrix}
    \cR _{L} (\st,\theta) & \cR _{ N} (\st,\theta) 
\end{bmatrix} \gets
\begin{pmatrix}
\idmat _{\nt} \\ & \idmat_{\nh}
\end{pmatrix}
\cR(\st,\theta)$ \hfill {\footnotesize $\rhd$\sl component blocks}

\item \label{palg.graph-tor.intdyn} $\intdyn _{\sigma, (l)}(\st,\theta) \gets \cR _{L, (l)}(\st,\theta)$ \hfill {\footnotesize $\rhd$\sl $l$-th internal dynamic update}

\item \label{palg.graph-tor.comp} $Z(\st,\theta) \gets R\bigr (\cR _{ L}(\st,\theta), \theta \bigr)$\hfill {\footnotesize $\rhd$\sl composition of jets up to order $l$}

\item $Z _{(l)}(\st, \theta) \gets Z _{(l)}(\st, \theta) - \Upsilon_{N,(l)} (\st, \theta)$\hfill {\footnotesize $\rhd$\sl subtraction on the $l$th homogeneous polynomial}

\item [\theenumi\textsuperscript{bis})] Transform $Z _{(l)}(\st, \theta)$ to complex. See Section~\ref{sec.num-real}\hfill {\footnotesize $\rhd$\sl {\sc optional} if all is being computed in real}

\item \label{palg.graph-tor.cohom} For every multiindex $\alpha \in \N ^{\nt}$ such that $|\alpha| = l$:
\begin{enumerate}
\itemsep=0em
\renewcommand{\theenumii}{\theenumi.\arabic{enumii}}
\renewcommand{\labelenumii}{\theenumii}
\item 
Let $\lambda_{L,\sigma,1},\ldots,\lambda_{L,\sigma,\nt}$ be the complex 
diagonal entries of $\Lambda_{L,\sigma}$ and set
$\beta_{\sigma,\alpha}\bydef
\prod_{i=1}^{\nt}\lambda_{L,\sigma,i}^{\alpha_i}$. Find
$w_{l,\alpha}(\theta)$ so that
$\Lambda _{N,\sigma} w _{l,\alpha}(\theta)
-w _{l,\alpha}(\theta + \rho _\sigma)\beta_{\sigma,\alpha}
=z _{l,\alpha}(\theta)$
\hfill {\footnotesize $\rhd$\sl cohomological eq.}

For this, an option is to expand the functions in Fourier series, see Lemma~\ref{lem.cohom}.

\end{enumerate}

\item [\theenumi\textsuperscript{bis})] Transform $\graphfun _{(l)}(\st, \theta)$ to real. See Section~\ref{sec.num-real}\hfill {\footnotesize $\rhd$\sl {\sc optional} if all is being computed in real}

\item \label{palg.graph-tor.rot} $R _{(l)} (\st, \theta) \gets \graphfun _{(l)}(\st, \theta + \rho _\sigma)$ \hfill {\footnotesize $\rhd$\sl rotation on the $l$th homogeneous polynomial}

\end{enumerate}
\end{palg}

Algorithm~\ref{palg.graph-tor} computes the graph \eqref{eq.GT.W-FPS}
recursively in the Taylor degree. At degree $l$, all lower-order coefficients are
already known, and the only new unknown is $\graphfun_{(l)}$. The linear solves
use the same reduced normal dynamics at every order, so the phase-dependent
matrices $C(\theta)V$ and the factorizations of
$C(\theta+\rho_\sigma)V$ can be prepared once and reused.

Step~\ref{palg.graph-tor.scaling} rescales the basis $V_L$ uniformly. This
does not change the geometric manifold; it only changes the scale of the local
coordinates $\st$. Numerically, a suitable scale helps avoid very large or very
small Taylor coefficients. Section~\ref{sec.scaling} describes the heuristic
used to choose $\mathfrak p$. In practice, we first run
Algorithm~\ref{palg.graph-tor} with $\mathfrak p=1$, estimate a better scale
from the resulting graph, and rerun the computation with that value. This makes
statements about a normalized $\nt$-dimensional domain in the $\st$ variables
more meaningful.

Step~\ref{palg.graph-tor.jet} evaluates $\poincmap _\sigma$ on truncated multivariate
polynomial jets in the $\nt$ selected variables. If $\poincmap _\sigma$ is a
Poincar\'e map generated by a flow or semi-flow, these jets can be propagated by
high-order variational or jet-transport techniques
\cite{GimenoJJMZ23,Forrier23,GimenoJZ22}.
At degree $l$, it is enough to propagate jets truncated at that same degree;
there is no need to carry higher-order coefficients through the nonlinear map.
After substituting $\sh=\graphfun(\st,\theta)$, the components of
$\cR$ determine the nonlinear remainders in
\eqref{eq.GT.adapted-dynamics} through
\begin{equation*} \label{eq.fg}
\begin{split}
   \remainL _\sigma \bigl(\st, \graphfun( \st, \theta), \theta \bigr)
   &= \cR _{L} (\st, \theta) - \Lambda _{L,\sigma} \st ,\\
   \remainN _\sigma \bigl(\st, \graphfun( \st, \theta), \theta \bigr)
   &= \cR _{N} (\st, \theta)
      - \Lambda _{N,\sigma} \graphfun( \st, \theta) .
\end{split}
\end{equation*}

Step~\ref{palg.graph-tor.intdyn} updates the internal dynamics:
$\cR _{L}$ contains the whole internal dynamics $\intdyn _{\sigma}$ up to
order $l$, so we can optionally update the $l$th degree of $\intdyn _{\sigma}$ at this point.

Step~\ref{palg.graph-tor.comp} composes the jets.
The graph invariance equation \eqref{eq.GT.graph-invariance-abstract} requires
\[
\cR_{N}(\st,\graphfun(\st,\theta),\theta)
=\graphfun(\cR_{L}(\st,\graphfun(\st,\theta),\theta),
\theta+\rho_\sigma),
\]
and this step forms the right-hand side from the coefficients already known at lower degree. 
The mismatch with $\cR_{N}$ is then the
forcing term for the degree-$l$ cohomological equation, and also provides a
natural residual check for the previously computed orders. 

Step~\ref{palg.graph-tor.cohom} grows combinatorially with the order $l$
through the $\binom{\nt+l-1}{l}$ monomials of degree $l$. Because
$\Lambda_{L,\sigma}$ is diagonal, substitution of the linear internal map
into the monomial $s^\alpha$ contributes the scalar factor
$\beta_{\sigma,\alpha}=\prod_i
\lambda_{L,\sigma,i}^{\alpha_i}$. A non-diagonalizable selected block would
mix monomials of the same degree and leads to the block formulation of
Section~\ref{sec.solve.cohom}. With a diagonal selected block, different
multiindices $|\alpha|=l$ are independent; after a Fourier transform, the
Fourier modes are independent as well. See Lemma~\ref{lem.cohom}.

Step~\ref{palg.graph-tor.rot} involves 
a fixed rotation by $\rho$. If we mesh the
periodic coefficients in a uniform grid of size $M$ in $\T ^d$, the rotation is
performed in $O(M \log M)$ by standard Fast Fourier Transform, see
Section~\ref{sec.ftrep}

Finally, most operations in Algorithm~\ref{palg.graph-tor} are independent
for each phase point $\theta$ and can therefore be parallelized. Synchronization
is needed only when one stage uses data assembled over the full phase grid; for
example, the Fourier/cohomological solve cannot start until the jet data have
been computed at all phase points. Section~\ref{sec.numtests} describes
the residual and validity checks that can be inserted between these stages.

\begin{rmk}[Parameter dependence]
Algorithm~\ref{palg.graph-tor} can also be extended to formal expansions in a
system parameter. One possibility is to augment the phase space with a trivial
parameter equation. Another is to introduce the parameter as an additional
formal variable in the jets. In the latter approach, the input data and the
spectral quantities must be expanded consistently in that variable. These
expansions can be obtained systematically from the Implicit Function Theorem
\cite[Sec.~2.3]{GimenoJJZ25}, together with standard perturbation results for
simple eigenvalues such as \cite{Sun85}. Different truncation orders may then
be used for the state and parameter variables.
\end{rmk}

\begin{rmk}[Additional structure]
The basic algorithm uses smoothness, a reducible invariant torus, and a
$B$-invariant spectral splitting. If the map has additional structure, such as
symplectic or conformally symplectic geometry, some linear-algebra and
cohomological steps can often be simplified. We do not use those
structure-specific reductions in the current work, and leave it for future work.
\end{rmk}

\section{Solvability of the Cohomological Equations and Non-Resonance Conditions}
\label{sec.solve.cohom}

The only non-local step in Algorithm~\ref{palg.graph-tor} is Step~\ref{palg.graph-tor.cohom}, where the coefficients of the
homogeneous term $\graphfun _{(l)}$ are obtained by solving a family of cohomological
equations over the torus. These equations encode the compatibility between the
rotation $\theta \mapsto \theta+\rho _\sigma$ on the base and the linear normal dynamics
given by the eigenvalues of the Floquet matrix. Consequently, their solvability is
controlled by resonance relations between the internal frequencies $\rho$ and the
products of Floquet multipliers that appear at each Taylor order.

There are two distinct issues to separate. First, one must exclude exact resonances (i.e. multi-index exponents that make a denominator zero), since they may prevent the formal solvability of the homological equations at a
given order. Second, even in the absence of exact resonances, small divisors may
appear when the corresponding denominators become arbitrarily small. The latter
phenomenon affects the regularity estimates of the solution and, ultimately, the
convergence or effective accuracy of the computed expansion. In the hyperbolic
case, where the relevant Floquet multipliers have different moduli, the denominators
are uniformly bounded away from zero and no Diophantine condition may be needed. In
contrast, when the competing multipliers have the same modulus, one needs a
Diophantine-type non-resonance condition adapted to the products of Floquet
multipliers generated by the multiindices of the Taylor expansion.

We now make this discussion explicit in Fourier space, through Neumann-series representations, and in block form through Sylvester equations. These give a mode-by-mode description of the cohomological equation appearing in Step~\ref{palg.graph-tor.cohom} and provide the precise non-resonance conditions required by the graph-parameterization algorithm.

\subsection{Diagonal Formulation}

Assume first that the selected and complementary Floquet blocks in
\eqref{eq.GT.spectral-blocks} are diagonal (possibly over $\C$). At a fixed Taylor
multiindex, the action of the selected variables contributes a single product
of selected Floquet multipliers. It is notationally convenient to denote that product by
$\beta$ and keep it separate from the complementary multiplier $\mu$.

\begin{lem}[Diagonal Fourier solver] \label{lem.cohom}
Let $\rho\in\R^d$ and $\mu,\beta\in\C$. Consider
\begin{equation}\label{eq.cohom-graph}
    \mu w(\theta)-w(\theta+\rho)\beta=z(\theta),
    \qquad \theta\in\T^d,
\end{equation}
where $z$ and $w$ have Fourier expansions
\[
z(\theta)=\sum_{k\in\Z^d}\hat z_k e^{2\pi \I \innerp{k}{\theta}},
\qquad
w(\theta)=\sum_{k\in\Z^d}\hat w_k e^{2\pi \I \innerp{k}{\theta}}.
\]
Then \eqref{eq.cohom-graph} is equivalent, mode by mode, to
\begin{equation}\label{eq.cohom-graph-fourier}
\bigl(\mu-e^{2\pi \I\innerp{k}{\rho}}\beta\bigr)\hat w_k
=
\hat z_k,
\qquad k\in\Z^d.
\end{equation}
If $|\beta|\neq|\mu|$, then
\[
\bigl|\mu-e^{2\pi \I\innerp{k}{\rho}}\beta\bigr|
\ge \bigl||\mu|-|\beta|\bigr|,
\qquad k\in\Z^d,
\]
so no small divisors occur and $w$ is uniquely determined.
\end{lem}

\paragraph{Equal-modulus case.}
If $|\mu|=|\beta|\neq0$, a Diophantine bound of the form
\[
\left|e^{2\pi \I\innerp{k}{\rho}}-\frac{\mu}{\beta}\right|
\ge \gamma |k|^{-\tau},
\qquad k\in\Z^d\setminus\{0\},
\]
controls the nonzero Fourier modes. The zero mode must be treated separately;
in particular, the exact resonance $\mu=\beta$ requires a compatibility or
normal-form condition. 

For Algorithm~\ref{palg.graph-tor}, let
$\Lt=\diag(\lambda_i)_{i=1}^{\nt}$ and
$\Lh=\diag(\mu_j)_{j=1}^{\nh}$ be the diagonal blocks in
\eqref{eq.GT.spectral-blocks}. For a multiindex
$\alpha=(\alpha_1,\ldots,\alpha_{\nt})\in\N^{\nt}$, set
\[
    \beta^\alpha\bydef
    \lambda_1^{\alpha_1}\cdots\lambda_{\nt}^{\alpha_{\nt}}.
\]
At Taylor order $l=|\alpha|$, the scalar denominator is therefore
\begin{equation}\label{non-ded}
    \mu_j-e^{2\pi\I\innerp{k}{\rho}}\beta^\alpha,
    \qquad k\in\Z^d,
    \quad |\alpha|=l.
\end{equation}
The coefficient is uniquely determined whenever the quantity in
\eqref{non-ded} is nonzero. If $|\mu_j|\neq|\beta ^\alpha|$, it is uniformly
bounded away from zero for all Fourier modes.

Note that the Diophantine bound the non-zero denominator \eqref{non-ded} is strongly related with the condition in \cite[Lemma 3.1]{GimenoJNO21}.

\begin{rmk}[A case computable to all orders] \label{rmk.center-hyper}
For a center--hyperbolic splitting with $|\lambda_i|=1$ in the selected block
and $|\mu_j|\neq1$ in the complementary block, one has
$|\beta ^\alpha|=1$ for every $\alpha$. Hence \eqref{non-ded} is uniformly
nonresonant at every Taylor order. If $\poincmap$ is sufficiently smooth, the
graph coefficients can therefore be continued to arbitrary formal power order.
\end{rmk}

\subsection{Series Formulation}
The same scalar equation also has a useful operator representation. This is a
standard way to obtain regularity estimates for cohomological equations in KAM
theory; see \cite{Llave01,LlaveGJV05}. We keep to the unequal-modulus case,
which is the one used in the numerical examples.

\begin{pro} \label{pro.series-solution}
Let $\rho\in\R^d$ and $\mu,\beta\in\C$, and let $\mathcal X$ be a Banach space of
periodic functions on $\T^d$ such that
\[
T_{\pm\rho}[x](\theta)\bydef x(\theta\pm\rho)
\]
are isometries on $\mathcal X$. Consider \eqref{eq.cohom-graph}. Then
\begin{enumerate}
\item If $|\beta|<|\mu|$, then for every $z\in\mathcal X$ there exists a unique
solution $w\in\mathcal X$ given by
\begin{equation}\label{eq.forward-series-solution}
w(\theta)=\sum_{j\ge0}\mu^{-j-1}\beta^j z(\theta+j\rho).
\end{equation}

\item If $|\mu|<|\beta|$, then for every $z\in\mathcal X$ there exists a unique
solution $w\in\mathcal X$ given by
\begin{equation}\label{eq.backward-series-solution}
w(\theta)=-\sum_{j\ge1}\mu^{j-1}\beta^{-j} z(\theta-j\rho).
\end{equation}
\end{enumerate}
In both cases the series converges absolutely in $\mathcal X$, and
\begin{equation}\label{eq.estimate}
\|w\|_{\mathcal X}\le
\frac{1}{\bigl||\beta|-|\mu|\bigr|}\,\|z\|_{\mathcal X}.
\end{equation}
\end{pro}

\begin{proof}
Equation~\eqref{eq.cohom-graph} is
$(\mu\idop-\beta T_\rho)w=z$. If $|\beta|<|\mu|$, then
$
w-\mu^{-1}\beta T_\rho w=\mu^{-1}z
$, 
and
\[
\|\mu^{-1}\beta T_\rho\|_{\mathcal X\to\mathcal X}
=\frac{|\beta|}{|\mu|}<1.
\]
The Neumann series gives
\[
w
=\mu^{-1}(\idop-\mu^{-1}\beta T_\rho)^{-1}z
=\sum_{j\ge0}\mu^{-j-1}\beta^jT_\rho^jz,
\]
which is \eqref{eq.forward-series-solution}, together with the first version
of \eqref{eq.estimate}. If $|\mu|<|\beta|$, rewrite the equation as
\[
w=\beta^{-1}\mu T_\rho^{-1}w-\beta^{-1}T_\rho^{-1}z.
\]
Since $\|\beta^{-1}\mu T_\rho^{-1}\|=\|\beta^{-1}\mu T_{-\rho}\|=|\mu|/|\beta|<1$, a second Neumann
series gives \eqref{eq.backward-series-solution} and the other version of
\eqref{eq.estimate}. The same estimates applied to the homogeneous equation
give uniqueness.
\end{proof}

\subsection{Matrix Formulation}

The diagonal formulation is enough for the implementation used in this paper,
but it is useful to record what changes when a Floquet block is not diagonal.
At Taylor degree $l$, let $M\in\C^{m\times m}$ denote the complementary
Floquet block. Let $S_l\in\C^{N_l\times N_l}$ be the matrix induced by the
selected block on homogeneous polynomials of degree $l$; equivalently,
$S_l$ represents the substitution $\st\mapsto\Lambda_L\st$ in the chosen monomial
basis. If $\Lambda_L$ is diagonal, the eigenvalues of $S_l$ are precisely the
products $\beta ^\alpha$ with $|\alpha|=l$.

For matrix-valued coefficient functions
$W,Z\colon\T^d\to\C^{m\times N_l}$, the block cohomological equation has the
form
\begin{equation}\label{eq.matrix-cohom-graph}
    M W(\theta)-W(\theta+\rho)S_l=Z(\theta).
\end{equation}
Writing Fourier series and comparing modes gives
\begin{equation}\label{eq.sylvester-mode}
  M\widehat W_k
  -e^{2\pi\I\innerp{k}{\rho}}\widehat W_k S_l
  =\widehat Z_k,
  \qquad k\in\Z^d.
\end{equation}
Thus every Fourier coefficient is obtained from a Sylvester equation. Define
\[
  \mathcal S_k(Y)
  \bydef
  M Y-e^{2\pi\I\innerp{k}{\rho}}Y S_l.
\]
The standard Sylvester criterion \cite{Rosenblum1956,GolubV2013} yields
\begin{equation}\label{eq.spectrum-condition}
  \mathcal S_k \text{ is invertible}
  \quad\Longleftrightarrow\quad
  \Spec(M)\cap
  e^{2\pi\I\innerp{k}{\rho}}\Spec(S_l)=\varnothing.
\end{equation}
Equivalently,
\[
  a\neq e^{2\pi\I\innerp{k}{\rho}}b
  \qquad
  \text{for every }a\in\Spec(M),\ b\in\Spec(S_l).
\]
When both blocks are diagonal this reduces exactly to the scalar denominator
in \eqref{non-ded}.

If \eqref{eq.spectrum-condition} fails for a mode $k$, solvability requires
$\widehat Z_k\in\Range(\mathcal S_k)$; if the kernel is nontrivial, the
solution is not unique. This is the natural formulation for repeated spectral
blocks and Jordan structure. 

\begin{pro}[Formal Sylvester solution]\label{pro.formal-solution}
Assume \eqref{eq.spectrum-condition} for every $k$ in the Fourier
representation. Then the formal solution of \eqref{eq.matrix-cohom-graph} is
\[
W(\theta)=
\sum_{k\in\Z^d}
\mathcal S_k^{-1}(\widehat Z_k)
 e^{2\pi\I\innerp{k}{\theta}}.
\]
Its regularity depends on quantitative estimates for
$\|\mathcal S_k^{-1}\|$ as $|k|\to\infty$.
\end{pro}

\begin{rmk}[Numerical computation] \label{rmk.sylvester}
Sylvester equations like \eqref{eq.sylvester-mode} are normally solved directly by standard numerical linear algebra. Vectorization is useful for displaying the operator. By the identity $\vecop(AXB)=(B^\top\otimes A)\vecop(X)$ |$\otimes$ being the Kronecker product|, equation
\eqref{eq.sylvester-mode} becomes
\[
\left(
\idmat_{N_l}\otimes M
-e^{2\pi\I\innerp{k}{\rho}}S_l^\top\otimes\idmat_m
\right)
\vecop(\widehat W_k)
=\vecop(\widehat Z_k).
\]
The spectral condition is again \eqref{eq.spectrum-condition}. See
\cite{BartelsW1972,LancasterT1985,GolubV2013} for standard Sylvester solvers.
\end{rmk}

\section{Hermitian Structure and Real Formulation}
\label{sec.hermitian}

The graph construction of Section~\ref{sec.graph} is most transparent in spectral
coordinates.  Since the Floquet matrix $B$ is real, nonreal Floquet
multipliers occur in complex-conjugate pairs and it is natural to represent
the corresponding invariant subspaces by conjugate eigenvectors.  This leads
to complex Taylor symbols even though the original map $\poincmap$, the
invariant torus, and the Floquet transformation are real.  A direct
implementation of the complex graph equations could then appear to require
evaluating $\poincmap _\sigma$ at complex phase-space points.  Such evaluations are
unnecessary and, for a merely smooth real map, need not even be defined.

The purpose of this section is to make precise the real structure hidden by
the complex spectral coordinates.  We first describe a graded
realification/complexification isomorphism for multivariate polynomials.  We
then fix a canonical ordering of complex-conjugate Floquet directions and
introduce the associated pair-swap involutions.  Finally, we prove that the
graph $\graphfun$ and the induced reduced dynamics $\intdyn _\sigma$ inherit this real structure.  Consequently, the complete graph computation can be carried out
using real Taylor variables and real evaluations of the original map while
retaining the diagonal or block-diagonal spectral information supplied by the
complex Floquet basis. Algorithm~\ref{palg.graph-tor} incorporates the additionally  \textsuperscript{bis} steps.

In this section ``Hermitian'' refers to invariance under conjugation combined with the exchange of paired variables. It is the polynomial analogue of conjugate-symmetric structures used in Hermitian analysis \cite{Angelo2013}.

\subsection{Realification and Complexification of Multivariate Polynomials}
\label{sec.hermitian.real-complex}

For polynomials whose coefficients are complex and satisfy an additional condition there is a transformation that naturally encodes the same information with only real coefficients. The transformation acts individually in homogeneous polynomials.

We begin at the level of scalar polynomials and with two variables $(\st _1, \st _2)$. 
Let $\C[\st_1, \st_2]$ be the graded ring of polynomials with complex coefficients. Define the map:
\begin{equation} \label{eq.involu}
 p ^\star(\st _1, \st _2) \bydef \overline{p (\overline{\st _2}, \overline{\st _1})}.
\end{equation}
Clearly, $p ^\star$ is an involution, meaning that $p ^\star \circ p ^\star = \idop _{\C[\st _1, \st _2]}$, and coefficientwise \eqref{eq.involu} is read as
\[
p (\st _1, \st _2) = \sum _{|\alpha _1 + \alpha _2| \geq 0} p _{\alpha _1, \alpha _2} \st _1 ^{\alpha _1} \st_2 ^{\alpha _2} \longmapsto p ^\star(\st _1, \st _2) = \sum _{|\alpha _1 + \alpha _2| \geq 0} \overline{p _{\alpha _2,\alpha _1 }} \st_2 ^{\alpha _2} \st _1 ^{\alpha _1} .
\]
Then we consider a subset of $\C[\st _1, \st _2]$ consisting in those that are fixed by the involution transformation:
\begin{equation} \label{eq.calH}
\mathcal{H} \bydef \{ p \in \C [\st _1, \st _2] \colon p ^\star = p \}
\end{equation}
The set $\mathcal{H}$ is a graded real algebra, meaning it is closed under $+$, $\cdot$, and real scalar operations, and it can be seen as an homogeneous-wise condition. By definition,
\begin{equation} \label{eq.Hcarac}
p \in \mathcal{H} \Longleftrightarrow p _{\alpha _1, \alpha _2} = \overline{p _{\alpha _2,\alpha _1 }} \qquad \text{ for all } |\alpha _1 + \alpha _2| \geq 0.
\end{equation}
In particular, if $p \in \mathcal{H}$ whenever $|\alpha _1| = |\alpha _2|$, necessarily $p _{\alpha _1, \alpha _2} \in \R$.

\begin{dfn}[Hermitian polynomial]
 A complex-coefficient polynomial $p$ is called Hermitian with respect to the conjugate-pair structure if $p ^\star = p$. The set of such polynomials is denoted by $\mathcal{H}$.
\end{dfn}

\smallskip

Consider the map described by the complex generators $(\st _1, \st _2)$ and the real ones $(a, b)$,
\begin{align}
\mathtt{R}\colon  \C [\st _1, \st _2] &\to \C [a,b], \qquad \st _1 \mapsto \frac{a + \I b}{2}, \quad \st _2 \mapsto  \frac{a - \I b}{2} , \qquad 
\text{equivalently } \mathtt{R} = \frac{1}{2}
\begin{pmatrix}
 1 & \I \\ 1 & -\I
\end{pmatrix} \label{eq.R}
\intertext{and its inverse}
\mathtt{C} \colon \C [a,b] &\to \C [\st _1, \st _2], \qquad a \mapsto \st _1 + \st _2, \quad b \mapsto  -\I(\st _1- \st _2)\qquad 
\text{equivalently } \mathtt{C} = 
\begin{pmatrix}
 1 & 1 \\ -\I & \I
\end{pmatrix}. \nonumber
\end{align}

The factor $\frac{1}{2}$ is a normalization choice. In our case, it will match numerical conventions in  Section~\ref{sec.krawtchouk}.

\begin{pro}
 $\mathtt R$ is a graded $\C$-algebra isomorphism and $\mathtt C=\mathtt R^{-1}$.
\end{pro}
\begin{proof}
 Since $\C[\st_1,\st_2]$ is a polynomial ring, specifying the images of the generators
$\st_1,\st_2$ defines a unique $\C$-algebra homomorphism $\mathtt R$; likewise $\mathtt C$
is uniquely determined by the images of $a,b$. The compositions give the identity operators are straightforward.

Finally, $\mathtt R$ and $\mathtt C$ are graded because $\st_1,\st_2$ are sent to
linear (degree $1$) polynomials in $a,b$, and $a,b$ are sent to linear polynomials
in $\st_1,\st_2$, so homogeneous degree is preserved.
\end{proof}

The following lemma is crucial to establish the notion of realification and complexification.
\begin{lem} \label{lem.compat}
For every $p\in\C[\st_1,\st_2]$,
$\overline{\mathtt R(p)}=\mathtt R(p^\star)
$.
\end{lem}
\begin{proof}
Evaluating at $(a,b)$ and using $s_1=\frac{1}{2}(a+\I b)$, $\st_2=\frac{1}{2}(a-\I b)$ gives
\[
\overline{\mathtt R(p)}(a,b)
=\overline{\,p\biggl(\frac{a+\I b}{2},\frac{a-\I b}{2}\biggr)\,}
=\overline{\,p\biggl(\overline{\frac{a-\I b}{2}},\overline{\frac{a+\I b}{2}}\biggr)\,}
=\overline{\,p(\overline{\st_2},\overline{\st_1})\,}\Big|_{\st_1=\frac{1}{2}(a+\I b),\ \st_2=\frac{1}{2}(a-\I b)}
=\mathtt R(p^\star)(a,b). \qedhere
\]
\end{proof}

\begin{pro}\label{pro.iso}
The restriction of $\mathtt R$ on $\mathcal{H}$ gives a graded real algebra isomorphism:
\begin{equation} \label{eq.frakRC}
\begin{aligned} 
\mathfrak{R} \colon \mathcal{H} &\to \R[a,b] & \mathfrak{C} \colon \R[a,b] &\to \mathcal{H} \\
p(\st _1, \st_2)& \mapsto p\biggl(\frac{a + \I b}{2}, \frac{a - \I b}{2}\biggr) &
q(a, b)& \mapsto q\bigl(\st _1 + \st _2, -\I(\st _1-\st _2)\bigr) 
\end{aligned}
\end{equation}
meaning that $\mathfrak{R} \circ \mathfrak{C} = \idop_{\R[a,b]}$ and $\mathfrak{C} \circ \mathfrak{R} = \idop_{\mathcal{H}}$.
\end{pro}

\begin{proof} \
\begin{itemize}
    \item $\mathfrak{R}(\mathcal{H}) \subset \R [a,b]$. Indeed, $p^\star=p$, by Lemma~\ref{lem.compat},
$\overline{\mathtt R(p)}=\mathtt R(p^\star)=\mathtt R(p)$ and $\mathfrak{R}(p)$ has real coefficients. 

    \item $\mathfrak{C}(\R[a,b]) \subset \mathcal{H}$. Let $q\in\R[a,b]$ and $p\bydef \mathtt C(q)\in\C[\st_1,\st_2]$. Then $\mathtt R(p)=\mathtt R(\mathtt C(q))=q$. By Lemma~\ref{lem.compat},
\[
\mathtt R(p^\star)=\overline{\mathtt R(p)}=\overline{q}=q=\mathtt R(p).
\]
Since $\mathtt R$ is injective, $p^\star=p$.
Therefore $q=\mathtt R(p)\in \mathtt R(\mathcal H)$.

    \item $\mathfrak R$ and $\mathfrak C$ are graded because $\mathtt R$ and $\mathtt C$
send degree $1$ generators to degree $1$ polynomials, hence preserve total degree. \qedhere
\end{itemize}
\end{proof}

Lemma~\ref{lem.compat} serves to prove that the $\C$-algebra isomorphism $\mathtt R$ restricts to an $\R$-algebra isomorphism $\mathfrak{R}$ giving commutative diagrams of ring of homomorphisms.
\[
\begin{tikzcd}
p \arrow[d,mapsto] &\C[\st_1,\st_2] 
  \arrow[r, "\mathtt R"] \arrow[d,"(\cdot)^\star"'] &
\C[a,b]  \arrow[d,"\overline{(\cdot)}"] & \mathtt R(p) \arrow[d,mapsto] \\
p ^\star & \C[\st_1,\st_2] 
  \arrow[r, "\mathtt R"] &
\C[a,b]  & \overline{\mathtt R(p)}
\end{tikzcd} \quad \text{ yield } \quad 
\begin{tikzcd}
\C[\st_1,\st_2]
  \arrow[r, shift left=.6ex, "\mathtt R"]
&
\C[a,b]
  \arrow[l, shift left=.6ex, "\mathtt C"]
\\
\mathcal H
  \arrow[u, hook]
  \arrow[r, shift left=.6ex, "\mathfrak R"]
&
\R[a,b]
  \arrow[u, hook]
  \arrow[l, shift left=.6ex, "\mathfrak C"]
\end{tikzcd}
\]

\begin{rmk}[Multivariate extension] \label{rmk.multi-extend}
 The same construction extends componentwise to $q$ pairs of variables, namely $\st _j = (\st _{j, 1}, \dotsc, \st _{j,q})$ for $j \in \{1,2\}$, $a = (a _1, \dotsc, a _q)$ and $b = (b _1, \dotsc, b _{q})$. 
 The construction also remains valid after adjoining additional variables that are
left fixed by the involution (e.g. polynomials in $\C[\st_1,\st_2,\st_3]$ with $\st _3^\star = \st _3$).
See Section~\ref{sec.canonical-conjugate-pairs} for details.
\end{rmk}

\begin{rmk}[Coordinate extension] \label{rmk.coord-extend}
 All the arguments apply when the coefficients have multiple coordinates, that is, for instance to $\R ^n[a,b]$ and $\C ^n[\st _1, \st _2]$. In that setting, the maps
$\mathtt R$, $\mathtt C$, $\mathfrak R$, and $\mathfrak C$ are applied to
each coordinate separately. For readability, we have written the scalar-valued version.
\end{rmk}

\subsection{Canonical Convention for Conjugate Pairs}
\label{sec.canonical-conjugate-pairs}

We now connect the polynomial construction with the Floquet coordinates of Section~\ref{sec.graph}. Since $B$ is real, every nonreal eigenvalue occurs together with
its complex conjugate.  To obtain a real invariant manifold, the selected
spectral subspace $E_L$ and its complete complement $E_N$ are required to be
closed under complex conjugation.  In particular, selecting one member of a
nonreal eigenpair while excluding the other would not define a real spectral
subspace and is not allowed in the real formulation yielding to immediate resonance conditions when solving the cohomological equation.

\smallskip

Let $q_L$ and $r_L$ denote, respectively, the number of nonreal conjugate
pairs and the number of real directions in $E_L$, so that
$$
    n_L=2q_L+r_L.
$$
Likewise let $q_N$ and $r_N$ describe $E_N$, with $    n_N=2q_N+r_N$.
Within each block of $V = (v _j)$ and $\Lambda = (\lambda _j)$ in \eqref{eq.GT.spectral-blocks} we adopt the following canonical ordering for $j = 1, \dotsc, q _L + q _N$.

\begin{dfn}[Canonical ordering]\label{dfn:caor}
Let us consider the following assumptions for a given ordering
\begin{enumerate}
\renewcommand{\theenumi}{\bf A\arabic{enumi}}
\renewcommand{\labelenumi}{{\theenumi})}
 \item \label{A1} All the non-real pair are placed before the real directions;
 
 \item \label{A2} Each non-real conjugated pair is stored in adjacent columns as $(v _j, \overline{v _j})$, with the corresponding Floquet multipliers ordered as and $(\lambda _j, \overline{\lambda _j})$;
 
 \item Real eigenvectors are represented by real columns and are
    therefore fixed under complex conjugation;
 
 \item \label{A4} The Taylor symbols associated with each non-real conjugate pair in $\Et$ are ordered accordingly as $(\st_j,\overline{\st_j})$, for $j=1,\dotsc,q_L$.
 
    Despite the conjugate notation, $s_j$ and $\overline{s_j}$ are treated
    as independent formal Taylor symbols.
\end{enumerate}
When \ref{A1}--\ref{A4} are satisfied, we say to be in canonical ordering form. 
\end{dfn}

A different external ordering may be used, but it should be converted once to this canonical form before applying realification operators or modifying the operators accordingly. Notice \ref{A4} avoids pair-dependent permutation during jet arithmetic. 

\smallskip

For notational purposes, let us introduce the pair-swap (involutived) matrices for canonical forms:
\begin{equation} \label{eq.Jmat}
 J _2 \bydef 
 \begin{pmatrix}
  0 & 1 \\ 1 & 0
 \end{pmatrix}, \quad 
 J _L \bydef \diag(\underbrace{J _2,\dotsc,J _2}_{q _L\text{ copies}}, \idmat _{r _L}), \quad 
 J _N \bydef \diag(\underbrace{J _2,\dotsc,J _2}_{q _N\text{ copies}}, \idmat _{r _N}), \quad 
 J \bydef
 \begin{bmatrix}
  J _L & J _N
 \end{bmatrix}.
\end{equation}
Thus, \eqref{eq.involu} can in general be written as $p^\star(\st) = \overline{p (J _L \overline{\st})}$.
The notation \eqref{eq.Jmat} shows that the canonical form is also preserved in directed skew-products because $\overline{V} = V J$, $\overline{\Lambda} = J \Lambda J$, and for directed matrices $\overline{\Lambda _\sigma} = J \Lambda _\sigma J$. 

The matrices \eqref{eq.Jmat} generalizes \eqref{eq.involu} for maps:
\begin{dfn} \label{dfn.Hermitian.fun}
 Given a function $W\colon \C ^{\nt} \times \T ^d \to \C ^{\nh}$. Define:
 \begin{equation*}
  W ^\star (\st, \theta) \bydef J _N \overline{W(J _L \overline{\st}, \theta)}.
 \end{equation*}
 $W$ is then called Hermitian if, and only if, $W^\star = W$. That is, when $W  (J _L \overline{\st}, \theta) = J _N \overline{W(\st, \theta)}$.
\end{dfn}

\subsection{Hermitian Symmetry of the Graph} \label{sec.graph.herm}

We want to show that the graph equation \eqref{eq.GT.K-graph} preserves the real structure. In other words, it satisfies Definition~\ref{dfn.Hermitian.fun}.

\smallskip

Following the canonical conventions in Section~\ref{sec.canonical-conjugate-pairs}. Introduce the operators

\begin{equation}\label{eq.RC-op}
 \mathtt{R} _L \bydef \diag(\underbrace{\mathtt{R},\dotsc,\mathtt{R} }_{q _L\text{ copies}}, \idmat _{r _L}), \quad 
 \mathtt{R} _N \bydef \diag(\underbrace{\mathtt{R},\dotsc,\mathtt{R} }_{q _N\text{ copies}}, \idmat _{r _N}), \qquad 
 \mathtt{C} _L \bydef \mathtt{R} _L^{-1}, \quad 
 \mathtt{C} _N \bydef \mathtt{R} _N^{-1}.
\end{equation}
And consider $V _{L,\R} \bydef \Vt \mathtt{R} _L \in \R ^{n \times \nt}$ and $V _{N,\R} \bydef \Vh \mathtt{R} _N \in \R ^{n \times \nh}$.  Since the Floquet transformation $C(\theta)$ is real, the corresponding real adapted frame is just:
\begin{equation}
 Q _{L,\R} (\theta) \bydef C(\theta)V _{L,\R} \qquad 
 \text{and} \qquad
 Q _{N,\R} (\theta) \bydef C(\theta)V _{N,\R} .
\end{equation}

\begin{rmk}[Compatibility with spectral selection]
\label{rem:HR-spectral-selection}
The complete-complement principle of Section~\ref{sec.graph} is compatible with the real
structure only when both $\Et$ and $\Eh$ are conjugation invariant.  Thus, if
$\Et$ contains a complex eigendirection, it contains its conjugate as well;
all unselected conjugate pairs remain together in $\Eh$.  This condition is independent of the time orientation $\sigma$.
\end{rmk}

\smallskip

Theorem~\ref{thm:HR-graph-hermitian} shows that $\graphfun$ in \eqref{eq.GT.graph} subject to \eqref{eq.GT.normalization} maps the real selected slice into the real complementary
slice. The induced dynamics $\intdyn_\sigma$ maps the real selected slice into itself. Thus, both are Hermitian in the sense of Definition~\ref{dfn.Hermitian.fun}.

\begin{thm}[Hermitian symmetry of the graph]
\label{thm:HR-graph-hermitian}
Assume that the spectral selections $\Et$ and $\Eh$ are closed under complex conjugation and canonically ordered according to Definition~\ref{dfn:caor} of Section~\ref{sec.canonical-conjugate-pairs}.
Assume that the cohomological equations occurring in Algorithm~\ref{palg.graph-tor} are uniquely solvable at every computed order.  Then the computed $\graphfun$ and $\intdyn _\sigma$ are Hermitian, that is,
\begin{equation*}
\label{eq.HR-W-f-Hermitian}
    \graphfun(J_L\overline{\st},\theta)
    =J_N\overline{\graphfun(\st,\theta)}, \qquad \text{and} \qquad
    \intdyn _\sigma(J_L\overline{\st},\theta)
    =J_L\overline{\intdyn _\sigma(\st,\theta)}.
\end{equation*}
In particular, coefficients of $\graphfun$ in \eqref{eq.GT.W-FPS} satisfy for all order $l$ and all multiindex $\alpha \in \N ^{\nt}$ such that $|\alpha | = l$,
\begin{equation*} 
w _{l, J _L \alpha}(\theta) = J _N \overline{w _{l,\alpha}(\theta)}.
\end{equation*}                                                                                                            
\end{thm}

Before proving the theorem, let us show the following lemma:
\begin{lem}\label{lem.cR.Hermitian}
 Let $\cR _\sigma$ be the directed adapted map from \eqref{eq.GT.reduced-map} and let $J$ in \eqref{eq.Jmat}. Then $\cR _\sigma$ is Hermitian, i.e. 
 \begin{equation*}
  \cR _\sigma (J _L \overline{\st}, J _N \overline{\sh}, \theta) = J\overline{\cR _\sigma (\st, \sh, \theta)}.
 \end{equation*}
\end{lem}
\begin{proof}
 Since $\poincmap _\sigma$, $\tor$, and $C$ are real and $\overline{V} = V J$ from \eqref{eq.Jmat}, then its Taylor expansion satisfy order-by-order the conditions for $\cR _\sigma$.
\end{proof}

Lemma~\ref{lem.cR.Hermitian} infers different conclusions. If $\poincmap _\sigma$ admits a real-analytic complexification, the equalities in the lemma are identities of analytic maps in a complex neighborhood of the torus.  
If $\poincmap _\sigma$ is only smooth, then the equality in the lemma is understood as identify of formal Taylor jets obtained from the real derivatives of $\poincmap _\sigma$.  This last formal interpretation is sufficient for the finite-order jet transport used by the Algorithm~\ref{palg.graph-tor} and does not require evaluating the original map at a genuinely complex phase-space point.

\begin{proof}[Proof of Theorem~\ref{thm:HR-graph-hermitian}]
 Define $\widetilde \graphfun \bydef \graphfun ^\star$ and $\widetilde \intdyn _\sigma \bydef \intdyn ^\star_\sigma$. Apply complex conjugation to \eqref{eq.GT.graph-invariance-abstract}, substitute $\st \mapsto J _L \overline{\st}$ and use Lemma~\ref{lem.cR.Hermitian}. Since $\rho _\sigma$ is real, 
 \begin{equation*}
  \cR _{N,\sigma}  \bigl(\st, \widetilde \graphfun (\st, \theta), \theta \bigr) = \widetilde \graphfun \bigl( \widetilde \intdyn _\sigma(\st, \theta), \theta + \rho _\sigma \bigr),
 \end{equation*}
 with $\widetilde \intdyn _\sigma(\st, \theta) = \cR _{L,\sigma}\bigl(\st, \widetilde \graphfun (\st, \theta), \theta \bigr)$. 
 Thus $(\widetilde \graphfun, \widetilde \intdyn _\sigma)$ satisfies the same directed invariance equation as $(\graphfun, \intdyn _\sigma)$ with same normalization conditions and the linear part of $\widetilde \intdyn _\sigma $ is also $\Lambda _{L,\sigma}$. So, order-by-order the uniqueness of cohomological equations implies $\widetilde \graphfun = \graphfun$ and $\widetilde \intdyn _\sigma = \intdyn _\sigma$.
\end{proof}

\subsubsection{Real Graph and Real Reduced Dynamics}
Theorem~\ref{thm:HR-graph-hermitian} yields an explicit real formulation. Let $u \in \R ^{\nt}$ and set $\st = \mathtt{R}_L u$. Use \eqref{eq.RC-op} and define:
\begin{equation}
 \graphfun _\R (u, \theta) \bydef \mathtt{C} _N \graphfun (\mathtt{R} _L u, \theta) \qquad \text{and} \qquad 
 \intdyn _{\sigma, \R} (u, \theta) \bydef \mathtt{C} _L \intdyn _\sigma (\mathtt{R} _L u, \theta).
\end{equation}

\begin{cor}[Equivalent real graph formulation]
\label{cor:HR-real-graph}
Under the hypotheses of Theorem~\ref{thm:HR-graph-hermitian}, the functions
$\graphfun _\R$ and $\intdyn _{\sigma,\R}$ have real Taylor coefficients and take real values
for real $(u,\theta)$.  Conversely,
\begin{equation*}
\label{eq.HR-complex-from-real-graph}
    \graphfun(\st,\theta)
    =
    \mathtt{R}_N W_\R(\mathtt{C}_L \st,\theta) \qquad \text{and} \qquad 
    \intdyn_\sigma(\st,\theta)
    =
    \mathtt{R}_L \intdyn _{\sigma,\R}(\mathtt{C}_L \st,\theta).
\end{equation*}
\end{cor}
Hence the complex Hermitian formulation and the real formulation contain exactly the same polynomial information. Similarly, a parameterization formulation similar to \eqref{eq.GT.K-graph} can be derived and no holomorphic extension assumption is required for jet transport.

\subsection{Practical Implementation of the Real Formulation}
\label{sec.num-real}
Section~\ref{sec.graph.herm} shows that the complex Hermitian formulation and the
real formulation contain the same Taylor information; see
Theorem~\ref{thm:HR-graph-hermitian} and
Corollary~\ref{cor:HR-real-graph}. The purpose of the real formulation is
therefore computational rather than representational: the complex Floquet
basis is retained where it is advantageous, in particular in the
cohomological equations, whereas the nonlinear map $\poincmap _\sigma$ and its jets
are evaluated only at real phase-space points.

\smallskip

For real reduced variables $u\in\R^{\nt}$, the identities
of Corollary~\ref{cor:HR-real-graph} allow the parameterization
\eqref{eq.GT.K-graph} to be written entirely in the real adapted frame as
\begin{equation}
\label{eq.real-param-implementation}
    \param_{\R}(u,\theta)
    =
    \tor(\theta)
    + Q_{L,\R}(\theta)u
    + Q_{N,\R}(\theta)\graphfun_{\R}(u,\theta).
\end{equation}
Consequently, $\param_{\R}(u,\theta)$ is real for real $(u,\theta)$, and the
jet-transport step of Algorithm~\ref{palg.graph-tor} can be performed using
real Taylor variables and real evaluations of $\poincmap _\sigma$.

The complex representation is needed only for the linear-algebraic part of
the computation. In particular, before solving the cohomological equations
in Step~\ref{palg.graph-tor.cohom} in Algorithm~\ref{palg.graph-tor}, the relevant homogeneous coefficient
blocks are transformed from real Taylor variables to the canonical complex
conjugate-pair variables of Section~\ref{sec.canonical-conjugate-pairs}. After the
cohomological equations have been solved in the complex Floquet coordinates,
the resulting coefficients are transformed back to the real representation
before the next jet-transport step. Since these changes of representation
act only on the Taylor variables, they commute with the dependence on
$\theta$, with the rotation $\theta\mapsto\theta+\rho_\sigma$, and with
the Fourier transform in the base variable.

Under the canonical ordering, the change of
Taylor variables factorizes into independent transformations, one for each
non-real conjugate pair. The real Floquet directions are left unchanged.
This pairwise structure is what makes Krawtchouk matrices useful in the
implementation: the full multivariate transformation never needs to be
formed as a dense tensor-product matrix.

\subsubsection{Krawtchouk Transformations} \label{sec.krawtchouk}

We first consider a single non-real conjugate pair. Let $(a,b)$ denote the
corresponding real Taylor variables and $(s,\overline{s})$ the associated
formal complex variables. Recall from \eqref{eq.R} the transformations $\mathtt{R}$ and $\mathtt{C}$.

Consider a homogeneous polynomial of local degree $k$,
\[
    q(a,b)=\sum_{q=0}^{k} X[q]\,a^{k-q}b^q.
\qquad
\text{after complexification,}
\quad
    \mathtt{C}(q)(s,\overline{s})
    =
    \sum_{p=0}^{k}Y[p]\,s^{k-p}\overline{s}^{\,p}.
\]
The coefficient of $s^{k-p}\overline{s}^{\,p}$ in $(s+\overline{s})^{k-q}(s-\overline{s})^q$
is a Krawtchouk coefficient. This gives the following definition originally introduced in \cite{Krawtchouk1929}, see also \cite[pg. 115]{FeinsilverK2005}.

\begin{dfn}
\label{dfn.Krawtchouk}
For $k\geq0$, the Krawtchouk matrix
$\kop\in\Z^{(k+1)\times(k+1)}$ is defined by
\begin{equation}
\label{eq.kraw-definition}
    \kop[p,q]
    =
    \sum_{j=0}^{\min(p,q)}
    (-1)^j
    \binom{q}{j}
    \binom{k-q}{p-j},
    \qquad
    p,q\in\{0,\dotsc,k\}.
\end{equation}
When the dependence on the local degree needs to be emphasized, we write
$\kop=\kop^{(k)}$.
\end{dfn}

Equivalently, for every fixed $q$,
\begin{equation*}
\label{eq.kraw-generating}
    \sum_{p=0}^{k}\kop[p,q]\,z^p
    =
    (1+z)^{k-q}(1-z)^q.
\end{equation*}
Thus $\kop[p,q]$ is precisely the coefficient produced by the binomial
expansion associated with the $q$-th power of the difference
$s-\overline{s}$. In particular,
\begin{equation}
\label{eq.kraw-involution}
    \kop^2=2^k\idmat,
\end{equation}
so that the same Krawtchouk kernel appears in both directions of the
coefficient transformation, up to the diagonal phase and normalization
factors.

Direct evaluation of \eqref{eq.kraw-definition} involves alternating sums
and may suffer substantial cancellation in floating-point arithmetic.
Moreover, repeatedly evaluating binomial coefficients is unnecessary.
A convenient exact construction uses the recurrence
\begin{equation}
\label{eq.kraw-recu}
\begin{aligned}
    \kop[p,0] &= \binom{k}{p},
    \qquad
    \kop[0,q] = 1,
    \\
    \kop[p,q]
    &=
    \kop[p,q-1]
    -\kop[p-1,q-1]
    -\kop[p-1,q],
    &&
    1\leq p,q\leq k.
\end{aligned}
\end{equation}
Hence one may precompute the family of Krawtchouk matrices required up to a
prescribed Taylor order. For a truncation at total degree $N$, storing all
matrices of degrees $0,\dotsc,N$ requires
$
    \sum_{k=0}^{N}(k+1)^2
    =
    O(N^3)
$ 
integer entries. Alternatively, the matrices can be generated on demand. The following lemma simplifies the storage:

\begin{lem} \label{lem.kop-property}
 For every local degree $k\geq 0$, $\kop[k-p,q] = (-1) ^q \kop[p,q]$ for $p,q\in\{0,\dotsc,k\}$.
\end{lem}
\begin{proof}
 Consider $\widetilde \kop [p,q] \bydef (-1)^q \kop [k-p,q]$ and let us prove it satisfies \eqref{eq.kraw-recu}.
 \begin{equation*}
     \widetilde\kop [p,0] = \kop[k-p,0] = 
     \binom{k}{k-p} = \binom{k}{p} = \kop [p,0], \quad \text{ and } \quad
      \widetilde\kop[0,q] = (-1) ^q \kop [k,q] = 1
 \end{equation*}
 where we have used $\kop[k,q] = (-1)^{q}$ treating individual cases. Finally, 
 \begin{multline*}
     \widetilde\kop [p,q-1] - \widetilde\kop [p-1, q-1] - \widetilde\kop [p-1, q] = \\(-1) ^{q}  
     \Bigl(- \kop [k-p,q-1] + \kop [k - p+1,q-1] - \kop[k-p+1,q] \Bigr)
 \end{multline*}
 Using \eqref{eq.kraw-recu} at $[k-p+1,q]$ and rearranging yield \eqref{eq.kraw-recu} with $\widetilde \kop$.
\end{proof}

Lemma~\ref{lem.kop-property} implies that only the rows
$0\leq p\leq\lfloor k/2\rfloor$ need to be computed explicitly from
\eqref{eq.kraw-recu}; the remaining rows follow by symmetry. This reduces
both storage and arithmetic by approximately a factor of two, see Figure~\ref{fig.kraw-mat}. If $k$ is
even, the middle row satisfies $\kop[k/2,q]=0$ for every odd $q$.

\begin{figure}[ht]
\centering
\[
\kop = \left(
\begin{array}{cccc}
\renewcommand{\arraystretch}{1.5}
  1 & 1 & \cdots & 1 \\
  \binom{k}{1} & \kop[1,1] & \cdots & \kop[1,k] \\
  \vdots & \vdots & & \vdots \\
  \binom{k}{m} & \kop[m,1] & \cdots & \kop[m,k]
\\ \hline
  \binom{k}{m+1} & \kop[m+1,1] & \cdots & \kop[m+1,k] \\
  \vdots & \vdots & & \vdots \\
  \binom{k}{k} & \kop[k,1] & \cdots & \kop[k,k]
\end{array} \right)
\quad
\begin{array}{l}
\begin{aligned}
p&=0,\ldots,m, \text{ with } m = \lfloor k/2 \rfloor\\
&\text{ rows computed using~\eqref{eq.kraw-recu}}
\end{aligned}
\\[8ex]
\begin{aligned}
p&=m+1,\ldots,k, \text{from }\\
&\kop[k-p,q]=(-1)^q\kop[p,q]
\end{aligned}
\end{array}
\]
\caption{Compact construction of the Krawtchouk matrix of local degree $k$ using \eqref{eq.kraw-recu} and Lemma~\ref{lem.kop-property}.} \label{fig.kraw-mat}
\end{figure}

\paragraph{Pairwise coefficient transforms.}
Let $X[q]$ be the coefficient of $a^{k-q}b^q$, and let $Y[p]$ be the
coefficient of $s^{k-p}\overline{s}^{\,p}$. The coefficient-level
complexification and realification maps are
\begin{equation}
\label{eq.kraw-pair-transforms}
\begin{aligned}
    Y[p]
    &=
    \sum_{q=0}^{k}
    \kop[p,q](-\I)^q X[q],
    &&\text{real to complex},
    \\
    X[p]
    &=
    2^{-k}\I^p
    \sum_{q=0}^{k}
    \kop[p,q]Y[q],
    &&\text{complex to real}.
\end{aligned}
\end{equation}
These are the coefficient realizations of $\mathtt{C}$ and
$\mathtt{R}$, respectively. They are inverse to one another, consistently
with \eqref{eq.kraw-involution} and with
$\mathtt{C}=\mathtt{R}^{-1}$.

\begin{pro} \label{pro.hermi}
 If $X[q] \in \R$ for all $q =0, \dotsc, k$, then $Y[k-p] = \overline{Y[p]}$ for all $p = 0, \dotsc, k$.
\end{pro}
\begin{proof}
Using Lemma~\ref{lem.kop-property}, $(-1)^q(-\I)^q = \overline{(-\I)^q}$, $\kop[p,q] \in \Z$, and $X[q] \in \R$,
\begin{equation*}
        Y[k-p] 
        = \sum _{q=0}^k \kop[k-p,q](-1)^q X[q] 
        = \sum _{q=0}^k \kop[p,q]\overline{(-\I)^q} X[q] =\overline{Y[p]}.
    \qedhere
\end{equation*}
\end{proof}

\paragraph{Global transformation.}
In canonical ordering of Section~\ref{sec.canonical-conjugate-pairs}.
For the $q_L$ non-real conjugate pairs in the selected space, write the real
Taylor variables pairwise as
$
    (a_1,b_1),\dotsc,(a_{q_L},b_{q_L})
$, 
followed, when $r_L>0$, by the real Taylor variables associated with the
real Floquet directions. Fix local degrees
$k_1,\dotsc,k_{q_L}$. If
$
    X[q_1,\dotsc,q_{q_L};\delta]
$ denotes the coefficient of
\[
    \prod_{j=1}^{q_L}
    a_j^{k_j-q_j}b_j^{q_j}
\]
times a monomial of multiindex $\delta$ in the $r_L$ real directions,
and 
$
    Y[p_1,\dotsc,p_{q_L};\delta]
$ denotes the coefficient of
\[
    \prod_{j=1}^{q_L}
    s_j^{k_j-p_j}\overline{s_j}^{\,p_j}
\]
times the same real-direction monomial, then the global real-to-complex
transformation factorizes as
\begin{equation}
\label{eq.kraw-global}
Y[p_1,\dotsc,p_{q_L};\delta] =
\sum_{q_1=0}^{k_1}\cdots\sum_{q_{q_L}=0}^{k_{q_L}}
\left(
    \prod_{j=1}^{q_L}
    \kop^{(k_j)}[p_j,q_j](-\I)^{q_j}
\right)
X[q_1,\dotsc,q_{q_L};\delta].
\end{equation}
Thus the full multivariate operator is a tensor product of pairwise
Krawtchouk kernels. In practice, however, \eqref{eq.kraw-global} should be
implemented by applying the one-pair kernels successively; explicitly
forming the tensor-product matrix is unnecessary.

For a complex Taylor multiindex $\alpha\in\N^{\nt}$, let
$
    \alpha^\star \bydef J_L\alpha
$,
that is, $\alpha^\star$ exchanges the exponents within every non-real
conjugate pair and leaves the real directions unchanged. Applying
Proposition~\ref{pro.hermi} pair by pair gives the scalar symmetry
\begin{equation}
\label{eq.kraw-global-hermitian}
    Y[\alpha^\star]=\overline{Y[\alpha]}.
\end{equation}
For vector-valued coefficients, the corresponding output pair exchange is
handled by $J_N$. In particular, for the graph coefficients this is
precisely the coefficient symmetry associated with
Theorem~\ref{thm:HR-graph-hermitian}. 

\begin{ex}
For $k\in\{0,1,2\}$, the corresponding Krawtchouk matrices are
\[
    \begin{pmatrix}1\end{pmatrix},
    \qquad
    \begin{pmatrix}
        1 & 1\\
        1 & -1
    \end{pmatrix},
    \qquad
    \begin{pmatrix}
        1 & 1 & 1\\
        2 & 0 & -2\\
        1 & -1 & 1
    \end{pmatrix}.
\]
As a two-pair example, consider the real polynomial $q(a_1,b_1,a_2,b_2)=a_1^2b_2$.

Its total degree is $3$, whereas its local degrees are $2$ in the first
pair and $1$ in the second. Applying the pairwise complexification gives
\begin{equation*}
    a_1^2b_2
    =
    -\I s_1^2s_2
    +\I s_1^2\overline{s_2}
    -2\I s_1\overline{s_1}s_2
    +2\I s_1\overline{s_1}\,\overline{s_2} 
    -\I\overline{s_1}^{\,2}s_2
    +\I\overline{s_1}^{\,2}\overline{s_2}.
\end{equation*}
The coefficients occur in conjugate-paired positions exactly as predicted by
\eqref{eq.kraw-global-hermitian}.
\end{ex}

\paragraph{Placement in Algorithm~\ref{palg.graph-tor}.}
The real implementation is incorporated in the algorithm as an optional step. The graph coefficients are stored in the real representation $W_{\R}$. The
parameterization \eqref{eq.real-param-implementation} is assembled in the
real adapted frame, and the nonlinear jet transport in
Step~\ref{palg.graph-tor.jet} is carried out using real arithmetic. At the
homogeneous degree currently being solved, the coefficient block required by
Step~\ref{palg.graph-tor.cohom} is complexified by the pairwise transforms
\eqref{eq.kraw-pair-transforms}--\eqref{eq.kraw-global}. The
cohomological equations are then solved in the complex Floquet coordinates,
where the selected linear dynamics has the spectral form used in
Section~\ref{sec.solve.cohom}. The newly computed graph coefficients are
finally realified and stored again in $W_{\R}$ before proceeding to the
next Taylor order.

There are two practical consequences. First, the original map
$\poincmap _\sigma$ is never evaluated at a complex phase-space point. 
Second, the diagonal or block-diagonal complex spectral representation is retained exactly where it is useful, namely in the cohomological solver. 
The Krawtchouk transformations therefore act only as interfaces between the real jet-transport representation and the complex spectral representation.

The Hermitian identities also provide inexpensive implementation checks:
after complexification one may verify \eqref{eq.kraw-global-hermitian}, while after realification the recovered coefficients should be real up to the prescribed numerical tolerance, see Test~\ref{test6-conj}.

\section{Local Center Invariant Manifold Reduction}
\label{sec.center-reduction}

The graph construction of Section~\ref{sec.graph} provides, in addition to a
local invariant manifold, the nonlinear dynamics induced on the selected
graph coordinates; see \eqref{eq.GT.induced-dynamics} and
\eqref{eq.GT.f-series}.  For a center spectral selection this induced map is
a low-dimensional nonlinear model of the dynamics in a neighborhood of the
original reducible invariant torus.  The purpose of the present section is
to make this reduction explicit and to describe how it can be used to
construct, detect, and continue nearby invariant tori.  In particular, this
is the form of the graph construction that will be used in the numerical
experiments.

Throughout this section we use the forward orientation $\sigma=+1$ and take
\begin{equation}
    \Et=E^c,
    \qquad
    \Eh=E^s\oplus E^u,
    \qquad
    n_c\bydef\dim E^c=\nt,
    \label{eq.center-selection}
\end{equation}
or, more generally, let $E^c$ denote the selected center spectral block and
let $\Eh$ be its complete spectral complement.  We assume that the
solvability conditions required by the graph construction are satisfied for
this spectral selection, see Lemma~\ref{lem.cohom}.  Algorithm~\ref{palg.graph-tor} then produces a local center graph together with its induced internal dynamics.  We denote the corresponding parameterization and reduced map by $\param^c$ and
$\intdyn_c$, respectively:
\begin{equation*}
    \param^c(\st,\theta)
    =\tor(\theta)
      +\adframe_L(\theta)\st
      +\adframe_N(\theta)\graphfun^c(\st,\theta),
      \qquad
    \intdyn_c(\st,\theta)
    =\cR_{L,+1}\bigl(\st,\graphfun^c(\st,\theta),\theta\bigr).
\label{eq.center-param-reduced-map}
\end{equation*}
The graph normalization gives
$    \graphfun^c(0,\theta)=0$,
$    \dop_{\st}\graphfun^c(0,\theta)=0$,
$    \dop_{\st}\intdyn_c(0,\theta)=\Lambda_c$,
where $\Lambda_c\bydef\Lambda_L$ for the selection
\eqref{eq.center-selection}.  The defining invariance equation is therefore
\begin{equation}
    \poincmap\bigl(\param^c(\st,\theta),\theta\bigr)
    =
    \param^c\bigl(\intdyn_c(\st,\theta),\theta+\rho\bigr).
    \label{eq.center-invariance}
\end{equation}

For real systems, the center block is selected together with its complex
conjugates as required by Remark~\ref{rem:HR-spectral-selection}.  The real
formulation of Section~\ref{sec.num-real} then gives real center coordinates
and a real reduced map.  In the remainder of this section we use this real
formulation and suppress the subscript $\R$ to avoid cumbersome notation:
\begin{equation*}
    \st\in\mathcal D_c\subset\R^{n_c},
    \qquad
    \param^c\colon\mathcal D_c\times\T^d\to \R^n,
    \qquad
    \intdyn_c\colon\mathcal D_c\times\T^d\rightarrow \R^{n_c}.
    \label{eq.center-real-convention}
\end{equation*}
The complex Hermitian formulation is equivalent by
Theorem~\ref{thm:HR-graph-hermitian} and
Corollary~\ref{cor:HR-real-graph}.

Equation~\eqref{eq.center-invariance} can be written as a semiconjugacy of
skew-products.  Define
\begin{equation*}
    \expoincmap(x,\theta)
    \bydef\bigl(\poincmap(x,\theta),\theta+\rho\bigr),
    \qquad
    \mathcal F_c(\st,\theta)
    \bydef\bigl(\intdyn_c(\st,\theta),\theta+\rho\bigr),
    \qquad
    \mani^c(\st,\theta)
    \bydef\bigl(\param^c(\st,\theta),\theta\bigr).
\label{eq.center-skew-products}
\end{equation*}
Then
$    \expoincmap\circ\mani^c
    =
    \mani^c\circ\mathcal F_c$.
Thus every trajectory of the reduced skew-product that remains inside
$\mathcal D_c\times\T^d$ is mapped by $\mani^c$ to a trajectory of the full
skew-product, up to the truncation and numerical errors of the computed
graph.  This identity is the basic reduction used in what follows.

The general principle of analyzing invariant dynamics on a center manifold
is classical.  Within the parameterization method, an invariant-manifold
embedding and the dynamics conjugate to its restriction are naturally treated
together \cite{CabreFLl05,HaroCLMF16}; for center manifolds of
non-hyperbolic fixed points of maps, the simultaneous construction of the
center manifold and its conjugate center dynamics is developed explicitly in
\cite{BergHWR2020,Sijbrand85}.  In celestial mechanics, center-manifold
and normal-form reductions have long been used to describe nearby periodic
and quasi-periodic dynamics and to construct invariant tori
\cite{Jorba99a,JorbaM99,CastellaJ00}.  High-order
invariant-manifold computations attached directly to invariant tori and
related center structures have also been used in computational dynamics; see
\cite{GimenoJNO21,BarcelonaHM2024}.  These works provide the relevant
background, but the setting here is different: the reference object is
already a reducible quasi-periodic invariant torus.  The same Fourier--Taylor
graph construction used in the rest of the paper returns both its high-order
local center embedding $\param^c$ and the explicit quasi-periodically forced
reduced map $\intdyn_c$.  Nearby invariant objects can therefore be sought
first in dimension $n_c+d$ rather than in the ambient dimension $n+d$, and
only afterwards lifted to the full phase space.

\begin{rmk}[Local character of the reduction]
\label{rmk:center-local}
The parameterization $\param^c$ is local.  No global center manifold is
required here, and a reduced trajectory or invariant torus is relevant to
the full system only while its center coordinates remain in a region where
the truncated graph is accurate, see Section~\ref{sec.numtests}.
\end{rmk}

\subsection{Higher-Dimensional Invariant Tori in the Center Dynamics}
\label{sec.center-initial-guesses}

Let $\vartheta\in\T^m$ denote the additional internal angles of an invariant
torus sought inside the reduced center dynamics, and let $\omega\in\T^m$ be
the corresponding rotation vector.  Here $m$ counts the new angular
directions of that invariant object; it is not the center-state dimension
$n_c=\nt$.  For example, for one elliptic complex-conjugate center pair one
has $n_c=\nt=2$, while an invariant curve in the center variables introduces
one additional angle, $m=1$.  An invariant torus of the reduced center
dynamics is represented by a map
$    g\colon\T^m\times\T^d\to\mathcal D_c
$ satisfying
\begin{equation}
    \label{eq.center-reduced-torus}
    \intdyn_c\bigl(g(\vartheta,\theta),\theta\bigr)
    =
    g(\vartheta+\omega,\theta+\rho).
\end{equation}
The variables $\theta$ contain the frequencies already present in the
original invariant torus, while $\vartheta$ accounts for the additional
frequencies generated by the center dynamics.  The following elementary
identity explains why solving \eqref{eq.center-reduced-torus} is useful in
the full system.

\begin{pro}[Lifting principle]
\label{pro:center-lifting}
Assume that $\param^c$ and $\intdyn_c$ satisfy
\eqref{eq.center-invariance}, and let $g$ satisfy
\eqref{eq.center-reduced-torus}.  Define
\begin{equation}
    \label{eq.center-lifted-param}
    \param^{c,g}(\vartheta,\theta)
    \bydef
    \param^c\bigl(g(\vartheta,\theta),\theta\bigr).
\end{equation}
Then
$
    \poincmap\bigl(\param^{c,g}(\vartheta,\theta),\theta\bigr)
    =
    \param^{c,g}(\vartheta+\omega,\theta+\rho)
$.
\end{pro}

\begin{proof}
Using first \eqref{eq.center-invariance} and then
\eqref{eq.center-reduced-torus},
\begin{align*}
    \poincmap\bigl(\param^{c,g}(\vartheta,\theta),\theta\bigr)
    &=
    \poincmap\bigl(\param^c(g(\vartheta,\theta),\theta),\theta\bigr)
    =
    \param^c\bigl(\intdyn_c(g(\vartheta,\theta),\theta),\theta+\rho\bigr)
    =
    \param^c\bigl(g(\vartheta+\omega,\theta+\rho),\theta+\rho\bigr)
    \\
    &=
    \param^{c,g}(\vartheta+\omega,\theta+\rho). \qedhere
\end{align*}
\end{proof}

Proposition~\ref{pro:center-lifting} separates the computation into a
low-dimensional nonlinear problem and an algebraic lifting step.  No
additional nonlinear solve is required after $g$ has been computed.  If the
components of $(\omega,\rho)$ are rationally independent, the dynamics on
the parameter torus is quasi-periodic.  Equation~\eqref{eq.center-reduced-torus}
also covers resonant rotation vectors, although individual trajectories are
then not necessarily dense in the parameter torus.  If the original torus
of the Poincar\'e skew-product has dimension $d$ and the reduced object
introduces $m$ independent angles, the lifted torus has dimension $d+m$ for
the discrete system.  When the Poincar\'e map removes one phase from an
underlying flow, the corresponding flow torus typically has one additional
dimension.

When such an additional angular direction is created from an elliptic center
pair in a parameter family, the local picture is closely related to the
classical quasi-periodic Hopf (torus/Neimark--Sacker) bifurcation mechanism;
see \cite{BraaksmaB1987,VitoloBS2011}.  

\smallskip

For numerical purposes, the linear center dynamics gives the simplest
initial approximation to \eqref{eq.center-reduced-torus}.  Suppose that the
selected real center bundle has dimension $n_c=2q$ and consists of $q$
non-real conjugate Floquet pairs on the unit circle.  In the corresponding
complex center variables,
\begin{equation*}
    z_j^+
    =
    \lambda_j z_j+\mathcal O(\|z\|^2),
    \qquad
    \lambda_j=e^{2\pi\I\nu_j},
    \qquad
    j=1,\ldots,q,
    \label{eq.center-linear-complex}
\end{equation*}
for rotation numbers $\nu_j\in\T$.  Reducibility makes this linear part
independent of $\theta$.  If the initial guess uses all $q$ elliptic pairs as
independent oscillatory directions, then $m=q$; for sufficiently small
prescribed amplitudes $a_j>0$,
\begin{equation*}
    z_j(\vartheta)=a_j e^{2\pi\I\vartheta_j},
    \qquad
    \omega_j^{(0)}=\nu_j,
    \qquad
    j=1,\ldots,q,
    \label{eq.center-product-seed}
\end{equation*}
provides a zeroth-order $q$-angle product-torus seed.  For one conjugate
pair, $n_c=\nt=2$ and $m=1$, and canonical real center coordinates give
\begin{equation}
    g_0(\vartheta,\theta)
    =
    a
    \begin{pmatrix}
        \cos(2\pi\vartheta)\\
        \sin(2\pi\vartheta)
    \end{pmatrix},
    \qquad
    \omega_0=\nu_c.
    \label{eq.center-circle-seed}
\end{equation}
If a different realification, scaling, or orientation is used, the
corresponding fixed linear transformation must be applied to
\eqref{eq.center-circle-seed}.  
The circle or product-torus ansatz is
specific to elliptic conjugate pairs.  Real center multipliers, resonant
blocks, or more degenerate center spectra require an initial guess adapted to
the corresponding local dynamics.

A more accurate initial approximation may be obtained from a normal form of
the computed reduced map.  Near an elliptic object, such a calculation can
provide approximate action--angle variables, amplitude-dependent
frequencies, and information about resonant terms.  
Classical computational uses of
normal forms and center reductions for constructing invariant tori can be
found in \cite{Jorba99a,JorbaM99}; for an explicit jet-based
normal-form algorithm for Poincar\'e maps, see \cite{GimenoJJZ25}.

An alternative initial guess can be reconstructed from a trajectory of the
reduced skew-product,
\begin{equation*}
    \label{eq.center-reduced-orbit}
    \st_{k+1}=\intdyn_c(\st_k,\theta_k),
    \qquad
    \theta_{k+1}=\theta_k+\rho,
\end{equation*}
starting from $\st_0$ inside the reliable domain of the center
parameterization.  If the orbit appears quasi-periodic, frequency analysis
can be used to estimate $\omega$, after which the samples can be fitted to a
truncated multidimensional Fourier representation 
$
    g_0(\vartheta,\theta)
    $.
Weighted Birkhoff averages combined with the parameterization method provide
one trajectory-based route to rotation numbers and Fourier coefficients for
invariant circles \cite{BlessingM2024,SanderM2020,GomezMS2010}; recent work also develops
reconstruction of higher-dimensional invariant tori from individual
trajectories \cite{RuthKB2026,RandolphP2026}. 

Independently of how the initial pair $(g_0,\omega_0)$ is obtained, it may be
corrected by solving \eqref{eq.center-reduced-torus},  
plus $m$ phase conditions because translations of the additional angles leave the parameterized torus unchanged. 
Related correction and continuation formulations for quasi-periodic invariant tori
are developed in \cite{DieciL95,SchilderOV05}, while
parameterization-based algorithms exploiting invariant bundles and
reducibility are discussed in \cite{CanadellH2017,HaroCLMF16}.  

Once the reduced equation has converged, the full-space approximation follows
directly from \eqref{eq.center-lifted-param} and may, if necessary, be used as
an initial guess for a full-space correction.

\section{Numerical Implementation} \label{sec.numtools}

The implementation of Algorithm~\ref{palg.graph-tor} uses two complementary
representations of the same finite Fourier--Taylor data. Fourier coefficients
are convenient for phase rotations and for the cohomological equations,
whereas values on a uniform phase grid are better suited to nonlinear
composition and jet transport. When the Poincar\'e map is generated by the flow
of an ODE, the polynomial action required in
Step~\ref{palg.graph-tor.jet} is obtained by propagating multivariate jets.
The same data structures also provide inexpensive a posteriori diagnostics for
Fourier resolution, Taylor truncation, Hermitian symmetry, the invariance
defect, and scaling of the selected spectral coordinates to obtain a numerically meaningful local domain. 
Closely related Fourier--Taylor and jet-transport strategies are used in
high-order computations of invariant manifolds of tori; see, for example,
\cite{GimenoJNO21,GimenoJJMZ23}.

Throughout this section we retain the orientation index $\sigma$. Thus
$\poincmap_\sigma$, $\rho_\sigma$, $\Lambda_{L,\sigma}$,
$\Lambda_{N,\sigma}$, and $\intdyn_\sigma$ refer to the directed dynamics of
Section~\ref{sec.inv-eq} and to the quantities used in
Algorithm~\ref{palg.graph-tor}.

\subsection{Fourier--Taylor Representations}\label{sec.ftrep}

The functions manipulated by Algorithm~\ref{palg.graph-tor} depend on the
selected coordinates $\st\in\C^{\nt}$ and on the angle
$\theta\in\T^d$. It is therefore natural to represent them by truncated
Fourier--Taylor series, with Taylor coefficients in $\st$ that are periodic
functions of $\theta$.  
In computations, we store only finitely many Taylor orders and Fourier modes:
we truncate the Taylor expansion at order $N$ and choose
$\bar M=(\bar M_1,\dots,\bar M_d)\in\N^d$, which defines the Fourier index set
\[
    \mathcal I = \{k=(k_1,\dots,k_d)\in\Z^d :
    -\bar M_j\leq k_j\leq \bar M_j,\ j=1,\dots,d\}.
\]
For a quantity with $q$ output components, write $F\colon \C^{\nt}\times\T^d\to\C^q$ as 
\[
    F(\st,\theta)
    \approx
    \sum_{l=0}^{N} F_{(l)}(\st,\theta), \quad 
    F_{(l)}(\st,\theta)
    =\sum_{|\alpha|=l}F_{l,\alpha}(\theta)\st^\alpha, \quad 
    F_{l,\alpha}(\theta) = 
    \sum _{k\in \mathcal I}
    \widehat F_{l,\alpha,k}
    e ^{2\pi \I \innerp{k}{\theta}}\st ^\alpha
\]
where $\alpha\in\N^{\nt}$ is a Taylor multi-index and $k\in \mathcal I$ is a Fourier
multi-index with elements in $\Z^d$. This representation amounts to storing the finite set of coefficients
$\{\widehat F _{l,\alpha,k}\colon |\alpha|\leq N,\ k\in \mathcal I\}$. Note that, the number of Taylor coefficients in $\nt$ variables up to degree $N$ is
$\binom{\nt + N}{\nt}$.

Equivalently, each Taylor coefficient $F_{l,\alpha}(\theta)$ can be stored in grid
form. Let $M=(M_1,\dots,M_d)$, with $M_j=2\bar M_j+1$. 
We define the
grid-index set and the corresponding uniform grid by
\[
    I = \{i=(i_1,\dots,i_d)\in\Z^d : 0\leq i_j<M_j,\ j=1,\dots,d\},
    \qquad
    \mathcal G = \{\theta_i=(i_1/M_1,\dots,i_d/M_d) : i\in I\}.
\]
In this grid form, we
store the values $\{F_{l,\alpha,i}\colon |\alpha|\leq N,\ i\in I\}$, where
\begin{equation*}
    F_{l,\alpha,i} \bydef F_{l,\alpha}(\theta_i),
    \qquad \theta_i \in \mathcal G.
\end{equation*}
This representation also defines, at each grid point, a Taylor series
\begin{equation}\label{eq.taygrid}
    F_i(s) \bydef \sum_{0\leq l \leq N} F_{(l),i}(s)
    = \sum_{0\leq l \leq N}\sum _{|\alpha|=l} F_{l,\alpha,i}s^\alpha,
    \qquad i \in I,
\end{equation}
which we call the Taylor-grid representation.

The Fourier--Taylor and Taylor-grid representations are related by the discrete
Fourier transform, which is efficiently evaluated by the Fast Fourier Transform
in $\mathcal O(M_{\rm grid}\log M_{\rm grid})$ operations for each Taylor coefficient, where
$M_{\rm grid}=\prod_{j=1}^d M_j$ is the total number of grid points.

For genuinely real grid data, such as the quantities obtained after the
realification described in Section~\ref{sec.num-real}, the angular Fourier
coefficients satisfy
$\widehat F_{l,\alpha,-k}=\overline{\widehat F_{l,\alpha,k}}$.
A real-to-complex transform may therefore store only a half-spectrum; with the
last angular direction used for the real transform, the number of stored
frequency slots per scalar Taylor coefficient is
$(\bar M_d+1)\prod_{j=1}^{d-1}M_j$. For complex-valued functions, no such
Hermitian condition is available, and all complex coefficients must be stored.

The choice of representation follows the algebra of
Algorithm~\ref{palg.graph-tor}. Steps~\ref{palg.graph-tor.cohom} and
\ref{palg.graph-tor.rot} are naturally carried out in Fourier space, while
Steps~\ref{palg.graph-tor.jet} and \ref{palg.graph-tor.comp} are carried out on
the phase grid with truncated multivariate polynomials. In particular,
Step~\ref{palg.graph-tor.jet} requires the action of the directed Poincar\'e map
$\poincmap_\sigma$ on the polynomial parameterization. The required polynomial
composition is obtained by jet transport when $\poincmap _\sigma$ is generated
by an ODE flow.


\subsection{Formal Power Series and Jet transport}
Step~\ref{palg.graph-tor.jet} in Algorithm~\ref{palg.graph-tor} involves the composition of possibly a Poincar\'e map and the torus up to high order. To achieve this, 
one has to integrate numerically the variational equations (associated to the ODE) up to the desired order. The most efficient technique to do so is the so-called Jet Transport (JT). 

\subsubsection*{Automatic Differentiation (AD)}
Forward AD evaluates partial derivatives of a computer program's output with respect to its
inputs. The primary requirement is that the program consists of a sequence of compositions
of elementary functions (specifically, those admitting local Taylor expansions). Far from
being numerical or symbolic differentiation, AD occupies a semi-analytic middle ground: it
computes exact numerical values for derivatives by systematically applying algebraic
differentiation rules (such as the chain rule), combining the computational efficiency of
numerical methods with the accuracy of symbolic approaches. The only source of numerical
inaccuracies are standard floating-point round-off errors.

The underlying mechanism relies on Taylor's theorem. Consider a univariate function $g$. By
treating a small displacement $\delta$ as a formal algebraic symbol satisfying
$\delta^{N+1}=0$, we can express the truncated Taylor expansion as a polynomial of degree
$N$:
$$
g(x+\delta) = \sum_{k=0}^{N} \frac{g^{(k)}(x)}{k!} \delta^k.
$$
The vector of coefficients $(x_0, x_1, \dots, x_N)$ represents a jet (or a generalized dual
number), where the $k$-th order coefficient corresponds to the normalized derivative
$g^{(k)}(x)/k!$. This algebraic structure naturally generalizes to multivariate functions
$g\colon \mathbb{R}^m \to \mathbb{R}$, where evaluation is carried out over the algebra of
multivariate truncated polynomials or power series up to a given order $N$. 

\subsubsection*{Jet Transport (JT): AD applied on numerical steppers}
Consider a general Cauchy problem:
$$
\begin{cases}
    \dot{\mathbf{x}} = X(\mathbf{x},t), \\
    \mathbf{x}(t_0) = \mathbf{x}_{0},
\end{cases}
$$
integrated via a numerical stepper (explicit or implicit) of the form $\mathbf{x}_{n+1} = \mathbf{x}_{n} + \Phi(h, \mathbf{x}_{n}, \mathbf{x}_{n+1})$, 
which generates an approximate discrete trajectory $\{\mathbf{x}_n\}_{n\geq 0}$ with step
size $h$. By applying Automatic Differentiation (AD) to the initial conditions
(effectively evaluating the stepper over a high-order algebra of dual numbers)
the numerical procedure yields a high-order Taylor expansion of the flow map with respect
to $\mathbf{x}_0$. This technique, widely adopted since the mid-2000s, has found
significant utility in fields like astrodynamics.

A key question is the convergence and error propagation of these derivatives, given that AD
evaluates the derivatives of the \textit{discrete stepper} rather than the continuous flow.
Initial foundational work established that applying 
first-order AD to Runge--Kutta schemes approximates the first variational equation with the
same order of convergence as the underlying trajectory \cite{Walther2007}. This result was later
generalized to arbitrary orders of the variational equations and also for Taylor methods \cite{GimenoJJMZ23}, and for generalized linear methods \cite{Forrier23}.
Crucially, these generalizations prove that high-order JT is algorithmically equivalent to
applying the same numerical stepper to the corresponding higher-order variational equations
after eliminating algebraic redundancies. In this context, algorithmic equivalence
signifies that, for identical inputs, both approaches perform the exact same set of
arithmetic operations, differing at most in their execution sequence.

\subsection{Scaling Factor and Local Validity Domain} \label{sec.scaling}

A finite Taylor expansion is useful only on a bounded neighborhood of the
origin in the selected coordinates, and the apparent size of that neighborhood
depends on the normalization of the selected Floquet basis. The scaling in
Step~\ref{palg.graph-tor.scaling} exploits this freedom without changing the
geometric manifold: only $\Vt$ is rescaled, while the complementary basis
$\Vh$ is left unchanged.


There are several ways to numerically estimate the accuracy of these expansions and to determine the domain of validity where they provide reliable approximations.
Here we review the criteria that will be used later to validate our computations.

To monitor the accuracy of the computed Fourier-Taylor series $F\colon \C^{\nt} \times
\T^d \rightarrow \C^{q}$,
we introduce the following quantities for each degree $l$.
\begin{equation}\label{eq.norminf1}
    a_\infty^F(l)
    \bydef
    \sup_{\theta \in \T^d}
    \max_{|\alpha|=l}
    \|F_{l,\alpha}(\theta)\|_\infty, 
    \qquad \text{and} \qquad
    a_{1}^F(l)
    \bydef
    \sup_{\theta \in \T^d} \sum_{|\alpha|=l} \|F_{l,\alpha}(\theta)
    \|_\infty.
\end{equation}
The first quantity monitors individual coefficient growth. The second gives a
simple bound for a homogeneous term: if $\|s\|_\infty\leq r$, then
\begin{equation}\label{eq.num-homogeneous-bound}
    \sup_{\theta\in\T^d}
    \|F_{(l)}(s,\theta)\|_\infty
    \leq r^l a_1^F(l).
\end{equation}
\newcommand{\scfac}{\mathfrak{p}}
Both quantities \eqref{eq.norminf1} can be numerically evaluated
using the Taylor-grid representation in \eqref{eq.taygrid}.

\smallskip

For the graph, let $\mathcal U\subset\C^{\nt}$ denote a candidate local domain
in the selected coordinates and set
$
    r\bydef\sup_{s\in\mathcal U}\|s\|_\infty
$. 
Using \eqref{eq.num-homogeneous-bound} for the last retained term of
\eqref{eq.GT.W-FPS}, a convenient tail-based radius estimate is
\begin{equation}\label{eq.rmax}
    r_{\rm tail}
    \bydef
    \left(
      \frac{\varepsilon_{\tol}}{a_1^\graphfun(N)}
    \right)^{1/N},
\end{equation}
whenever $a_1^\graphfun(N)\neq0$. Here $\varepsilon_{\tol}$ is a prescribed
threshold for the contribution of $\graphfun_{(N)}$. Equation~\eqref{eq.rmax}
is a practical truncation indicator rather than a convergence theorem; lower
orders, Fourier resolution, and the invariance defect should be checked as
well.

A first run with $\mathfrak p=1$ can therefore be used to estimate a useful
coordinate scale. Setting
\begin{equation*}\label{eq.num-scaling-choice}
    \mathfrak p=r_{\rm tail},
    \qquad
    V_L\mapsto r_{\rm tail}V_L
\end{equation*}
and recomputing the expansion brings the effective domain in the rescaled
$s$-coordinates to order one. This is precisely the scaling operation in
Step~\ref{palg.graph-tor.scaling}; no rescaling of $V_N$ is involved.

The domain can also be assessed directly from the parameterization equation
\eqref{eq.GT.directed-PM}. Given
$x_0=\param(s_0,\theta_0)$, let $(s_n,\theta_n)$ be generated by the internal
dynamics \eqref{eq.GT.induced-dynamics} and the rotation $\rho_\sigma$, and let
$(x_n,\theta_n)$ be the corresponding orbit of the directed skew-product of
Section~\ref{sec.inv-eq}. We use the one-step invariance defect along the
reduced orbit,
\begin{equation}\label{eq.inverror}
    e_I(s_0,\theta_0,n)
    \bydef
    \left\|
      \poincmap_\sigma\bigl(\param(s_n,\theta_n),\theta_n\bigr)
      -\param\bigl(\intdyn_\sigma(s_n,\theta_n),
                    \theta_n+\rho_\sigma\bigr)
    \right\|_\infty,
\end{equation}
and the accumulated orbit discrepancy
\begin{equation}\label{eq.orbiterror}
    e_O(s_0,\theta_0,n)
    \bydef
    \|\param(s_n,\theta_n)-x_n\|_\infty.
\end{equation}
The first quantity measures the residual of \eqref{eq.GT.directed-PM}
directly; the second measures how this local residual accumulates under
iteration. In strongly hyperbolic dynamics, long orbit comparisons may be
dominated by amplification of integration and roundoff errors, so
\eqref{eq.inverror} is usually a more reliable local indicator. 
See also \cite{tesilebihan}.

Let $\mathcal Q_R\subset\{s\in\C^{\nt}\colon \|s\|_\infty\leq R\}$ be a finite mesh
of selected coordinates, let $\mathcal G$ be the angular grid defined in Section~\ref{sec.ftrep},
and let $n_{\rm iter}$ be the number of iterates used in the test. For
$m\in\{I,O\}$, define
\begin{equation*}\label{eq.aggregate-domain-error}
    E_m(R)
    \bydef
    \max_{s_0\in\mathcal Q_R,\,\theta_0\in\mathcal G}
    \max_{0\leq n\leq n_{\rm iter}}
    e_m(s_0,\theta_0,n).
\end{equation*}
The largest tested $R$ for which $E_m(R)\leq\varepsilon_{\tol}$ provides an
a posteriori estimate of the practical domain associated with the selected
error indicator. In applications we use this test together with the coefficient
information in \eqref{eq.norminf1}--\eqref{eq.rmax}, rather than interpreting
either criterion alone.

\subsection{Numerical Validity Tests} \label{sec.numtests}

We describe different tests that help us to numerically validate the
computations made by the different steps of Algorithm~\ref{palg.graph-tor} and to assess the quality of the obtained graph $\graphfun$.
Once the graph $\graphfun$ has been computed, we test the invariance equations
\eqref{eq.GT.graph-invariance-abstract} and \eqref{eq.GT.directed-PM} and the quality of the Fourier--Taylor representation used to store
the solution.
The algorithm~\ref{palg.graph-tor} is not an iterative method for
the invariance equation \eqref{eq.GT.graph-invariance-abstract}; instead each order is obtained by matching
powers in the
manifold coordinates and solving the corresponding cohomological equation. Therefore,
the tests below are not stopping criteria, but a posteriori checks of
the computed parameterizations, although some of them can be used during the iterative procedure. 

The following test measures the defect of the computed torus parameterization
and of the Floquet change in Algorithm~2.3 in \cite{GimenoJNO21}, and is used as
a stopping criterion. Moreover, in our computations, equation~\eqref{eq.GT.graph-invariance-abstract} can be
verified order by order through its grid representation.

\begin{test}[Defect test] \label{test1}
Let $\mathcal{A}$ be a mesh in $\T ^d$. A function
$z\colon \T^d \rightarrow \R^n$ is said to verify the equation
$\mathcal{I}(z(\theta))=0$ with tolerance $\tol$ if, and only if,
$
 \max \limits_{\theta\in \mathcal{A}}
  \ \|\mathcal{I}(z(\theta))\| \leq \tol
$.
\end{test}
Test~\ref{test1} will be used with different equation residuals: the torus or Floquet
defect, the order-by-order graph defect, and the cohomological residual. This
keeps the tolerance convention independent of the particular equation being
checked.

\begin{test}[Remeshing test] \label{test2}
Let $\mathcal{A}$ be a mesh in $\T ^d$. A periodic function $z$ is said to
verify the equation $\mathcal{I}(z(\theta))=0$ with tolerance $\tol$
and shift $\gamma \in \T ^d$ if, and only if,
$
 \max \limits_{\theta\in \mathcal{A}}
  \ \|\mathcal{I}(z(\theta + \gamma))\| \leq \tol
$.
\end{test}
The remeshing test detects errors that can be hidden by evaluating only at the
collocation grid. The shift is applied spectrally, so no additional
interpolation formula is required.

\smallskip

The next test checks whether the truncated Fourier representation is
compatible with the chosen mesh size. We assume that the functions involved are
analytic, so their Fourier coefficients decay exponentially. We distinguish two
cases, depending on whether the representation is real or complex.

\begin{test}[Real Fourier-tail test] \label{test3-foureal}
A truncated real Fourier representation given by
\[
x(\theta) = \frac{x_{0}}{2} + \sum _{|k| = 1}^N x ^{(c)}_k \cos \innerp{k}{\theta} +  x ^{(s)}_k \sin \innerp{k}{\theta}, \qquad \theta \in \T ^d\text{, } k \in \N ^{d}.
\]
is said to verify Test~\ref{test3-foureal} with tolerance $\tol$ if, and
only if, for $k \in \N^d$
$
 \max \limits_{|k | \in \{N-1, N\}} \|(x _k^{(c)}, x _k^{(s)})\|_2 \leq \tol
$.
\end{test}

\begin{test}[Complex Fourier-tail test] \label{test4-fou}
A truncated complex Fourier representation given by
\[
x(\theta) = \sum _{|k| = -N}^N \widehat x _k e ^{2\pi \I \innerp{k}{\theta}}, \qquad \theta \in \T ^d\text{, } k \in \Z ^{d}.
\]
is said to verify Test~\ref{test4-fou} with tolerance $\tol$ if, and
only if, for $k \in \N^d$
$
 \max \limits_{|k | \in \{-N,-N+1,N-1, N\}} |\widehat x _k| \leq \tol
$.
\end{test}

Similarly, we can check whether the order of the Taylor truncation is large enough by
looking at the size of the last-order coefficients. 
For a Fourier-Taylor
representation, the same test can be applied pointwise in the grid variable,
that is, to the Taylor-grid representation.

\begin{test}[Taylor-tail test] \label{test5-taylor}
A truncated Taylor representation given by
\[
    x(s) = \sum_{|\alpha|=0}^N x_\alpha s^\alpha, \quad \alpha \in \N^{\nt}
\]
is said to verify Test~\ref{test5-taylor}
with tolerance $\tol$ if, and only if,
$
 \max \limits_{|\alpha|=N} |x _\alpha| \leq \tol
$.
\end{test}

For computations that use the real formulation of Section~\ref{sec.num-real},
the Taylor coefficients must also preserve the conjugate-pair structure
established in Section~\ref{sec.graph.herm}. The relevant index involution is
already defined by $\alpha^\star\bydef J_L\alpha$; see
\eqref{eq.Jmat} and \eqref{eq.kraw-global-hermitian}.

\begin{test}[Hermitian-symmetry test]\label{test6-conj}
Assume the canonical ordering of
Section~\ref{sec.canonical-conjugate-pairs}. The graph coefficients satisfy the
Hermitian-symmetry test with tolerance $\tol$ if
\begin{equation*}\label{eq.num-hermitian-W-test}
    \max_{2\leq l\leq N}
    \max_{|\alpha|=l}
    \max_{\theta\in\mathcal G}
    \left\|
      w_{l,\alpha^\star}(\theta)
      -J_N\overline{w_{l,\alpha}(\theta)}
    \right\|_\infty
    \leq\tol.
\end{equation*}
Writing
$\intdyn_{\sigma,(l)}(s,\theta)
 =\sum_{|\alpha|=l}f_{\sigma,l,\alpha}(\theta)s^\alpha$,
the analogous test for $\intdyn _{\sigma,(l)}$ replaces $J_N$ by $J_L$.
\end{test}
If Test~\ref{test6-conj} is evaluated directly in Fourier space, complex
conjugation also reverses the angular Fourier index: the graph relation becomes
$\widehat w_{l,\alpha^\star,k}
 =J_N\overline{\widehat w_{l,\alpha,-k}}$ up to the chosen tolerance. 

\smallskip

Finally, the expected order of the first omitted Taylor term can be checked from
the invariance defect. Let $\param^{\leq N}$ be the parameterization associated
with the computed $\graphfun ^{\leq N}$ and $\intdyn _\sigma ^{\leq N}$, and define the residual of
\eqref{eq.GT.directed-PM} by
\begin{equation}\label{eq.num-final-invariance-residual}
    \epsilon _N(\st,\theta)
    \bydef
    \poincmap_\sigma\bigl(\param^{\leq N}(\st,\theta),\theta\bigr)
    -\param^{\leq N}\bigl(
       \intdyn_\sigma^{\leq N}(\st,\theta),
       \theta+\rho_\sigma
     \bigr).
\end{equation}

\begin{test}[Scaling and coefficient-growth test]\label{test7-growth}
Choose a nonzero $\st_1$ in the tested domain. Let $\st_2=\st_1/2$ and 
$
    c _i
    \bydef
    \max_{\theta\in\mathcal G}
    \|\epsilon _N(\st_i,\theta)\|_\infty
$, for $i = 1,2$.
If the dominant contribution to \eqref{eq.num-final-invariance-residual} is the
first omitted homogeneous term, then
$c_i\approx C\|s_i\|_\infty^{N+1}$ and
\[
    \frac{\log c_1-\log c _2}{\log 2}
    \approx N+1.
\]
The comparison is meaningful only while the invariance defect remains above other error floors, such as time-integration, Fourier-truncation, and roundoff errors. So this test will pass for suitable $\st _1$ and low enough $N$.
\end{test}

\subsubsection{List of tests on Algorithm~\ref{palg.graph-tor}}
\label{sec.numtests-palg}

During the implementation of Algorithm~\ref{palg.graph-tor}, several tests can
be performed to check that the previous steps are accurate to catch implementation errors at the degree where they occur. They complement the final invariance check, the remeshing test, and the truncation order tests.

First, after Step~\ref{palg.graph-tor.comp}, the already determined homogeneous
orders must satisfy the graph equation. Since $Z$ represents the composed
right-hand side at this stage, for every $q<l$ one should have
\begin{equation}\label{eq.num-composition-residual}
    Z_{(q)}(s,\theta)
    =\cR_{N,\sigma,(q)}(s,\theta).
\end{equation}
Thus the grid residual associated with
\eqref{eq.num-composition-residual} should satisfy Test~\ref{test1}. A failure
at an order $q<l$ normally indicates an inconsistency in previously stored
coefficients, in the jet composition, or in a representation change.

Second, after Step~\ref{palg.graph-tor.cohom}, each newly computed coefficient
can be substituted back into the cohomological equation used to obtain it.
Test~\ref{test1}, either on the phase grid or mode by mode after the Fourier
transform. 

\smallskip

Table~\ref{tab:numtests-algorithm} summarizes tests perform in the implementation of the algorithm.

\begin{table}[htbp]
\centering
\small
\renewcommand{\arraystretch}{1.75}
\begin{tabular}{@{}p{0.13\textwidth}|p{0.19\textwidth}|p{0.63\textwidth}@{}}
\hline
{\bf Stage} & {\bf Test or residual} & {\bf Object on which it is applied} \\
\hline
Step~\ref{palg.graph-tor.scaling}
& \eqref{eq.norminf1}--\eqref{eq.rmax}, Test~\ref{test7-growth}
& graph $\graphfun^{\leq N}$ and selected basis $V_L$ to choose
  $\mathfrak p$ before recomputing the expansion. \\

Step~\ref{palg.graph-tor.jet}
& Tests~\ref{test3-foureal} or \ref{test4-fou}, and
  Test~\ref{test5-taylor}
& Transported jet $\cR_\sigma$ (or its $L$ and $N$ blocks); checks Fourier
  resolution and the Taylor degree after applying $\poincmap_\sigma$. \\

Step~\ref{palg.graph-tor.intdyn}
& Test~\ref{test6-conj} when real formulation
& New homogeneous internal term $\intdyn_{\sigma,(l)}$.
  \\

Step~\ref{palg.graph-tor.comp}
& Test~\ref{test1}
& Lower-order composition defect \eqref{eq.num-composition-residual}; verifies that previously computed graph   coefficients still satisfy the invariance equation after recomposition. \\

Step~\ref{palg.graph-tor.cohom}
& Test~\ref{test1} and, when appropriate, Test~\ref{test6-conj}
& Cohomological residual for each $w_{l,\alpha}$ and the newly computed graph
  coefficients; checks both the linear solve and Hermitian symmetry. \\

Step~\ref{palg.graph-tor.rot}
& Test~\ref{test3-foureal} or Test~\ref{test4-fou}
& Rotated homogeneous term
  $R_{(l)}(s,\theta)=\graphfun_{(l)}(s,\theta+\rho_\sigma)$; checks that the
  shifted coefficient remains resolved in the retained Fourier box. \\

After degree $N$
& Tests~\ref{test1}, \ref{test2}, and \ref{test7-growth}
& $\param^{\leq N}$, $\graphfun^{\leq N}$, and
  $\intdyn_\sigma^{\leq N}$ through the directed invariance residual
  \eqref{eq.num-final-invariance-residual}. \\
\hline
\end{tabular}
\caption{Numerical checks performed at the indicated stages of
Algorithm~\ref{palg.graph-tor}.}
\label{tab:numtests-algorithm}
\end{table}

\section{Examples} \label{sec.examples}

In this section, we apply the method to several models to illustrate the
different computational settings discussed above: the computation of
parameterizations using complex-valued or real-valued arithmetic, and the use of
a multiple-shooting implementation.

The explorations presented in Sections~\ref{sec.pendulum} and \ref{sec.sssm} have
been carried out on a workstation equipped with an 8-core AMD Ryzen 7 3700X
processor (16 SMT threads, maximum frequency 4.4~GHz) and 16~GB of RAM, under Debian
GNU/Linux 13 (kernel 6.12).
The real-versus-complex comparison
of Section~\ref{sec.solarsail} was run on a second workstation with an
AMD Ryzen 9 5900X processor with 10 threads.
The source
code was written in C, compiled with GCC 14.2.0, and linked against OpenBLAS
0.3.29 (LAPACK interface) and FFTW 3.3.10. The implementation also uses
OpenMP~4.5 in most of the computations. 
The numerical integrations were performed as follows. Sections~\ref{sec.pendulum} and \ref{sec.solarsail} use Taylor integration (Taylor package \cite{GimenoJZ22}) with tolerance $10^{-16}$, whereas Section~\ref{sec.sssm} uses a
Runge-Kutta-Verner method of order 8(9) with jet transport and tolerance $10^{-14}$ \cite{Verner78}.
All reported times are elapsed (wall-clock) seconds.
The code of Section~\ref{sec.solarsail} was obtained with the help of OpenAI model 
\textit{GPT-G.1 Sol} using this paper and related documentation as prompts. 

\subsection{Quasi-periodically forced pendulum}\label{sec.pendulum}
The system describes the movement of a
quasi-periodically forced pendulum coupled with a simple equation. It extends the example in \cite{GimenoJNO21} with a coordinate $z$ to illustrate the center manifold reduction.
\begin{equation}
\begin{aligned}
\dot{x} &=y \\
\dot{y} &=-\alpha \sin x + \varepsilon \zeta (\theta_0, \theta_1, \theta _2) \\
\dot{z} &=-\beta z + \varepsilon \eta (x, \theta_0,\theta_1, \theta _2), \\
\dot{\theta} _j &= \omega _j/(2\pi), \ \  j=0,1,2
\end{aligned}
\label{eq.FP}
\end{equation}
where $x,y \in \R$, and $\alpha$ is a parameter
. For $\theta_j \in \T$ and $\varepsilon$
accounts for the weight of the forcing function $\zeta$:
\begin{equation*}
    \zeta (\theta_0, \theta_1, \theta _2)=\left[ 5 + \sum_{j=0}^{2}
      \cos 2\pi \theta_j \right]^{-1} \qquad
    \eta(x,\theta _0, \theta_1, \theta _2) = \frac{1 - \cos(x)}{\zeta (\theta_0,\theta_1, \theta _2)}.
\end{equation*}
As frequencies we have chosen, 
\begin{equation*}
\omega_0=1, \ \ \omega_1=\sqrt{2}, \ \ \omega_2=\sqrt{3}.
\label{eq.FP_freq}
\end{equation*}
The parameter values in \eqref{eq.FP} are $\alpha=\beta=\frac{1}{2}$.
The system \eqref{eq.FP} perturbs the Hamiltonian at $\varepsilon = 0$ as
\begin{equation}
\label{eq.energy-like}
 H_{\varepsilon}(x,y;\theta _0, \theta _1, \theta _2) = \frac{y^2}{2} - \alpha \cos (x) - \varepsilon x\zeta(\theta _0, \theta _1, \theta _2).
\end{equation}

To reduce the dimension of \eqref{eq.FP}, we take $\poincmap$ to be the return map to
the section $\theta_0=0$. Thus, $\poincmap$ is the Poincar\'e map associated with the
frequency $\omega_0$, or equivalently the map obtained by integrating for
$T=2\pi/\omega_0$ units of time. For the resulting discrete system, we
apply Algorithm~2.3 in \cite{GimenoJNO21} with $\varepsilon=0.015$ to compute a
torus, together with its Floquet transformation and Floquet matrix, near
$x=y=0$. The torus has dimension $d=2$ under $\poincmap$ and is a response torus to the quasi-periodic perturbation whose inner dynamics is prescribed by the frequencies $\omega_1$, $\omega_2$. Dynamically, it is 
equivalent to the fixed point at the origin of the unperturbed pendulum.
Algorithm~2.3 in \cite{GimenoJNO21} is run using 32 Fourier modes in each
angular variable and a prescribed Newton tolerance of $10^{-10}$. After four
iterations, the obtained torus and Floquet transformation satisfy
\eqref{eq.error0}, according to Test~\ref{test1} and
Test~\ref{test2}, with errors below $10^{-12}$. 

\begin{figure}[ht]
    \centering
    \begin{tabular}[c]{@{}c@{\hspace{0.015\textwidth}}c@{\hspace{0.015\textwidth}}c@{}}
        \includegraphics[width=0.20\textwidth,page=3]{pendol-crop.pdf} &
        \includegraphics[width=0.35\textwidth,page=1]{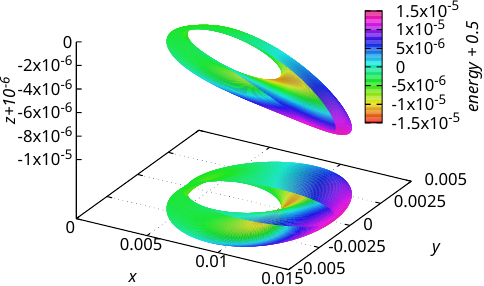} &
        \includegraphics[width=0.20\textwidth,page=4]{pendol-crop.pdf}
    \end{tabular}
    \caption{Invariant response torus of \eqref{eq.FP}, shown in two planar projections and 3D coordinates, colored by \eqref{eq.energy-like}.}
    \label{fig:pendol}
\end{figure}

With this choice of parameters, the spectrum of the matrix $B$ is as follows
\[
\begin{split}
    \lambda_{1,2}
    &= -2.663122512970447 \times 10^{-1}
    \pm \I 9.638868112019690 \times 10^{-1},
    \qquad
    \lambda_3 = 4.321391826377226 \times 10^{-2},
\end{split}
\]
which indicates that the linear behavior around the torus is center$\times$stable. Considering $v_j$ to be the eigenvector of $B$ associated to the eigenvalue $\lambda_j$, we take $V$ in \eqref{eq.GT.spectral-blocks} to be $V=\begin{bmatrix}
    v_1 & v_2 & v_3
\end{bmatrix}$, 
such that every $v_j$ has norm 1, and choose $n_L=2$, $n_N=1$.
With this ordering, we run Algorithm~\ref{palg.graph-tor} using the default
scaling factor, and thus the parameterization given by
\eqref{eq.GT.K-graph} locally provides the center manifold attached to the
invariant response torus described as a graph given by the truncated expansion
$\graphfun$, see Section~\ref{sec.center-reduction}.
In Table~\ref{tab:timevsorder} we present the total computational time
spent to compute the expansion $\graphfun$ up to different degrees $N$ of the Taylor series. The number of Fourier coefficients has been fixed to the original number used to parameterize the torus. Our
implementation has been parallelized and uses a total of $8$ threads.
In order to check performance of the parallelization, we include
Table~\ref{tab:timevsthreads} where an expansion $\graphfun$ of degree $N=16$ is
computed using a different number of threads.
Observe that doubling the number of threads does not halve the time (the speedup per doubling is between 1.3 and 1.4), because of the barriers between the steps.

\begin{table}[ht]
    \centering

    \begin{minipage}[t]{0.46\textwidth}
        \centering
        \small
        \begin{tabular}{r||c|c|c|c}
            {\bf Degree} & 8 & 10 & 15 & 20 \\ \hline
            Time (s) & 54.84 & 88.26 & 168.20 & 324.07
        \end{tabular}
        \caption{Computational time (wall-clock) to compute
        $\graphfun$ up to several degrees using a parallelized
        implementation with $8$ threads.}
        \label{tab:timevsorder}
    \end{minipage}
    \hfill
    \begin{minipage}[t]{0.52\textwidth}
        \centering
        \small
        \begin{tabular}{r||c|c|c|c}
            {\bf Num. of threads} & 1 & 2 & 4 & 8 \\ \hline
            Time (s) & 472.72 & 365.17 & 267.55 & 185.34
        \end{tabular}
        \caption{Computational time (wall-clock) to compute
        $\graphfun$ up to degree 16 versus the number of threads.}
        \label{tab:timevsthreads}
    \end{minipage}
\end{table}

\begin{figure}[ht]
    \centering
    \includegraphics[width=0.32\linewidth]{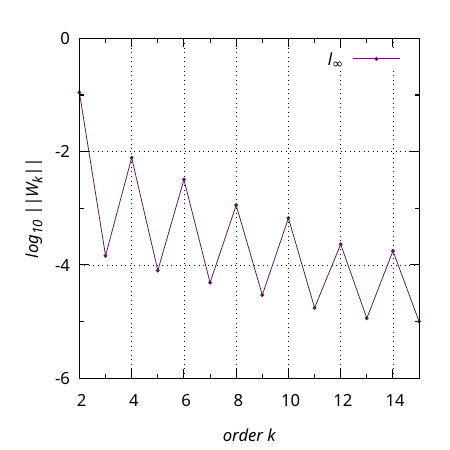}
    \includegraphics[width=0.32\linewidth]{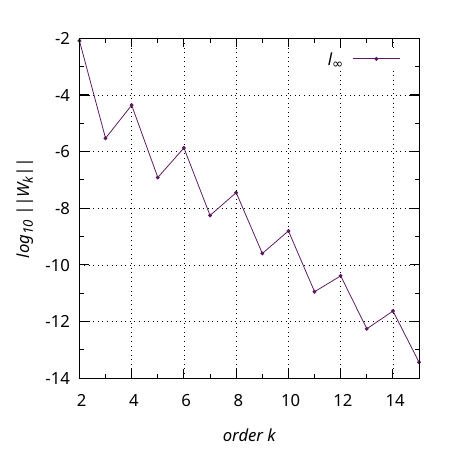}
    \caption{Decay of the coefficients of the graph-parameterization method accordingly to $a_\infty^\graphfun(l)$. From left to right: non-scaled parameterization and scaled using \eqref{eq.rmax}.}
    \label{fig:norms_subinf}
\end{figure}

\begin{figure}[ht]
    \centering
    \includegraphics[width=0.32\linewidth]{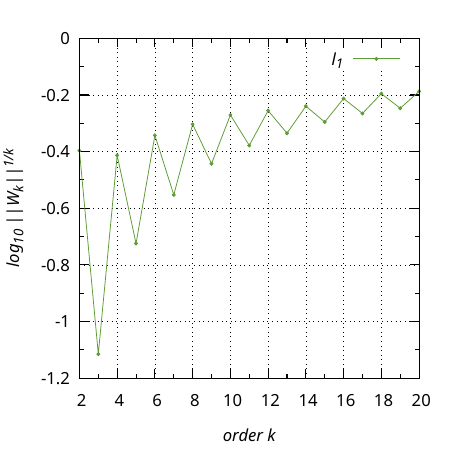}
    \includegraphics[width=0.32\linewidth]{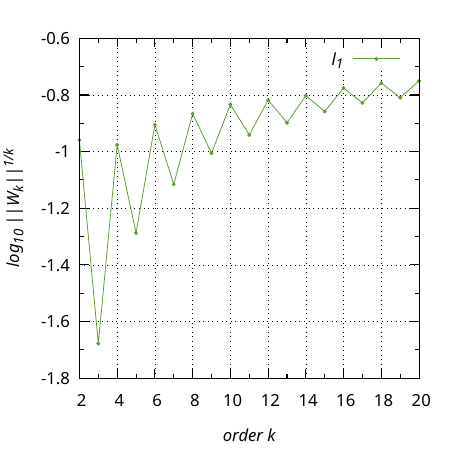}
    \caption{Estimated radius of convergence of the series $\graphfun$ using $(a_1^\graphfun(l))^{1/l}$. From left to right: non-scaled parameterization and scaled using \eqref{eq.rmax}. 
    }
    \label{fig:norms_sub1}
\end{figure}

In order to have an idea of the divergence of the expansion $\graphfun$, the decay of
its coefficients, and its radius of convergence, we compute the quantities \eqref{eq.norminf1}.
Results are presented in
Figures~\ref{fig:norms_subinf} and \ref{fig:norms_sub1} for a fixed degree $N=20$ of the Taylor coefficients. The left panels show
the quantities $a_\infty^\graphfun(l)$ and $(a_1^\graphfun(l))^{1/l}$, respectively, for the
expansions computed using the default scaling factor $\scfac=1$, whereas the
right panels show the same quantities but using a scaling factor $\scfac$
accordingly to \eqref{eq.rmax}, where $\varepsilon_\tol=10^{-15}$. 
Observe from Figure~\ref{fig:norms_subinf}, that the scaling factor produces a much faster decay of the coefficients: in this case, the last coefficients are of order of the machine precision. At the same time, from Figure~\ref{fig:norms_sub1}, we observe that the growth of the coefficients is mild and, up to this order, they should provide good approximations in a neighborhood around the invariant torus.

\begin{figure}[ht]
    \centering
    \begingroup
    \setlength{\tabcolsep}{2pt}
    \noindent\resizebox{\linewidth}{!}{%
      \begin{tabular}{@{}*{4}{c}@{}}
        \includegraphics[height=3cm]{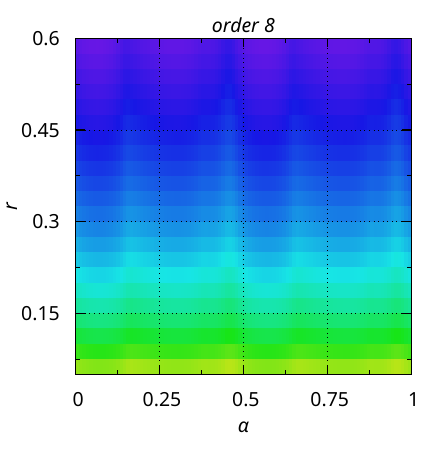} &
        \includegraphics[height=3cm]{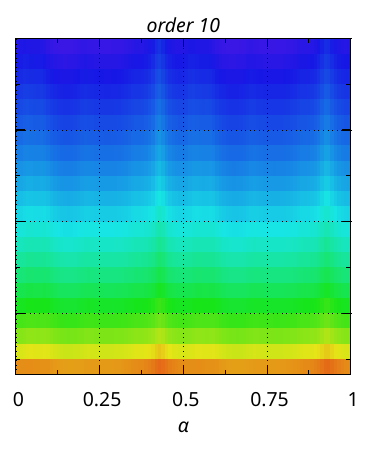} &
        \includegraphics[height=3cm]{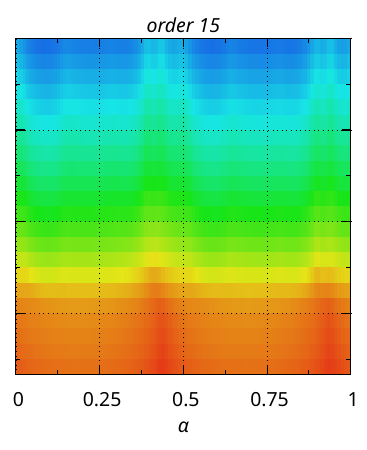} &
        \includegraphics[height=3cm]{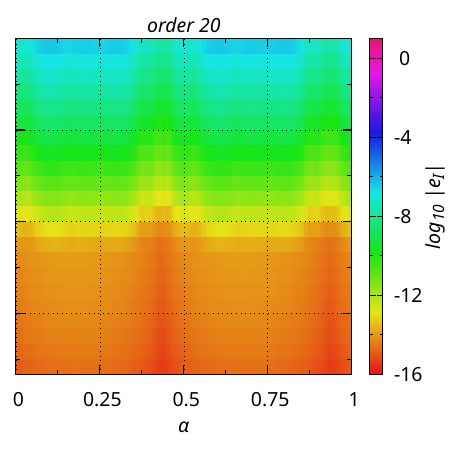}
    \end{tabular}%
}
\endgroup
    \caption{
    Estimated domain of validity of $\graphfun$ for several truncation degree
    $N$, using the default scaling factor. Colors represent
    $\max_{\theta_0\in\mathcal G}\log_{10} e_I(s_0,\theta_0,1)$, with $s_0=r e^{\I
    2\pi\alpha}$. The horizontal axis corresponds to the angular coordinate
    $\varphi$ and the vertical axis to the radius $r$. }
    \label{fig:fundomout}
\end{figure}

\begin{figure}[ht]
    \centering
    \includegraphics[scale=.65]{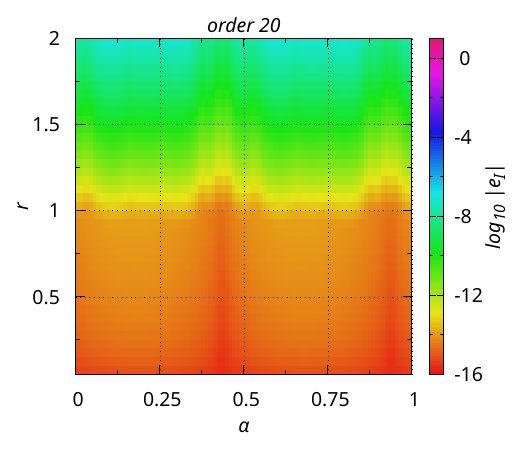}
    \caption{Estimated domain of validity of the rescaled expansion $\graphfun$ computed up to degree $N=20$.
    Colors
    represent $\max_{\theta_0 \in \mathcal G}\log_{10} e_I(s_0,\theta_0,n)$,
    whit $n=1$ and $s_0=re^{\I 2\pi \alpha}$. The horizontal axis corresponds
    to $\alpha$ and the vertical axis to $r$.}
    \label{fig:fundomresc}
\end{figure}

The results discussed provide only information on expansions $\graphfun$. However, for
practical purposes, we are interested in computing the practical domain, which
has been defined in Section~\ref{sec.scaling} as the largest open domain around
the invariant torus where the invariance errors \eqref{eq.inverror} or
\eqref{eq.orbiterror} are below a certain tolerance. Figure~\ref{fig:fundomout}
shows the results for several expansions $\graphfun$ computed up to different Taylor
orders of the evaluation of the error in the invariance equation using $n=1$.
The expansions in these cases have been scaled using the default scaling factor.
Each point in the plots shows the maximum of $e_I(s_0,\theta_0,n)$ over the
Taylor-grid, that is $\theta_0 \in \mathcal G$. The horizontal and vertical axes
represent the initial condition $s_0$, which is taken as $s_0=(re^{\I 2\pi
\alpha},re^{-\I2\pi\alpha})$, according to the complexification convention.
The results show how the practical validity domain increases with the order of the
computed expansions. In Figure~\ref{fig:fundomresc} we show the same quantity
for the rescaled expansion calculated up to degree $20$. As discussed in
Section~\ref{sec.scaling}, the validity domain, where the error is below
$\varepsilon_\tol=10^{-15}$, can be approximated by the unit polydisk.

\begin{figure}[ht]
    \centering
    \begingroup
    \setlength{\tabcolsep}{0pt}
    \begin{tabular}{@{}c@{\qquad}c@{}}
      \includegraphics[width=.45\linewidth]{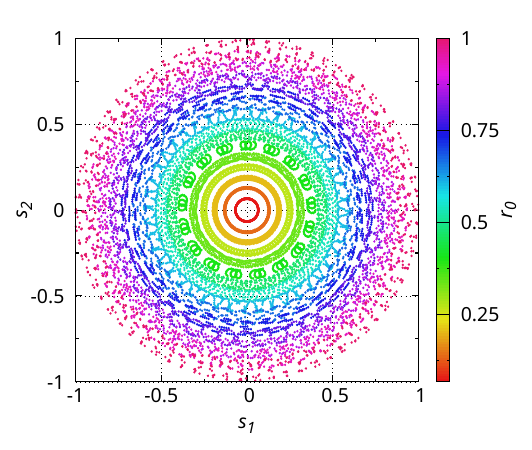} &
      \includegraphics[width=.45\linewidth]{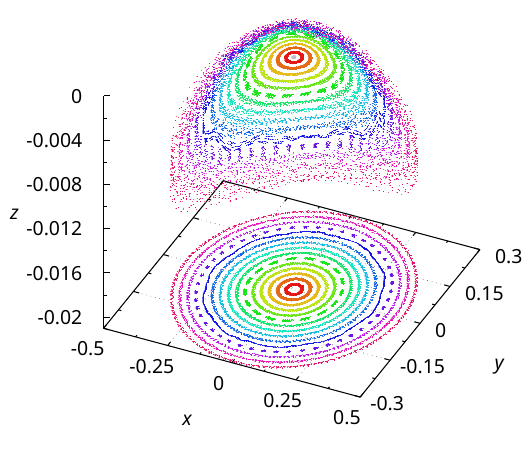}
    \end{tabular}
    \endgroup
    \caption{Iteration of the reduced map \eqref{eq.GT.induced-dynamics} for different initial conditions $s_0$, colored by initial radius. Left: reduced coordinates. Right: their lifting to the original phase space and its $x$-$y$ projection.}
    \label{fig:iterated}
\end{figure}

The estimation of the validity domain allows us to produce trajectories that
remain in the invariant manifold that describes \eqref{eq.GT.K-graph}. Taking
different initial conditions $s_0$, where $\|s_0\|_\infty < 1$, we generate
trajectories using the reduced dynamics given by $\cR _L$. Each of these
trajectories in the reduced coordinates $u$ can then be transformed into the
original coordinates of the system \eqref{eq.FP} using $\param(u,\theta)$.
Figure~\ref{fig:iterated} shows in different colors the different trajectories
obtained for each initial condition by iterating the reduced dynamics. The left panel shows the trajectories in the reduced coordinates, while the right panel shows the same trajectories in the $xyz$ coordinates. Each of these trajectories lies on the center manifold of the response torus, and thus belong to a trajectory on a $3$-dimensional torus, which is equivalent to the periodic orbits that surround the fixed point of the non-perturbed pendulum.

\subsection{An Oscillating Solar Sail around E-M system}\label{sec.solarsail}

Solar Radiation Pressure (SRP) is the pressure exerted by the light
upon a body. The impact of SRP is especially remarkable in those bodies 
which are highly reflecting and have a large area-to-mass ratio (for instance, 
dust \cite{GimenoJJN24}). 

A solar sail is a novel type of low thrust that consists, essentially, in equipping 
a spacecraft
with a highly reflecting surface. This is intended to improve the area-to-mass ratio 
and the reflectivity of the probe in order to maximize the acceleration due to SRP. The 
thrust produced by a sail is limitless (as the source of energy is the Sun) and smooth 
(as long as the shadows are neglected). From the dynamical perspective, a Solar Sail 
changes the geometry of the phase space (see \cite{McInnes99, FarresJ08a, FarresJ08b}). 

\subsubsection{The Augmented Bicircular Problem}
Let us consider a spacecraft moving in the Earth-Moon system. The basic model to 
describe the motion of the probe is the so-called Circular Restricted Three Body Problem (CRTBP), 
which assumes the Earth and the Moon (the primaries) 
to describe circular orbits around their common center 
of masses. The resulting system of ODE is Hamiltonian and has two degrees of freedom when the only 
motion which is considered is the one taking place in the Earth-Moon plane (the only situation 
that will be considered in this work). Usually, the CRTBP is written in synodic coordinates, 
meaning that the frame of reference is rotating so the primaries become fixed. In these 
coordinates, the Hamiltonian has five equilibrium points: Three of them; $L_1$, $L_2$ and $L_3$
are unstable; and two of them $L_4$ and $L_5$ are stable. 

There are several works pointing out the impact of the gravitational pull of the Sun on a massless particle in the Earth-Moon system, specially near the triangular points 
(see \cite{GomezJMS91a, SimoGJM95, JorbaCFJ18} and references therein). Moreover, the impact 
of Solar Radiation Pressure (SRP), the pressure exerted by the photons coming from 
the Sun upon a body, is also remarkable (see \cite{GimenoJJN24}). 
Notice that direct gravitational pull and the SRP act along the Sun-body line. Therefore, both
act as periodic perturbations (of the same period) of the CRTBP. When these periodic perturbations
are considered, the equilibrium points become periodic orbits with the period of the Sun. Those 
periodic orbits are usually called dynamical equivalents. 

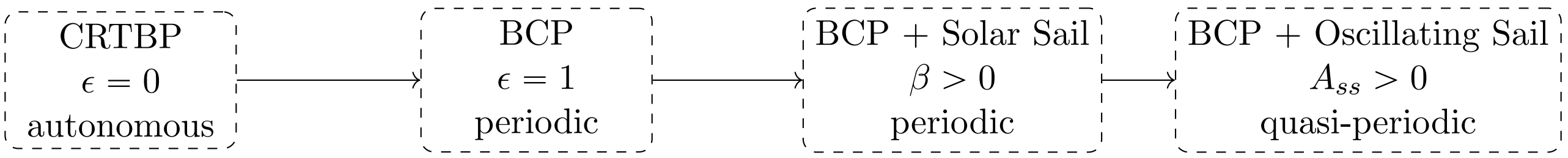
\begin{figure}[ht]
\begin{center}
\begin{tikzpicture}[
    node distance=4.5cm,
    every node/.style={
        draw,
        rectangle,
        rounded corners,
        dashed,
        align=center,
        minimum width=2.5cm,
        minimum height=1cm
    }
]

\node (crtbp) {CRTBP\\$\epsilon = 0$\\autonomous};
\node (bcp) [right of=crtbp] {BCP\\$\epsilon = 1$\\periodic};
\node (bcss) [right of=bcp] {BCP + Solar Sail\\$\beta > 0$\\periodic};
\node (bcposs) [right of=bcss] {BCP + Oscillating Sail\\$A_{ss} > 0$\\quasi-periodic};

\draw[->] (crtbp) -- (bcp);
\draw[->] (bcp) -- (bcss);
\draw[->] (bcss) -- (bcposs);

\end{tikzpicture}
 \end{center}
\caption{Hierarchy of models, including solar perturbations, solar radiation pressure, and oscillating sail orientation.}\label{fig:oss_hie}
\end{figure}

When the spacecraft is endowed with a solar sail, the effect of the SRP is maximized by 
increasing the so-called lightness number $\beta$. In this work we consider an oscillating 
solar sail. Namely, a sail that changes its orientation to the Sun periodically. 
A possible model for the motion of an oscillating sail is the one governed by the following 
set of ODE. 
\begin{align*}
\dot{x} &= v_x, \qquad 
\dot{y} = v_y, \\
\dot{v}_x &= 2 v_y + x - (1-\mu)\frac{x-\mu}{r_E^3} - \mu \frac{x+1-\mu}{r_M^3} 
- \epsilon m_s \frac{x-x_s}{r_S^3} - \epsilon \frac{m_s}{a_s^2}\cos(\omega_s t + \phi_0) + F_{ss,x}, \\
\dot{v}_y &= -2 v_x + y - (1-\mu)\frac{y}{r_E^3} - \mu \frac{y}{r_M^3} 
- \epsilon m_s \frac{y-y_s}{r_S^3} + \epsilon \frac{m_s}{a_s^2} \sin(\omega_s t +  \phi_0) + F_{ss,y}.
\end{align*}
Here, $r_{E}^{2} = (x-\mu)^2 + y^2$, $r_{M}^{2} = (x+1-\mu)^2 + y^2$, $r_{S}^{2} = (x-x_S)^2 + (y-y_S)^2$ are the distances from the sail to the Earth, Moon and Sun respectively. The position 
of the Sun is given by the periodic functions.
$x_S = a_S \cos(\omega_S t + \phi_0)$ and $y_S = -a_S \sin(\omega_S t + \phi_0)$
And the force due to SRP upon the Sail is 
\begin{align*}
F_{ss,x} &= -\beta \gamma \cos^2(\alpha) \left ( \cos(\alpha)\cos(\omega_s t+\phi_0) + \sin(\alpha)\sin(\omega_s t+\phi_0)\right),   \\
F_{ss,y} &= -\beta \gamma \cos^2(\alpha) \left ( \sin(\alpha)\cos(\omega_s t+\phi_0) - \cos(\alpha)\sin(\omega_s t+\phi_0) \right ). 
\end{align*}
Finally, $\alpha \in [-\pi/2, \pi/2]$ the orientation of the sail, depends on time as 
$$\alpha(t)=A_{ss} \frac{\pi}{2} \sin \left ( \omega_{ss}t + \psi_{0} \right) .$$
Notice that, for $A_{ss}=0$, the sail remains perpendicular to the Sun and the force 
exerted by the SRP becomes maximal. On the other hand, if $\alpha = \pm \frac{\pi}{2}$ 
the effect due to SRP is canceled. The amplitude $A_{ss} \in [0, 1]$  acts 
as a control parameter of the oscillation. If $A_{ss}=0$, the sail remains perpendicular 
to the Sun for all time. 
The parameter $\omega_{ss}$ is the frequency of oscillation of the sail. If the frequency 
of the sail is rationally independent of the angular frequency of the Sun, the system becomes 
a quasi-periodic perturbation of the CRTBP and, therefore, the equilibrium points become 
quasi-periodic tori (see \cite{JorbaV97a}). Obviously, the frequency of the sail can be chosen 
to be resonant and this choice would keep the model periodic. Both choices make sense depending 
on the context but this discussion is out of the scope of this paper. Here we will take 
$\omega_{ss}=(\sqrt{5}-1)/2$ i.e. the golden ratio, thus, the model is quasi-periodic for 
non-vanishing amplitude of the sail. Notice that the model can be autonomous, periodic or 
quasi-periodic depending on the choice of the parameters. This information is summarized in 
Figure~\ref{fig:oss_hie}. Moreover, Table~\ref{tab:oss_param} contains a summary of the 
parameters as well as the specific value of the physical constants used in the computations 
within this section. 
The idea of an Oscillating Solar Sail was introduced in \cite{heiligers2017}, where 
the extra frequency is taken to be resonant (so the model is periodic). The quasi-periodic 
case has been studied in \cite{chujo2024}.

\begin{table}[t]
\begin{center}
\begin{tabular}{llr}
\hline
Parameter & Description & Value\\
\hline
$\mu$ & Earth--Moon mass parameter & {\texttt 0.012150581623433623}\\
$m_s$ & mass of the Sun  & {\texttt 328900.54999999906}\\
$\beta$ & sail lightness number  & ~ \\
$\omega_s$ & Sun's angular frequency  & {\texttt 0.92519598551829646}\\
$a_s$ & Sun's mean distance  & {\texttt 388.81114302335106}\\
$\epsilon$ & homotopy parameter for Sun's mass  & ~\\
$\omega_{ss}$ & sail oscillation frequency  & 0.618033988749895\\
$A_{ss}$ & amplitude of sail oscillation  & ~\\
$\phi_0$ & Sun's initial phase  & ~\\
$\psi_0$ & sail's initial phase  & ~\\
\hline
\end{tabular}
\end{center}
\caption{Parameters of the Augmented BCP. The values given in the right column are physical constants that are fixed. The value of the parameters of the Sail depend on each case. The value of $\omega_{ss}$ 
is also fixed. See text for more details.}\label{tab:oss_param}
\end{table}

\subsubsection{The dynamical equivalents of $L_3$}
Hierarchy shown in Figure~\ref{fig:oss_hie} is translated on the nature of the equilibrium 
point $L_3$ of the RTBP. When $\epsilon>0$, the periodic perturbation due to the Sun is activated 
and $L_3$ is replaced by a periodic orbit of period $T_{S}=2\pi/\omega_{S}$. The equilibrium is 
continued as a periodic orbit until $\epsilon=1$ where the mass of the Sun is achieved. This 
periodic orbit is named the dynamical equivalent of $L_3$ in the BCP and it has been extensively 
studied in \cite{JorbaN20}. When $\beta>0$ and $A_{ss}=0$ the SRP perturbation is considered and the 
dynamical equivalent of $L_3$ changes. The perturbation remains periodic with period $T_{S}$ as 
SRP is parallel to solar gravitation but with opposite sense. In fact, the value $\beta=1$ cancels
out direct solar gravity. However, it does not cancel indirect gravitational effect. In 
\cite{JorbaCuscoEA2020, GimenoJJN24} it in shown that the impact of SRP in the Earth-Moon system 
is large even for very small values of $\beta$. We continued the dynamical equivalent 
of $L_3$ until $\beta=0.01$ and $A_{ss}=0$ (the sail remains perpendicular to the Sun line), 
which is, we recall, still a periodic orbit. 

Once the amplitude of the sail oscillation (or swing) $A_{ss}$ is set to be different
from zero, the system becomes a quasi-periodic perturbation of the RTBP (again, see Figure~\ref{fig:oss_hie}) and the dynamical equivalent has a trajectory filling densely 
an invariant torus. Its two basic frequencies are $\omega_{S}$ and $\omega_{ss}$. 
Small values of the swing also have large impact on the dynamical equivalent of $L_3$.
We have selected two different values of $A_{ss}$ to conduct our experiments. 

\begin{table}[]
\centering
\begin{tabular}{llcc}
\toprule
\textbf{$A_{ss}$}           & \textbf{$N$}            & \textbf{$\lambda_{u}$}                       & \textbf{$\lambda_{e}$}                                                   \\ \midrule
\multicolumn{1}{l}{0.025} & \multicolumn{1}{l}{16} & \multicolumn{1}{l}{3.03613889448374952e+00} & \multicolumn{1}{l}{8.12837328855963603e-01 $\pm \I$5.82490752560312797e-01} \\ \midrule
0.2                         & 40                      &      3.07621464827837521e+00                 & 8.18636511293320046e-01 $+\I$5.74311990452311072e-01                      \\ \bottomrule
\end{tabular}
\caption{For $\beta=0.01$ relevant parameters of the two tori. The integer $N$ is the 
truncation order of the torus. $\lambda_u$ is the unstable eigenvalue and $\lambda_{e}$ is 
the elliptic pair. The stable eigenvalue can be recovered from the unstable as $\lambda_{s}=\lambda_{u}^{-1}$.}
\label{tab:L3_NE_stab}
\end{table}

\begin{figure}[tbp]
    \centering

    \begin{tabular}{@{}c@{\hspace{0.01\linewidth}}c@{\hspace{0.01\linewidth}}c@{}}
    \multicolumn{3}{@{}c@{}}{
      \small $\beta=0.01,\quad A_{\mathrm{ss}}=0.025,\quad \deg=16$
    } \\[0.25em]
    \includegraphics[width=0.32\linewidth,page=1]{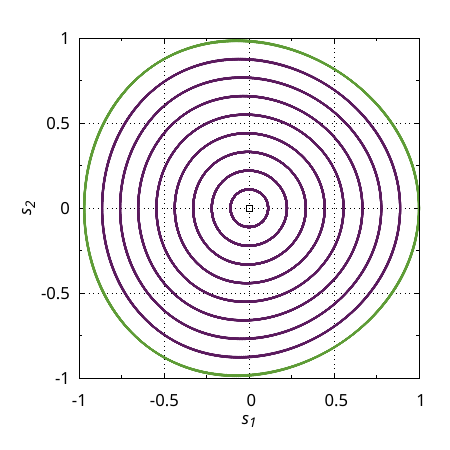} &
    \includegraphics[width=0.32\linewidth,page=2]{L3_OSS.pdf} &
    \includegraphics[width=0.32\linewidth,page=3]{L3_OSS.pdf}
  \end{tabular}
    \caption{Left: Iteration of the reduced map for different initial conditions. The outermost invariant curve is depicted in green. Center: Comparison between the original torus (black) and 
    the outermost invariant torus (green) projected onto the $x$-$y$ space. Right: Projection 
    onto the $v_{x} $-$ v_{y}$ variables.}
    \label{fig:OSS_L3}
\end{figure}

First, we have continued the dynamical equivalent of $L_3$ for $A_{ss}=0.025$. Relevant 
information of this particular torus can be found in the first row of Table~\ref{tab:L3_NE_stab}.
In the stroboscopic map, it can be represented as an invariant curve with internal 
rotation number $\rho=4.1971886370383951$ 
(this is straightforward from $\rho=2\pi \omega_{ss}/\omega_{S}$). It can be represented 
using a trigonometric polynomial of degree $M=16$. It is of stability type saddle $\times$ center, 
being $\lambda_u\approx3.04$. We have applied 
the algorithm in \cite{Jorba01} to compute the torus and its stability. We have applied Algorithm~\ref{palg.graph-tor} to recover the center manifold about this dynamical equivalent, 
up to degree 16. In Figure~\ref{fig:OSS_L3} we display the center manifold and the outermost invariant curve. In the Left panel, inner iterates of the reduced dynamics are shown. 
Notice that this set of concentric curves around the origin are expressed in coordinates of the center manifold $(s_{1}- s_{2})$. The black dot in the center corresponds 
to the original dynamical equivalent. The outermost green curve is the curve at 
initial condition $r\approx 0.9994$ (the domain of validity of the center manifold 
has been scaled to 1). The invariance equation error with a Fourier mesh of 64 
is about $10^{-10}$ in phase space. In the center and right panels we display the 
dynamical equivalent (an invariant curve) in purple and the outermost torus
(a 2-D torus) in green for phase space coordinates of the stroboscopic map. The 
center panel displays the projection onto $x$-$y$ coordinates while the right panel 
shows the projection onto the $v_{x}$-$v_{y}$ variables.  

\begin{figure}[ht]
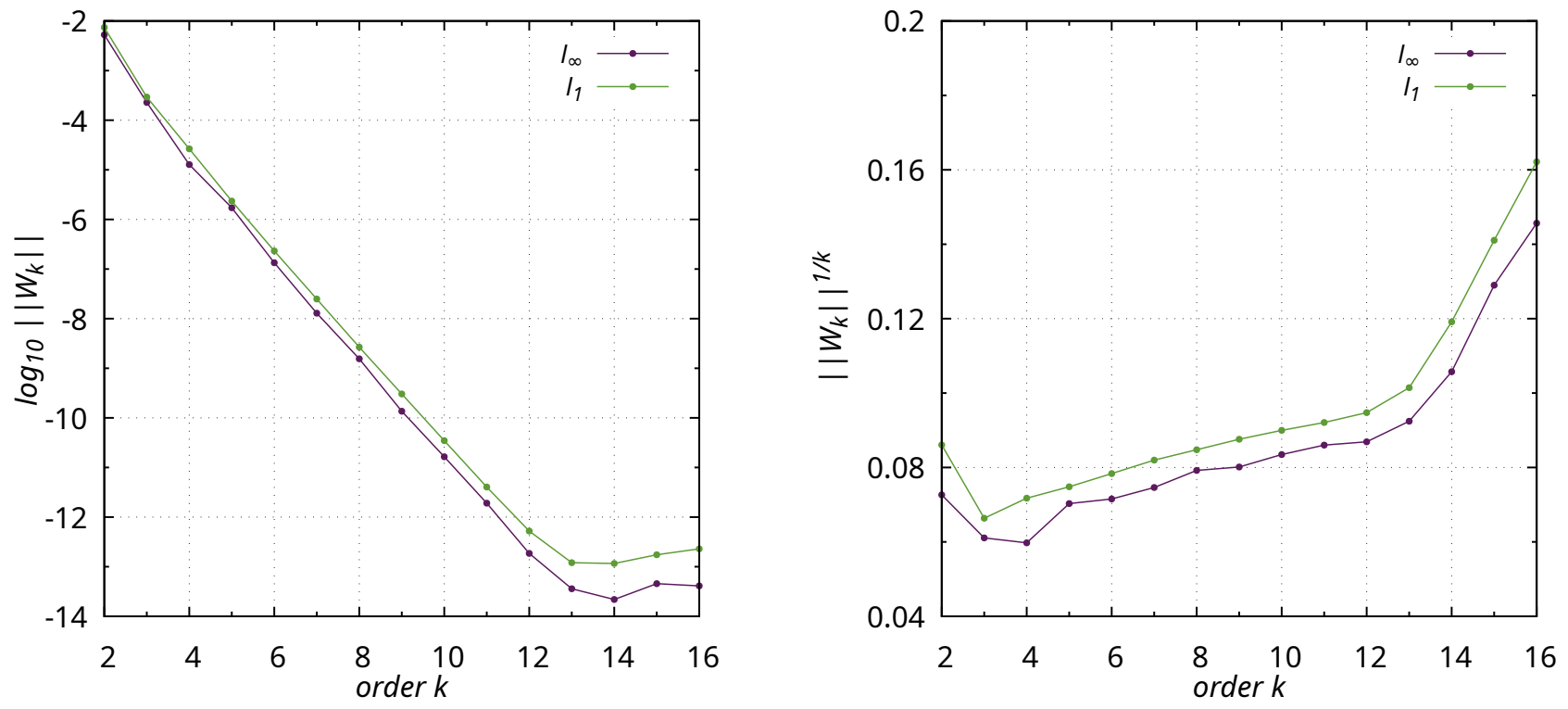

    \centering
\begin{tabular}{@{}c@{\hspace{0.02\linewidth}}c@{}}
    \includegraphics[width=0.32\linewidth,page=4]{L3_OSS.pdf} &
    \includegraphics[width=0.32\linewidth,page=5]{L3_OSS.pdf} 
\end{tabular}
    \caption{Order-wise coefficient norms \eqref{eq.norminf1} in real (left) and their $k$-th roots (right), for the solar-sail center manifold.} 
    \label{fig:OSS_L3_graph_norms}
\end{figure}

For the second conducted experiment we have selected a larger swing ($A_{ss}=0.2$). 
Due to the large impact of SRP on $L_3$ such a swing (despite having little impact 
on the stability of the orbit) deforms the torus approaching it to collision 
with the Earth. As a consequence, the number of Fourier coefficients needed to 
describe the dynamical equivalent of $L_3$ is larger $(M=40)$. This information 
can be read directly from the second row of Table~\ref{tab:L3_NE_stab}. This 
second experiment has the only purpose of comparing Algorithm~\ref{palg.graph-tor} in
its real and complex implementation. In particular, we demonstrate that the real 
implementation is computationally more efficient than the complex one. 

In Table~\ref{tab:OSS_timevsorder} we display the computation of the center manifold 
for degrees from 2 to 10. The first row RGC Time stands for the time 
needed to compute the Graph using the real method (it includes all the steps of
the algorithm), the second row CGC Time represents the time to compute the same 
graph but using complex numbers. The third column represents the ratio of the 
times CGC/RGC representing the speedup. 

\begin{table}[ht]
        \centering
        \small
        \begin{tabular}{r||c|c|c|c||r}
            \textbf{Degree} & 2 & 4 & 8 & 10 & \textbf{Total}\\ \hline
            RGC Time (s) & 0.109684 & 0.270083 & 1.161499 & 2.090302& 7.213675\\
            CGC Time (s) &0.187892& 0.519322 & 2.620454 & 4.984686& 16.170677\\
            CGC/RGC & 1.713 &  1.851 & 2.132 &  2.242& 2.2417\\
        \end{tabular}
        \caption{Computational time (wall-clock) to compute Real (RGC) and complex (CGC) $\graphfun$ using $10$ threads.}
        \label{tab:OSS_timevsorder}
    \hfill
\end{table}

In Table~\ref{tab:OSS_timevsthreads} we display computation times of the real and 
complex method of the center manifold up to degree 16 with respect to the number 
of threads. 
\begin{table}[ht]
        \centering
        \small
        \begin{tabular}{r||c|c|c|c}
            \textbf{Num. of threads} & 1 & 2 & 4 & 8 \\ \hline
            RGC Time (s) & 54.94 & 29.74 & 16.04 & 8.70\\ \hline 
            CGC Time (s) & 125.94 & 67.57 & 35.29 & 19.49
        \end{tabular}
        \caption{Computational time (wall-clock) to compute
        $\graphfun$ up to degree 16 compared to the number of threads.}
        \label{tab:OSS_timevsthreads}
\end{table}

Most of the time saved by the real method comes from the jet transport step. The 
real approach has an extra advantage: memory saving. Indeed, as 
complex Taylor coefficients storage is approximately twice the real ones.

\subsection{Simplified Solar System Model}\label{sec.sssm}
The Circular Restricted Three Body Problem (CRTBP) is a model for the motion of a test-particle 
that is gravitationally pulled by two massive bodies; the so called primary masses. In the model, 
the test-particle (sometimes called massless particle or probe) has a negligible mass so it does 
not affect the motion of the primaries, which move following a circle about their common center 
of mass. Notice that any other Keplerian motion is also allowed (i.e. ellipses, hyperbolas and 
parabolas), leading to more involved models. The book \cite{Szebehely67} is a standard reference 
that deals mostly with the circular problem and treats the elliptic as an extension. See also 
\cite{faintich1972, Barrabes2015}, as examples for the hyperbolic and parabolic cases. 

The CRTBP is a foundational model for most of the problems in celestial mechanics and astrodynamics. 
However, the presence of other masses does affect long-term behavior of the probe. And, while many-body 
simulations are available and are always used in the last stages of a project, intermediate 
models are way more useful to study the phase space of the problem: those are cheaper to evaluate 
and simple to study from the theoretical perspective. The most widely used model of this kind 
is the Bicircular Problem (BCP), a model that imposes the motion of a third primary which moves along 
a circular orbit together with center of masses of the first two primaries, see \cite{Huang60, SimoGJM95}. 
The BCP can be written as a periodic perturbation of the CRTBP, whose primary effect on 
the phase space of the CRTBP it that every quasi-periodic orbit gains a frequency: 
the equilibria become periodic orbits, the periodic orbits become quasi-periodic tori with two frequencies, etc.

G\'omez, Masdemont and Mondelo developed a methodology to produce quasi-periodic models with many
frequencies \cite{Mondelo01,GomezMM02}. They call these models Simplified Solar System Models (SSSM). 
These Quasi-periodic models, therefore, describe the motion of a test particle subjected to a quasi-periodic 
force resulting from several planets and other effects such as eccentricities, inclinations, etc. 
It can be written as time-dependent perturbations of the CRTBP. 

A case of special relevance is the Earth-Moon example (Section 4.2 in \cite{GomezMM02}). 
Here, the Earth-Moon CRTBP is established as a basis and the Sun is included as an external 
perturbation (such as in the case of the BCP). However, the effect of the Sun on the motion 
of the Earth and the Moon is considered here (this is not the case of the BCP) and 
their orbital elements move periodically according to the excitation due to the Sun. 
That is, besides the direct gravitational effect on the probe, the effect of the 
lunar eccentricity; the lunar inclination and the inclination between the 
orbital and equatorial planes are considered.

This leads to a model with five frequencies, namely:
\begin{enumerate}
    \item The mean longitude of the Moon, $\omega_1$, 
    \item the mean elongation of the Moon from the Sun, $\omega_2$ (this is 
    the same external frequency as in the BCP), 
    \item the mean longitude of the lunar perigee, $\omega_3$, 
    \item the longitude of the mean ascending node of the lunar orbits on the ecliptic, $\omega_4$, 
    \item and the Sun's mean longitude of perigee $\omega_5$.
\end{enumerate}

The equations of motion for the test particle in this Earth-Moon $SSSM$ are written as 
\begin{equation}\label{eq:em-sssm}
  \begin{cases}
        \ddot{x}=c_{1}^{i} + c_{4}^{i} \dot{x} + c_{5}^{i} \dot{y} + c_{7}^{i} x + c_{8}^{i} y + c_{9}^{i}z + c^{i}_{13} \frac{\partial \Omega^{i}}{\partial x},\\

        \ddot{y}=c_{2}^{i} - c_{5}^{i} \dot{x} + c_{4}^{i} \dot{y} + c_{6}^{i}\dot{z} -c_{8}^{i} x + c_{10}^{i} y + c_{11}^{i}z + c^{i}_{13} \frac{\partial \Omega^{i}}{\partial y},\\
        
        \ddot{z}=c_{3}^{i} - c_{6}^{i} \dot{y} + c_{4}^{i} \dot{z} + c_{9}^{i} x - c_{11}^{i} y + c_{12}^{i}z + c_{13}^{i} \frac{\partial \Omega^{i}}{\partial z}.

    \end{cases}
\end{equation}
being, 
\begin{align*}
    \Omega^{i}= \frac{1-\mu}{\sqrt{(x-\mu)^2 + y^2 + z^2}} &+ \frac{\mu}{\sqrt{(x-\mu+1)^2 + y^2 + z^2}}+ \frac{\mu_{S}}{\sqrt{(x-x^{i}_{S})^2 + (y-y^{i}_{S})^{2} + (z-z_{S}^{i})^2}}.
\end{align*}
Here, $\mu$ is the Earth-Moon mass parameter and $\mu_S$ is the mass of the Sun
as multiple of the total mass of the Earth-Moon system. The position of the Sun
is given by $(x_{S}^{i}, y_{S}^{i}, z_{S}^{i})$, and the coefficients
$c_{j}^{i}$, $i=1,\dots, 5$ and $j=1,\dots,13$ are quasi-periodic functions with
vector of basic frequencies $\nu=(\nu_1 , \dots, \nu_5)$, with $\nu_1 = \omega_1
$, $\nu_2= \omega_1 - \omega_3$, $\nu_3=\omega_1 - \omega_2 + \omega_4$, $\nu_4
= \omega_1 - \omega_5$ and $\nu_5 = \omega_5 -\omega_2$. We denote this model as
$SSSM_i$, where $i$ denotes how many of the perturbations are considered: $\nu_1, \dots, \nu_i$.
Therefore, the functions $c^{i}_{j}$ as well as the position of the Sun change
when different perturbations are considered. 

It is well known that the CRTBP has five equilibrium points: three of them are
collinear (aligned with the primaries), those are $L_{1,2,3}$. The remaining two
$L_{4,5}$ are the triangular points, which lie in the third vertex of an
equilateral triangle whose remaining two vertices are the primaries. The
collinear points are, for all the values of the masses, unstable and of type
center $\times$ center $\times$ saddle. The triangular points are center
$\times$ center $\times$ center if $\mu$ is below a critical value. When
periodic perturbations are activated, the equilibria are replaced by invariant
tori (with the same frequencies as the perturbations) and, if the perturbations
are small enough, have the same stability as the fixed points (see
\cite{JorbaV97a, JorbaV97b}).

As is well known, the collinear libration points $L_{1,2}$ exhibit strong instability.
To mitigate its effects, we adopt a
multiple-shooting strategy that splits the integration time of the
directed return map $P_{\sigma}$ into $\Nseg$ segments and replaces the single
invariant torus with a cyclic family $\{\tor_i\}_{i \in \Z_{\Nseg}}$. Accordingly,
computing the invariant manifold attached to the single torus is
replaced by computing a cyclic family of parameterizations
$\{W_i\}_{i \in \Z_{\Nseg}}$ together with the corresponding reduced dynamics
$\{f_{\sigma,i}\}_{i \in \Z_{\Nseg}}$ (see
Appendix~\ref{app.multiple-shooting}). 

\begin{figure}[ht]
    \centering
    \includegraphics[width=0.32\linewidth]{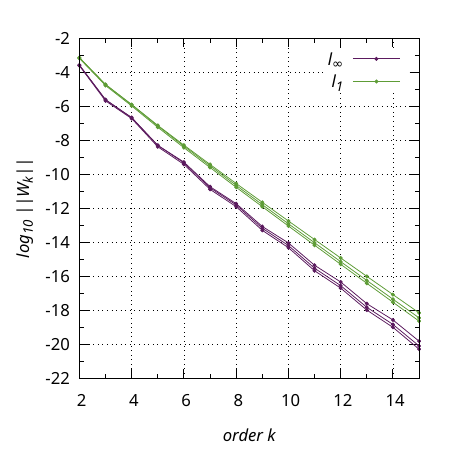}
    \caption{Norms of each
    parameterization section of the center manifold attached to the equivalent of $L_2$ of the $SSSM_2$ as in
    Equation~\eqref{eq.norminf1}.}
    \label{fig:SSSM_L2_graph_norms}
\end{figure}

We consider the $SSSM_2$ and the dynamical equivalent to $L_2$, which is a
two-dimensional invariant torus for the flow. This torus presents a Floquet
multiplier of order $10^{6}$. Since we want to compute its associated center
manifold ($n_L=4$), we choose to use the multiple-shooting version of
Algorithm~\ref{palg.graph-tor} with $\Nseg=4$. The original one-dimensional
invariant torus for the Poincar\'e map is represented on a grid of $M = 161$
points in the angular variables. The parameterizations are computed up to
several degrees using a multi-threaded implementation.
Table~\ref{tab:SSSM_timevsorder} reports the corresponding wall-clock times for
degrees up to 15. Table~\ref{fig:SSSM_L2_graph_norms} presents the norms
\eqref{eq.norminf1} across all orders of the $\Nseg$-segment parameterizations.

\begin{table}[ht]
        \centering
        \small
        \begin{tabular}{r||c|c|c|c}
            {\bf Degree} & 6 & 8 & 10 & 15 \\ \hline
            Time (s) & 459.25 & 791.90 & 1179.2 & 24850  \\ \hline
            Practical domain (km) & 1338 & 2141 & 2775 & 3667 
        \end{tabular}
        \caption{Computational time (wall-clock) to compute $\graphfun$ with $n_L=4$ attached to $L_2$ of the $SSSM_2$ with $\Nseg=4$
        up to several degrees using multiple-shooting and a complex-valued implementation parallelized over $16$ threads and their associated practical domain of validity reported in kilometers from the dynamical equivalent of $L_2$.
        }
        \label{tab:SSSM_timevsorder}
    \hfill
\end{table}

Finally, to assess the accuracy of the computed parameterizations, we repeat the
practical domain estimation described in Section~\ref{sec.scaling}.
Specifically, we use the error along the reduced orbit defined in
\eqref{eq.inverror}, with $n_{\rm iter}=1$ and tolerance
$\varepsilon_{\tol}=10^{-10}$. Rather than reporting the maximum distance $R$
from the origin in manifold coordinates, we report the maximum distance, in
kilometers, from the dynamical equivalent of $L_2$, as defined in
\eqref{eq.GT.K-graph}, for $\|s\|_\infty<R$. For the parameterization of
degree-$15$, the practical domain extends up to almost $3700$ km from the
original invariant torus.

\section*{Statements and Declarations}
The authors declare that they have no conflict of interest.

\section*{Acknowledgments}


The project has been supported by the Spanish grants
PID2021-125535NB-I00, PID2024-158570NB-100, and PID2025-174015NB-I00 (MICINN/AEI/FEDER, UE),
the Catalan grant 2025 SGR 00322, and by the Air Force Office of Scientific
Research under award number FA8655-24-1-7059.  The project that led to these results also received
the support of a fellowship from ``la Caixa'' Foundation (ID
100010434), the fellowship code is LCF/BQ/PR23/11980047.  This work
has been also funded through the Severo Ochoa and Mar\'ia de Maeztu
Program for Centers and Units of Excellence in R\&D
(CEX2020-001084-M).
\appendix

\bibliographystyle{alpha}
\addcontentsline{toc}{section}{References}
{\footnotesize \bibliography{ds,na,misc,pre}}

\section{Multiple Shooting for the Graph-parameterization}
\label{app.multiple-shooting}

The graph-parameterization of Section~\ref{sec.graph} is written in terms of a
single directed return map $\poincmap_\sigma$.  When this map is obtained by
integrating a strongly hyperbolic flow over a complete return time, however,
the evaluation of the map and, more importantly, the transport of
high-order jets over the whole return interval can lose numerical accuracy.
Multiple shooting replaces this long propagation by a cyclic sequence of
shorter flow maps.  It is therefore a numerical factorization of the same
invariance problem, rather than a different invariant-manifold construction.
The method is classical for boundary-value problems
\cite{AscherMR95}; multiple-shooting formulations for invariant tori and their
bundles were developed in \cite{HaroM2021}, and high-order invariant manifolds
of reducible quasi-periodic tori were treated in \cite{GimenoJNO21}.  A recent
multiple-shooting parameterization method for invariant manifolds of periodic
orbits is given in \cite{Kumar2026}.

\smallskip

We use the directed notation of Section~\ref{sec.inv-eq}.  Let
$N_{\mathrm{seg}}\geq 2$ be the number of shooting segments and set
\begin{equation}
    \label{eq.MS.delta}
    \delta_\sigma \bydef \frac{\rho_\sigma}{N_{\mathrm{seg}}}.
\end{equation}
The shooting sections are numbered in the direction of the dynamics selected
by $\sigma$.  Thus, for $\sigma=+1$ they are traversed in forward time, while
for $\sigma=-1$ the order is reversed and the segment maps are the
corresponding local inverse flow maps, see Appendix~\ref{app.stable}.  Let
\begin{equation*}
    \label{eq.MS.segment-skew-product}
    \expoincmap_{\sigma,i}(x,\theta)
    \bydef
    \bigl(\poincmap_{\sigma,i}(x,\theta),
          \theta+\delta_\sigma\bigr), 
    \qquad 
    \text{for }i\in\Z_{N_{\mathrm{seg}}} \bydef\Z/{N_{\mathrm{seg}}}\Z,
\end{equation*}
be the directed skew-product associated with the $i$th segment.  The indexing
is chosen so that one complete cycle recovers the directed return map,
\begin{equation}
    \label{eq.MS.factorization}
    \expoincmap_\sigma
    =
    \expoincmap_{\sigma,N_{\mathrm{seg}}-1}
    \circ\cdots\circ
    \expoincmap_{\sigma,0}.
\end{equation}
The angular arguments in \eqref{eq.MS.factorization} are advanced
successively by \eqref{eq.MS.delta}; hence the accumulated rotation after one
cycle is $\rho_\sigma$.

The invariant torus is represented on the shooting sections by a cyclic
family $\{\tor_i\}_{i\in\Z_{N_{\mathrm{seg}}}}$ satisfying
\begin{equation*}
    \label{eq.MS.torus-chain}
    \poincmap_{\sigma,i}\bigl(\tor_i(\theta),\theta\bigr)
    =
    \tor_{i+1}\bigl(\theta+\delta_\sigma\bigr),
    \qquad i\in\Z_{N_{\mathrm{seg}}}.
\end{equation*}
Define
$
    A_{\sigma,i}(\theta)
    \bydef
    \dop_x\poincmap_{\sigma,i}\bigl(\tor_i(\theta),\theta\bigr)
$.
In analogy with \eqref{eq.GT.directed-adapted-frame}, assume that along each
segment there is an adapted frame
$
  \adframe_i(\theta)
    =
    \begin{bmatrix}
       \adframe_{L,i}(\theta) & \adframe_{N,i}(\theta)
    \end{bmatrix}
$
with the same dimensions $\nt$ and $\nh$ as in
\eqref{eq.GT.LN-splitting}, and constant segment matrices
\begin{equation}
    \label{eq.MS.frame-invariance}
    \Lambda_{\sigma,i}
    =
    \begin{pmatrix}
      \Lambda_{L,\sigma,i} & 0\\
      0 & \Lambda_{N,\sigma,i}
    \end{pmatrix},
\qquad \text{such that} \quad
    A_{\sigma,i}(\theta)\adframe_i(\theta)
    =
    \adframe_{i+1}(\theta+\delta_\sigma)
    \Lambda_{\sigma,i}.
\end{equation}
We take $\adframe_0$ to agree with the frame used in
Section~\ref{sec.adapt-frame}.  Composing \eqref{eq.MS.frame-invariance} over
one complete cycle and comparing with
\eqref{eq.GT.directed-adapted-frame} gives
\begin{equation}
    \label{eq.MS.full-cycle-Lambda}
    \Lambda_\sigma
    =
    \Lambda_{\sigma,N_{\mathrm{seg}}-1}
    \cdots
    \Lambda_{\sigma,1}\Lambda_{\sigma,0},
\end{equation}
and therefore, blockwise,
\begin{equation}
    \Lambda_{L,\sigma}
    =
    \Lambda_{L,\sigma,N_{\mathrm{seg}}-1}\cdots
    \Lambda_{L,\sigma,0},
    \qquad
    \Lambda_{N,\sigma}
    =
    \Lambda_{N,\sigma,N_{\mathrm{seg}}-1}\cdots
    \Lambda_{N,\sigma,0}.
    \label{eq.MS.full-cycle-blocks}
\end{equation}
Notice that the factors in \eqref{eq.MS.full-cycle-Lambda} need not be equal
and, in general, their individual eigenvalues do not need to be 
$N_{\mathrm{seg}}$-th roots of the full-return multipliers.

\begin{rmk}[Root multipliers in the lifted formulation]
\label{rmk:MS.roots}
A root relation does hold for the lifted block-cyclic formulation of the
multiple-shooting cocycle.  More precisely, if $\widetilde\Lambda_\sigma$ is
the $nN_{\mathrm{seg}}\times nN_{\mathrm{seg}}$ block-cyclic matrix obtained
from the segment matrices $\Lambda_{\sigma,i}$, then
\[
    \mu\in\Spec(\widetilde\Lambda_\sigma)
    \quad\Longleftrightarrow\quad
    \mu^{N_{\mathrm{seg}}}
    \in
    \Spec\!\left(
       \Lambda_{\sigma,N_{\mathrm{seg}}-1}\cdots\Lambda_{\sigma,0}
    \right).
\]
This is the precise sense in which multiple-shooting multipliers are roots of
the single-return multipliers; see \cite[Sec.~4.1]{GimenoJNO21}.  A further
rescaling of the segment frames may sometimes be used to distribute a given
multiplier uniformly among the shooting intervals; such a balancing is an
additional normalization, not a consequence of
\eqref{eq.MS.frame-invariance}; see also
\cite[Appendix~A]{Kumar2026}.
\end{rmk}

For each segment, we use exactly the graph normalization of
Section~\ref{sec.norma-cond}.  Namely,
\begin{equation}
    \param_i(\st,\theta)
    =
    \tor_i(\theta)
    +\adframe_{L,i}(\theta)\st
    +\adframe_{N,i}(\theta)\graphfun_i(\st,\theta),
    \label{eq.MS.K-graph}
\end{equation}
where
$
    \graphfun_i(0,\theta)=0$,
 and 
$    \dop_\st\graphfun_i(0,\theta)=0$.
The segment analogue of the transformed map
\eqref{eq.GT.reduced-map} is
\begin{equation}
    \cR_{\sigma,i}(\st,\sh,\theta)
    \bydef
    \adframe_{i+1}(\theta+\delta_\sigma)^{-1}
    \Bigl[
      \poincmap_{\sigma,i}\bigl(
          \tor_i(\theta)
          +\adframe_{L,i}(\theta)\st
          +\adframe_{N,i}(\theta)\sh,
          \theta
      \bigr)
      -\tor_{i+1}(\theta+\delta_\sigma)
    \Bigr].
    \label{eq.MS.reduced-map}
\end{equation}
Writing
$\cR_{\sigma,i}=(\cR_{L,\sigma,i},\cR_{N,\sigma,i})$, the induced segment
internal dynamics is
\begin{equation*}
    \intdyn_{\sigma,i}(\st,\theta)
    \bydef
    \cR_{L,\sigma,i}
    \bigl(\st,\graphfun_i(\st,\theta),\theta\bigr),
    \label{eq.MS.induced-dynamics}
\end{equation*}
and the graph equations are
\begin{equation}
    \cR_{N,\sigma,i}
    \bigl(\st,\graphfun_i(\st,\theta),\theta\bigr)
    =
    \graphfun_{i+1}\bigl(
       \intdyn_{\sigma,i}(\st,\theta),
       \theta+\delta_\sigma
    \bigr),
    \qquad i\in\Z_{N_{\mathrm{seg}}}.
    \label{eq.MS.graph-invariance}
\end{equation}
Thus \eqref{eq.MS.graph-invariance} is simply the cyclic, segment-wise version
of \eqref{eq.GT.graph-invariance-abstract}.

\begin{pro}[Equivalence with the directed single-shooting equation]
\label{pro:MS.equivalence}
Assume that \eqref{eq.MS.graph-invariance} holds for all
$i\in\Z_{N_{\mathrm{seg}}}$.  Starting from
$(\st_0,\theta_0)=(\st,\theta)$, define recursively
\begin{equation}
    \st_{i+1}
    =\intdyn_{\sigma,i}(\st_i,\theta_i),
    \qquad
    \theta_{i+1}=\theta_i+\delta_\sigma,
    \qquad i=0,\ldots,N_{\mathrm{seg}}-1.
    \label{eq.MS.reduced-composition}
\end{equation}
Then $\param_0$ satisfies the directed parameterization equation
\eqref{eq.GT.directed-PM} for the full return map, with full-return internal
dynamics
\begin{equation}
    \intdyn_\sigma^{[N_{\mathrm{seg}}]}(\st,\theta)
    \bydef \st_{N_{\mathrm{seg}}}.
    \label{eq.MS.full-reduced-map}
\end{equation}
In particular, multiple shooting parameterizes the same local invariant
manifold as the corresponding single-shooting problem.
\end{pro}

\begin{proof}
By \eqref{eq.MS.graph-invariance} and the definition
\eqref{eq.MS.K-graph}, each segment satisfies
$
    \poincmap_{\sigma,i}
    \bigl(\param_i(\st_i,\theta_i),\theta_i\bigr)
    =
    \param_{i+1}(\st_{i+1},\theta_{i+1})
$.
Composing these identities for
$i=0,\ldots,N_{\mathrm{seg}}-1$, using the cyclic convention
$\param_{N_{\mathrm{seg}}}=\param_0$, and applying
\eqref{eq.MS.factorization} gives
\[
    \poincmap_\sigma\bigl(\param_0(\st,\theta),\theta\bigr)
    =
    \param_0\bigl(
       \st_{N_{\mathrm{seg}}},
       \theta+N_{\mathrm{seg}}\delta_\sigma
    \bigr).
\]
Since $N_{\mathrm{seg}}\delta_\sigma=\rho_\sigma$, this is precisely
\eqref{eq.GT.directed-PM} with the internal map
\eqref{eq.MS.full-reduced-map}.
\end{proof}

Proposition~\ref{pro:MS.equivalence} separates the mathematical and numerical
roles of multiple shooting.  Mathematically, the complete cycle is the same
return map and the same invariant graph.  Numerically, the nonlinear map and
its jets are propagated only over one segment at a time.  Hence a large
expansion or contraction accumulated over a full return is distributed among
shorter integrations.  

\smallskip

At the level of the
cohomological equations.  Expand each graph as in \eqref{eq.GT.W-FPS},
\begin{equation*}
    \label{eq.MS.W-FPS}
    \graphfun_i(\st,\theta)
    =
    \sum_{l\geq2}\graphfun_{i,(l)}(\st,\theta),
    \qquad
    \graphfun_{i,(l)}(\st,\theta)
    =
    \sum_{|\alpha|=l}w_{i,l,\alpha}(\theta)\st^\alpha.
\end{equation*}
Assume, as in the diagonal solver of Lemma~\ref{lem.cohom}, that the selected
segment block is diagonal,
\[
    \Lambda_{L,\sigma,i}
    =\diag(\lambda_{\sigma,i,1},\ldots,
           \lambda_{\sigma,i,\nt}),
\]
and, for a multiindex $\alpha\in\N^{\nt}$, set
\begin{equation*}
    \label{eq.MS.beta}
    \beta_{\sigma,i,\alpha}
    \bydef
    \prod_{q=1}^{\nt}
       \lambda_{\sigma,i,q}^{\alpha_q}.
\end{equation*}
At homogeneous degree $l$, the segment-wise counterpart of
Step~\ref{palg.graph-tor.cohom} has the form
\begin{equation}
    \label{eq.MS.cohom-segment}
    \Lambda_{N,\sigma,i}w_{i,l,\alpha}(\theta)
    -
    w_{i+1,l,\alpha}(\theta+\delta_\sigma)
    \beta_{\sigma,i,\alpha}
    =
    z_{i,l,\alpha}(\theta),
    \qquad |\alpha|=l.
\end{equation}
The functions $z_{i,l,\alpha}$ are known from the lower Taylor orders, exactly
as in the single-shooting recursion.

Define
$
    \mathcal M_{N,\sigma,i}
    \bydef
    \Lambda_{N,\sigma,i+N_{\mathrm{seg}}-1}
    \cdots
    \Lambda_{N,\sigma,i+1}
    \Lambda_{N,\sigma,i},
$ and $
    \beta_{\sigma,i,\alpha}^{[N_{\mathrm{seg}}]}
    \bydef
    \prod_{q=0}^{N_{\mathrm{seg}}-1}
       \beta_{\sigma,i+q,\alpha}
$, with segment indices taken modulo $N_{\mathrm{seg}}$.

\begin{pro}[Full-cycle reduction of the cohomological equations]
\label{pro:MS.cohom-cycle}
For fixed $i$, $l$, and $\alpha$ with $|\alpha|=l$, the coupled system
\eqref{eq.MS.cohom-segment} implies the full-cycle equation
\begin{equation}
    \label{eq.MS.cohom-cycle}
    \mathcal M_{N,\sigma,i}w_{i,l,\alpha}(\theta)
    -
    \beta_{\sigma,i,\alpha}^{[N_{\mathrm{seg}}]}
    w_{i,l,\alpha}(\theta+\rho_\sigma)
    =
    Z^{[N_{\mathrm{seg}}]}_{i,l,\alpha}(\theta),
\end{equation}
where
\begin{equation}
    \label{eq.MS.effective-rhs}
\begin{split}
    Z^{[N_{\mathrm{seg}}]}_{i,l,\alpha}(\theta)
    \bydef
    \sum_{q=0}^{N_{\mathrm{seg}}-1}
    &\left(
       \prod_{p=0}^{q-1}\beta_{\sigma,i+p,\alpha}
     \right)
     \left(
       \Lambda_{N,\sigma,i+N_{\mathrm{seg}}-1}
       \cdots
       \Lambda_{N,\sigma,i+q+1}
     \right)\\
    &\hspace{8mm}\times
     z_{i+q,l,\alpha}
     \bigl(\theta+q\delta_\sigma\bigr).
\end{split}
\end{equation}
Empty products in \eqref{eq.MS.effective-rhs} are understood as the identity.
If $i=0$ and the frame is normalized as in
\eqref{eq.MS.full-cycle-blocks}, then the linear operator in
\eqref{eq.MS.cohom-cycle} is exactly the full-return operator appearing in the
single-shooting cohomological equation.
\end{pro}

\begin{proof}
For $q=0,\ldots,N_{\mathrm{seg}}-1$, take
\eqref{eq.MS.cohom-segment} with index $i+q$ and angle
$\theta+q\delta_\sigma$.  Multiply it from the left by
\[
    \left(\prod_{p=0}^{q-1}\beta_{\sigma,i+p,\alpha}\right)
    \left(
      \Lambda_{N,\sigma,i+N_{\mathrm{seg}}-1}
      \cdots
      \Lambda_{N,\sigma,i+q+1}
    \right).
\]
When the resulting $N_{\mathrm{seg}}$ equations are added, every intermediate
term containing
$w_{i+q,l,\alpha}(\theta+q\delta_\sigma)$ cancels with the corresponding term
from the preceding equation.  The only terms that remain on the left are the
first term of the $q=0$ equation and the transported term of the
$q=N_{\mathrm{seg}}-1$ equation.  Using cyclicity of the segment index and
$N_{\mathrm{seg}}\delta_\sigma=\rho_\sigma$ gives exactly
\eqref{eq.MS.cohom-cycle}.  The uncancelled right-hand sides give
\eqref{eq.MS.effective-rhs}.
\end{proof}

Taking Fourier coefficients of \eqref{eq.MS.cohom-cycle} yields
\begin{equation*}
    \label{eq.MS.cohom-cycle-fourier}
    \left[
      \mathcal M_{N,\sigma,i}
      -
      e^{2\pi\I\innerp{k}{\rho_\sigma}}
      \beta_{\sigma,i,\alpha}^{[N_{\mathrm{seg}}]}
      \idmat_{\nh}
    \right]
    \widehat w_{i,l,\alpha,k}
    =
    \widehat Z^{[N_{\mathrm{seg}}]}_{i,l,\alpha,k}.
\end{equation*}
Consequently, the exact resonance condition is
\begin{equation}
    \label{eq.MS.resonance}
    e^{2\pi\I\innerp{k}{\rho_\sigma}}
    \beta_{\sigma,i,\alpha}^{[N_{\mathrm{seg}}]}
    \notin
    \Spec\bigl(\mathcal M_{N,\sigma,i}\bigr).
\end{equation}
For $i=0$, \eqref{eq.MS.full-cycle-blocks} identifies
$\mathcal M_{N,\sigma,0}$ with $\Lambda_{N,\sigma}$ and the product of the
selected segment multipliers with the corresponding full-return monomial
multiplier.  Therefore \eqref{eq.MS.resonance} is precisely the same
full-return non-resonance condition as in Section~\ref{sec.solve.cohom}.
Thus, multiple shooting does not remove a resonance.

In practice, Algorithm~\ref{palg.graph-tor} therefore requires only three
structural changes.  The jet-transport step
\ref{palg.graph-tor.jet} is performed segment by segment using
\eqref{eq.MS.reduced-map}; the cohomological solve
\ref{palg.graph-tor.cohom} is replaced by the cyclic system
\eqref{eq.MS.cohom-segment} (or, equivalently, by its full-cycle reduction
\eqref{eq.MS.cohom-cycle}); and the rotation in
Step~\ref{palg.graph-tor.rot} is $\delta_\sigma$ rather than $\rho_\sigma$ on
each segment.  All Taylor degrees, graph normalization conditions, real/complex
conventions, and a posteriori tests remain the same as in the single-shooting
algorithm.

\section{Backward Skew-Product for Stable Directions}
\label{app.stable}

Section~\ref{sec.inv-eq} introduces the orientation parameter
$\sigma\in\{+1,-1\}$ so that the forward and inverse skew-products are
handled by the same graph formulation.  Stable and unstable parameterizations under map inversion are standard in the parameterization-method; for invariant tori and their
whiskers see \cite{HaroLl06b}.  

Reversing the orientation does not define a
different local stable manifold and does not eliminate resonances from the
cohomological equations.  Its main effect is numerical: it reverses the
relative amplification of the selected stable directions and transverse
errors during the evaluation of the Poincar\'e map and its jets.  This is
the motivation for using inverse Poincar\'e maps in high-order stable
manifold computations; see, e.g., \cite{GimenoJNO21}.

\smallskip

Let $\torinv\bydef\{(\tor(\theta),\theta):\theta\in\T^d\}$.  Since the
directed map with $\sigma=-1$ is, by construction, the local inverse of the
one with $\sigma=+1$, one immediately has
\begin{equation}
\label{eq.stable-backward-set-equivalence}
    \mathcal W^s_{\mathrm{loc}}(\torinv;\expoincmap_{+1})
    =
    \mathcal W^u_{\mathrm{loc}}(\torinv;\expoincmap_{-1}).
\end{equation}
Indeed, the positive orbit of a point under $\expoincmap_{+1}$ is the same
sequence of points as the negative orbit under $\expoincmap_{-1}$.  Thus
$\sigma=-1$ changes the time orientation used to compute the manifold, not
the geometric invariant set obtained in the original phase space.  The
same observation is reflected in
\eqref{eq.GT.directed-reducibility}--\eqref{eq.GT.directed-blocks}: the
adapted bundles are unchanged, whereas the Floquet multipliers are replaced
by their reciprocals.

\smallskip

The numerical distinction is already visible at the linear level.

\begin{pro}[Relative amplification of transverse errors]
\label{pro:stable-backward-conditioning}
Consider a stable Floquet eigenmode with multiplier $\lambda_s$ and a
transverse eigenmode with multiplier $\mu$, where
$0<|\lambda_s|<1$.  In the reducible coordinates of
\eqref{eq.GT.directed-blocks}, write
\[
    v=a v_s+b v_\mu, \qquad a\neq0,
\]
and set $r(v)\bydef |b|/|a|$.  Then one forward and one backward linearized
return satisfy, respectively,
\begin{equation}
\label{eq.stable-relative-amplification}
    r(\Lambda v)
    =\frac{|\mu|}{|\lambda_s|}\,r(v),
    \qquad
    r(\Lambda^{-1}v)
    =\frac{|\lambda_s|}{|\mu|}\,r(v).
\end{equation}
Hence, if the transverse mode is center or unstable, $|\mu|\geq1$, forward
propagation amplifies its contribution relative to the selected stable
component, whereas backward propagation damps it.
\end{pro}

\begin{proof}
Under the forward cocycle the two eigencomponents are multiplied by
$\lambda_s$ and $\mu$, so their amplitude ratio is multiplied by
$|\mu|/|\lambda_s|$.  By \eqref{eq.GT.directed-reducibility}, the inverse
cocycle has the same eigenvectors with reciprocal multipliers.  The
corresponding ratio is therefore
$|\mu|^{-1}/|\lambda_s|^{-1}=|\lambda_s|/|\mu|$, proving
\eqref{eq.stable-relative-amplification}.
\end{proof}

This elementary calculation gives the mechanism relevant to jet transport.
Roundoff, integration error, Taylor truncation, and imperfect torus or
Floquet data introduce small components outside the selected bundle.  For
the full stable manifold, with $\Et=E^s$ and
$\Eh=E^c\oplus E^u$, the desired signal is contracted in forward time while
transverse center/unstable components are contracted less or amplified.
Backward time reverses this relative behavior.  After $j$ linearized
returns the two factors in
\eqref{eq.stable-relative-amplification} are simply raised to the power
$j$, which explains why the distinction becomes important for strongly
hyperbolic problems \cite{GimenoJNO21}.

The conclusion is not universal for every partial spectral selection.  If,
for example, a slow stable bundle is selected while a more strongly
contracting stable bundle is left in the complement, reversing time also
reverses their relative ordering.  The orientation should then be chosen
from the actual spectral ratios together with the non-resonance conditions
of Section~\ref{sec.solve.cohom}; see also
\cite{HaroLl05c,HaroLl06b}.  Thus $\sigma=-1$ is naturally adapted to the
full stable bundle, but should not be interpreted as an unconditional
preconditioner for arbitrary spectral subspaces.

\smallskip

Time reversal has a different effect on the cohomological equations: it
preserves their exact resonance set.

\begin{pro}[Invariance of the exact resonance set]
\label{pro:stable-backward-resonance}
Consider a scalar Fourier--Taylor coefficient of the diagonal
cohomological equation in Section~\ref{sec.solve.cohom}, and let
$\beta_\alpha$ denote the product of selected Floquet multipliers associated
with the multiindex $\alpha$, as in \eqref{non-ded}. Let
$d^+_{k,\alpha}$ be the scalar coefficient multiplying the Fourier unknown
in \eqref{eq.cohom-graph-fourier}, and let $d^-_{k,\alpha}$ be the
corresponding coefficient for $\sigma=-1$.  Then
\begin{equation}
\label{eq.stable-denominator-relation}
    d^-_{k,\alpha}
    =-
    \frac{d^+_{k,\alpha}}
    {e^{2\pi\I\innerp{k}{\rho}}\,\mu\,\beta_\alpha}.
\end{equation}
Consequently,
\begin{equation}
\label{eq.stable-resonance-equivalence}
    d^-_{k,\alpha}=0
    \quad\Longleftrightarrow\quad
    d^+_{k,\alpha}=0,
\end{equation}
and, for each fixed multiindex $\alpha$, the forward and backward
small-divisor estimates differ only by the nonzero factor
$(|\mu|\,|\beta_\alpha|)^{-1}$, which is independent of the Fourier mode
$k$.
\end{pro}

\begin{proof}
For $\sigma=-1$, the rotation changes from $\rho$ to $-\rho$ and
\eqref{eq.GT.directed-blocks} replaces each selected and complementary
multiplier by its inverse. Thus $\beta_\alpha$ is replaced by
$\beta_\alpha^{-1}$, and the coefficient corresponding to
\eqref{eq.cohom-graph-fourier} is
\[
    d^-_{k,\alpha}
    =\mu^{-1}
     -e^{-2\pi\I\innerp{k}{\rho}}\beta_\alpha^{-1}.
\]
Putting these terms over a common denominator gives
\eqref{eq.stable-denominator-relation}.  Local invertibility excludes zero
Floquet multipliers, so the prefactor never vanishes and
\eqref{eq.stable-resonance-equivalence} follows.  Taking absolute values in
\eqref{eq.stable-denominator-relation} gives the stated rescaling of the
small divisors.
\end{proof}

Proposition~\ref{pro:stable-backward-resonance} shows that $\sigma=-1$ does
not avoid the non-resonance assumptions represented by \eqref{non-ded}, nor
does it guarantee larger Fourier denominators.  The numerical gain described
above therefore occurs primarily in the evaluation and transport of the
nonlinear map and its jets, rather than through a different mathematical
solvability condition.  This dynamical reversal should also not be confused
with the two Neumann-series representations
\eqref{eq.forward-series-solution} and \eqref{eq.backward-series-solution}:
those are alternative representations of the inverse of a \emph{fixed}
cohomological operator, whereas $\sigma=-1$ changes the underlying
dynamical map before that operator is assembled.

\smallskip

Finally, the backward formulation requires local invertibility near
$\torinv$.  If $\poincmap$ is the time-$T$ map of a smooth ODE flow,
$\poincmap_{-1}$ is evaluated by integrating the same flow for time $-T$
(and applying the same jet-transport machinery), not by forming an inverse
matrix.  For a genuinely noninvertible map or for a semi-flow without
backward evolution, this formulation is unavailable.  Even when it exists,
backward integration can have a different computational cost or stiffness
profile.  Thus the two orientations are mathematically equivalent
descriptions of the same local stable manifold, but they are not
numerically interchangeable.

\end{document}